\documentclass[12pt]{amsart}
\usepackage{amsbsy}
\usepackage{graphicx,epsfig,subfigure,psfrag}
\begin{document}

\newtheorem{theorem}{Theorem}
\newtheorem{proposition}{Proposition}
\newtheorem{lemma}{Lemma}
\newtheorem{corollary}{Corollary}
\newtheorem{definition}{Definition}
\newtheorem{remark}{Remark}
\newcommand{\be}{\begin{equation}}
\newcommand{\ee}{\end{equation}}
\newcommand{\tex}{\textstyle}
\numberwithin{equation}{section} \numberwithin{theorem}{section}
\numberwithin{proposition}{section} \numberwithin{lemma}{section}
\numberwithin{corollary}{section}
\numberwithin{definition}{section} \numberwithin{remark}{section}
\newcommand{\ren}{\mathbb{R}^N}
\newcommand{\re}{\mathbb{R}}
\newcommand{\n}{\nabla}
\newcommand{\iy}{\infty}
\newcommand{\pa}{\partial}
\newcommand{\fp}{\noindent}
\newcommand{\ms}{\medskip\vskip-.1cm}
\newcommand{\mpb}{\medskip}
\newcommand{\BB}{{\bf B}}
\newcommand{\AAA}{{\bf A}}
\newcommand{\Am}{{\bf A}_{2m}}
\newcommand{\ef}{\eqref}
\newcommand{\eee}{{\mathrm e}}
\newcommand{\ii}{{\mathrm i}}
\renewcommand{\a}{\alpha}
\renewcommand{\b}{\beta}
\newcommand{\g}{\gamma}
\newcommand{\G}{\Gamma}
\renewcommand{\d}{\delta}
\newcommand{\D}{\Delta}
\newcommand{\e}{\varepsilon}
\newcommand{\var}{\varphi}
\renewcommand{\l}{\lambda}
\renewcommand{\o}{\omega}
\renewcommand{\O}{\Omega}
\newcommand{\s}{\sigma}
\renewcommand{\t}{\tau}
\renewcommand{\th}{\theta}
\newcommand{\z}{\zeta}
\newcommand{\wx}{\widetilde x}
\newcommand{\wt}{\widetilde t}
\newcommand{\noi}{\noindent}
\newcommand{\uu}{{\bf u}}
\newcommand{\UU}{{\bf U}}
\newcommand{\VV}{{\bf V}}
\newcommand{\ww}{{\bf w}}
\newcommand{\vv}{{\bf v}}
\newcommand{\WW}{{\bf W}}
\newcommand{\hh}{{\bf h}}
\newcommand{\di}{{\rm div}\,}
\newcommand{\inA}{\quad \mbox{in} \quad \ren \times \re_+}
\newcommand{\inB}{\quad \mbox{in} \quad}
\newcommand{\inC}{\quad \mbox{in} \quad \re \times \re_+}
\newcommand{\inD}{\quad \mbox{in} \quad \re}
\newcommand{\forA}{\quad \mbox{for} \quad}
\newcommand{\whereA}{,\quad \mbox{where} \quad}
\newcommand{\asA}{\quad \mbox{as} \quad}
\newcommand{\andA}{\quad \mbox{and} \quad}
\newcommand{\withA}{,\quad \mbox{with} \quad}
\newcommand{\orA}{,\quad \mbox{or} \quad}
\newcommand{\ssk}{\smallskip}
\newcommand{\LongA}{\quad \Longrightarrow \quad}
\def\com#1{\fbox{\parbox{6in}{\texttt{#1}}}}
\def\N{{\mathbb N}}
\def\A{{\cal A}}
\def\WW{{\cal W}}
\newcommand{\de}{\,d}
\newcommand{\eps}{\varepsilon}
\newcommand{\spt}{{\mbox spt}}
\newcommand{\ind}{{\mbox ind}}
\newcommand{\supp}{{\mbox supp}}
\newcommand{\dip}{\displaystyle}
\newcommand{\prt}{\partial}
\renewcommand{\theequation}{\thesection.\arabic{equation}}
\renewcommand{\baselinestretch}{1.2}

\title
{\bf Towards the KPP--problem  and
  ${\bf { log\, t}}$--front \\ shift for higher-order nonlinear PDEs
I.\\ Bi-harmonic and  other parabolic equations}

\author{
V.A.~Galaktionov}

\address{Department of Mathematical Sciences, University of Bath,
 Bath BA2 7AY, UK}
\email{masvg@bath.ac.uk}



  \keywords{KPP--problem, travelling wave, stability,  higher-order semilinear parabolic
  equations,  $\log t$-front shifting}
 \subjclass{35K55, 35K40, 35K65}
 \date{\today}




\begin{abstract}

The seminal paper by Kolmogorov, Petrovskii, and Piskunov of 1937
\cite{KPP} on the travelling wave propagation in the
reaction-diffusion equation
 \be
 \label{001}
 u_t=u_{xx}+u(1-u) \inB \re \times \re_+, \quad u_0(x)=H(-x) \equiv\{1\,\,
 \mbox{for}\,\, x<0; \,\,\, 0 \,\, \mbox{for} \,\, x \ge 0\}
 \ee
 (here $H(\cdot)$ is the Heaviside function),
 opened a new era in general theory of nonlinear PDEs and
  various applications. This paper became an {\em encyclopedia} of deep mathematical techniques
 and tools for nonlinear parabolic equations, which, in last seventy years, were further
 developed in hundreds of papers and in dozens of monographs.

  The KPP
 paper established the fundamental fact that, in \ef{001}, there
 occurs a travelling wave $f(x-\l_0 t)$, with the minimal speed
 $\l_0=2$, and, in the moving frame with the front shift $x_f(t)$
 ($u(x_f(t),t) \equiv \frac 12$), there is uniform convergence
 $ u(x_f(t)+y,t) \to f(y)$ as $ t \to + \iy$, where $ x_f(t) = 2t(1+o(1))$.
In 1983, by a probabilistic approach, Bramson \cite{Br} proved
that there exists an {\em unbounded} $\log t$-shift of the wave
front in the PDE problem \ef{001} and
 $
 x_f(t) = 2t - \frac 32 \, \log t(1+o(1))$ as $ t \to +\iy$.

Our  goal
is  to reveal some aspects of  KPP-type problems for higher-order
semilinear parabolic
PDEs,
including the {\em bi-harmonic} equation  and the {\em
 tri--harmonic}
one
 $$
 u_t= -u_{xxxx} +u(1-u) \andA
 u_t=u_{xxxxxx}+u(1-u),
 $$
  and other poly-harmonic PDEs up to tenth order.
 Two main questions to study are:

 (i) existence of travelling waves via any analytical/numerical methods, and

 (ii) a formal derivation of the
 $\log t$-shifting of moving fronts.

\end{abstract}

\maketitle



\setcounter{equation}{0}
\section{Introduction: the classic KPP--problem, known
remarkable results, and extended PDE models}
 \label{Sect1}
  \setcounter{equation}{0}











\subsection{The classic KPP--problem of 1937: convergence to TWs}

In the seminal  Kolmogorov--Petrovskii--Piskunov (KPP) problem of
1937 \cite{KPP}
 \be
 \label{1.1}
  u_t = u_{xx} +u(1-u) \inB \re \times \re_+,
\quad u(x,0)=u_0(x) \,\,\, \mbox{in\,\, $\re$},
 \ee
    with
the step (Heaviside) initial function
 \be
 \label{1.H}
u_0(x) =H(-x) \equiv \left\{ \begin{matrix}1, \quad x<0;\\ 0,
\quad x\ge 0,
\end{matrix}
\right.
 \ee
the solution was proved to converge to the so-called {\em minimal
travelling wave} (TW)  corresponding to
 \be
 \label{minSp1}
 \mbox{the minimal speed of TW propagation}: \quad \l_0=2.
  \ee
   Namely, looking for a TW profile $f(y)$ for
an arbitrary speed $\l>0$ yields
 \be
 \label{TW1}
 \left\{
  \begin{matrix}
 u_*(x,t)=f(y), \,\, y=x-\l t, \,\,\,\,\mbox{where $f$ solves the ODE problem}\ssk\ssk\ssk\\
  -\l f'=f''+f(1-f), \,\, y \in \re; \quad
 f(-\iy)=1, \,\, f(+\iy)=0. \qquad
 \end{matrix}
 \right.
  \ee
This 2nd-order ODE, on the phase-plane $\{f,f'\}$, by setting
$f'=P(f)$, reduces its order:
 \be
 \label{TW11}
  \tex{
  \frac {{\mathrm d}P}{{\mathrm d}f}=-\l- \frac {f(1-f)}P \, ,
}
 \ee
 and it was shown that there exists the {\em minimal} speed
 $\l_0=2$ and the corresponding
 {\em minimal} TW profile $f(y)$.
Using the natural normalization
 \be
 \label {norm1}
 \tex{
 f(0)= \frac 12,
 }
 \ee
 this minimal TW profile is defined uniquely.
  In addition,
  \be
  \label{norm2}
  f'(y) <0 \inB \re.
   \ee
The characteristic equation for the linearized operator in
\ef{TW1}, $\l=2$, has a multiple zero:
 \be
 \label{Norm21}
g'' + 2 g'+g=0 \andA g=\eee^{\mu y} \LongA (\mu+1)^2=0 \LongA
\mu_{1,2}=-1,
 \ee
 that yields the following asymptotic behaviour of $f(y)$:
 \be
  \label{norm23}
   f(y)=C_0 y \, {\mathrm e}^{-y}(1+o(1)) \asA y
  \to + \iy \whereA \mbox{$C_0>0$ is a constant.}
  \ee

 Concerning the relation between the
ODE TW problem \ef{TW11} for $\l_0=2$ and the PDE Cauchy one
\ef{1.1}, \ef{1.H}, the novel remarkable analysis in
\cite{KPP}\footnote{In particular, the original KPP-proof
crucially used a version of Sturmian zero set theorem of 1836,
which was re-discovered by the authors independently. Moreover, in
1937 (101 years after Sturm!), this was the first ever use of this
Sturm's First Theorem on non-increase of the number of
intersections of {\em all} the TWs for {\em any} $\l \ge \l_0$
with the solution $u(x,t)$ of \ef{1.1}, \ef{1.H} of 1D parabolic
PDEs since Sturm's original result in 1836. See further historic
comments on those amazing facts in \cite[p.~23]{GalGeom} and
around.} of convergence as $t \to +\iy$ of the solution of the
Cauchy problem \ef{1.1}, \ef{1.H} to the minimal TW \ef{TW1} was
performed in the TW moving frame. This was very essential, and not
in view of the obvious $x$-translational invariance of the
equation \ef{1.1}; see below. Eventually, using PDE methods
 and naturally defining the front location via
  \be
 \label{TW2}
  \tex{
 u(x_f(t),t) = \frac 12 \quad \mbox{for all\,\, $t \ge 0$},
 }
  \ee
 the
KPP-authors proved that
  the
TW front asymptotically moves like
 \be
 \label{1.2}
  x_f(t) = 2t -
g(t) \whereA g'(t) \to 0 \quad \mbox{as $t\to + \infty$}.
 \ee

Then the convergence result of \cite{KPP} takes the form:
 \be
 \label{TW3}
 u(x_f(t)+y,t) \to f(y) \asA t \to + \iy \quad \mbox{uniformly in
 $y \in \re$}.
 \ee

These classic results were proved for more general source terms
$Q(u)$ than the quadratic one in \ef{1.1}, satisfying some natural
restrictions. However, here and later on, we keep using this
simplest reaction term
\be
\label{FF1}
 Q(u)=u(1-u),
  \ee
bearing in mind and viewing KPP results as establishing, in
particular, a kind of a {\em structural stability} of the problem.
Concerning the ODE \ef{TW11}, with $\l_0=2$, the results of
\cite{KPP} imply the actual structural stability  of the dynamical
system \ef{TW11} on the phase-plane, at least, in a neighbourhood
of the heteroclinic connection $(0,0) \to (1,0)$, that generates
the TW profile $f(y)$. In particular, this means stability
relative any small perturbation in $C^1$ of the nonlinearity
\ef{FF1}. A similar much more delicate structural stability result
is obtained in \cite{KPP} for the PDE \ef{1.1} relative the source
$Q$. We plan to trace out similar (but indeed weaker)
 structural stability
properties for higher-order ODEs involved, so will keep dealing
with the reaction given by \ef{FF1}.


\subsection{Bramson's $\log t$-front drift}

We next point out   the next principal question arising within the
KPP ideology. Namely,
 it is the question on the actual behaviour of the TW shift $g(t)$
in \ef{1.2} for $t \gg 1$, which was not addressed in \cite{KPP}
and, moreover, it  was not  even mentioned therein whether $g(t)$
remains bounded or gets unbounded as $t \to +\iy$ (possibly, it
might not that important for the KPP-authors). Later on, it turned
out that this is a principle, difficult, and rather general
question in such PDE problems. Indeed, this is about a ``centre
subspace (manifold) drift" of general solutions of the KPP PDE
\ef{1.1} along  a one-parameter family of exact ODE TWs
$\{f(y+a),\,\, a \in \re\}$.

\ssk

This open problem was solved in 1983 (i.e., 46 years later), when
Bramson \cite{Br} (see also \cite{Ga}), by pure
probabilistic techniques based on the Feynman--Kac integral
formula together with sample path estimates for Brownian motion,
proved that, within the PDE setting \ef{1.1}, \ef{1.H},  there is
an {\em unbounded} $\log t$-shift of the moving TW
front\footnote{The exponent $k= \frac 32$ in \ef{1.3}
  seems deserve to be compared with the ``magic" $\frac 23$
 in Kolmogorov--Obukhov theory of local structure of
turbulence (discovered in 1941, i.e., just 4 years after the
KPP-paper). Namely, the famous dimensional
  {\em Kolmogorov--Obukhov power ``K-$41$"
law} (1941) \cite{Kolm41, Obu41}  for the energy spectrum of
turbulent fluctuations for wave numbers $k$ from the
{\em inertial range},
 \be
  \label{Kom}
  \tex{
 E(k)= C \e^{\bf \frac 23} \,\, k^{-\frac 53} \whereA \e= \frac 1{\rm Re}
 \, \langle |D \uu|^2 \rangle,
 }
 \ee
 describing the rate of dissipation of kinetic energy for high Reynolds numbers
  ${\rm Re}$ ($\langle \cdot
 \rangle$ is an  invariant measure of calculating expected
 values).}
 \be
 \label{1.3}
  g(t) = k \log t(1+o(1)) \asA t \to +\iy, \quad
\mbox{with\,\, ${\bf k=  {\frac 32}}$},
  \ee

  Thus, \ef{1.3} implies eventual, as $t \to +\iy$, {\em infinite} retarding of the
  solution $u(x,t)$ from the corresponding minimal TW (uniquely fixed by \ef{norm1}), thought the
  convergence \ef{TW3} takes place in the TW frame.



 \noi {\bf Remark: on PDE approaches to the KPP--2 problem.}
Several important results have been already proved in many
KPP-like studies during last seventy years by using PDE methods.
 For instance, sharp front propagation results have been obtained
 for various classes of {\em positive} initial data $u_0(x)$
 with different decay rates as $x \to + \iy$; see \cite{Cou11, Ebert2000,
 Hamel09, Yan06} as a source of further references and results. Moreover,
 it seems, it was not studied whether special types of pattern convergence to
TWs exist for initial data $u_0(x)$ of {\rm changing sign}. For
instance, if $u_0(x)$ has a fixed number $l \ge 1$ of zeros (sign
changing) that persist for all times $t \in \re_+$. In particular,
it is unclear still how (bearing in mind  Sturm's First Theorem of
1836 on nonincrease of the number of zeros; see \cite{GalGeom} for
the history and various applications) could this affect the rate
of convergence to the minimal TW and the corresponding $\log
t$-shift?

It seems, in this framework, there is no still a fully developed
PDE (non-probabilistic using ideas of Brownian motion) approach to
establishing such a precise expansion as in \ef{1.2}, \ef{1.3}.
\subsection{The main goal of the present paper: extensions of the KPP--problem
to higher-order semilinear parabolic PDEs}

The main goal of the present paper is to show that the
KPP--ideology is very wide and can be extended (along very similar
lines) to a variety of other more complicated higher-order
semilinear PDEs with the same source-type term. It is worth
mentioning  now, that, since for such PDEs, no Maximum Principle,
comparison, Sturm's, and other related properties of
order-preserving semigroups apply, one cannot expect in principle
so complete, detailed,   and beautiful results as in the classic
KPP--Bramson study of 1937--1983.

 Of course, there are already
many strong results, to be referred to later on, concerning a
TW-like analysis in Cahn--Hilliard and related higher-order
parabolic equations. However, our overall goal here and in
forthcoming papers is to
initiate a more general study, and, using {\em any mathematical
means} (including various analytic, formal, and, inevitably,
numerical methods), to show that such a general viewing of the KPP
ideas makes sense, and that many higher-order semi- and
quasilinear PDEs of different types inherit some  deep (but not
all) key features of this classic KPP analysis.

\subsection{On some known front propagation results for bi-harmonic and higher-order diffusion}
Note that, nowadays, there is a vast enough literature devoted to
front propagation features for fourth-order semilinear parabolic
equations such as the {\em Swift--Hohenberg} one
 \be
 \label{SH1}
 u_t=-(1+D_x)^2u + \a u-u^3.
 \ee
 One of the first and detailed such study was performed in 1990 in the monograph
 \cite{Collet90}. Further and more recent papers can be then traced
 out using the {\tt MathSciNet}.  From applications, such
KPP-type problems mean that higher-order diffusion is taking into
account; see references on various fourth- and  $2m$th-order
semilinear and quasilinear parabolic  PDEs in \cite{Bert01,
Collet90, Gl4, EGW1, Gal2m, GalCr}.
 In particular, as further generalizations, in \cite{Bert01} (see also references therein and in the current paper), TW profiles were studied for
   a class of
 quasilinear {\em thin film equations}
  \be
  \label{TFE1}
  u_t+(f(u))_x= (b(u)u_x)_x-(c(u)u_{xxx})_x \quad (c(u) \ge 0).
  \ee
Questions of  general {\em instabilities} of TWs in such {\em
semilinear} Cahn--Hilliard and other related parabolic PDEs (these
results are of a special attention in the present study) were also
already  treated in a number of papers; see \cite{Gao2004,
Strauss2004, Li2012} and references there in. We will quote and
use later on some other known results.

 However, it seems
that the study
 that is more and directly oriented to the KPP--ideology, including possible types of unbounded front shift, was not
 performed before.

\subsection{Layout of the paper}

Thus,  we will discuss some aspects of  KPP-type problems for
higher-order semilinear {\em parabolic} partial differential
equations (PDEs), with the same Heaviside initial data.

In Section \ref{S2}, as a natural and simplest KPP-like extension,
we consider the first model:
 the {\em semilinear bi-harmonic equation}, i.e.,
  a fourth-order semilinear heat or a reaction-diffusion
  equation (SHE--4 or RDE--4)
 \be
 \label{E4}
 u_t= -u_{xxxx} +u(1-u) \inB \re \times \re_+, \quad u(x,0)=H(-x)
 \inB \re.
  \ee
The corresponding TW with the speed of propagation $\l$ is then
governed by the following fourth-order ODE:
 \be
 \label{E5}
u_*(x,t)=f(y), \quad  y=x - \l t \LongA -\l f'=-f^{(4)}+f(1-f),
 \ee
 with the same singular boundary conditions at infinity:
 \be
 \label{BC1}
 f(y) \to 0  \andA f(y) \to 1 \asA y \to \pm \iy \quad
 \mbox{``maximally" exponentially fast.}
  \ee
This ``maximal" decay of $f(y)$ at infinity somehow includes some
kind of the remnants of a  ``minimality" of the possible TW
profiles, though, for such higher-order equations, any direct
specification of such a property is difficult to express
rigorously (and literally).

For any $\l \ne 0$, the problem \ef{E5}, \ef{BC1} is  of the
elliptic type, but {\em it is not a  variational one}. Therefore,
we cannot used  advanced methods for higher-order ODEs with
potential operators associated with homotopy-hodograph and other
approaches (see key examples in \cite{KKVV00, PelTroy, VV02})
and/or Lusternik--Schnirel'man and fibering theory (see
\cite{GMPSobI, GMPSobII} and references therein). Hence, the ODE
\ef{E5}, though looking rather simple, and, at least, simpler than
most of related fourth-order ODEs already studied in detail,
represents a serious challenge and cannot be tackled directly by
known tools of modern nonlinear analysis and operator theory.

The KPP--4 ODE problem \ef{E5}, \ef{BC1} on existence of a
heteroclinic connection $(0,0,0,0) \mapsto (1,0,0,0)$ in a
four-dimensional phase space, is more difficult that its
second-order KPP--2 counterpart \ef{TW1}, so we cannot easily
obtain which ``minimal" (if any) speed $\l_0$ occurs for the
Heaviside data $H(-x)$. Moreover, we will show that, in the usual
sense, such a ``minimal speed" to be understood in the usual
sense, is in fact nonexistent for many higher-order equations.
Therefore, we first study TW profiles for rather arbitrary $\l \in
\re$, and will denote them by
 $$
 f=f(y;\l).
 $$
Here, for the ODE problem \ef{E5}, \ef{BC1}, it is easy to
describe directly the exponential oscillatory bundles of orbits as
$y \to \pm \iy$. Such an analysis, even more general and complete,
is well established; see references in Section \ref{S2}. Next, by
a shooting technique, we can justify existence of such TWs for
sufficiently small $\l>0$. It turns out that, for such a KPP--4
problem (as well as for many others), the question of existence of
a kind of a ``minimal" speed becomes irrelevant. Moreover, we do
not think that, for such PDEs and ODEs without any traces of the
Maximum and Comparison Principles, the notions of ``minimal
solution/speed" can be properly and so easily defined. Actually,
often, these essentially do not make sense, so that this part of
the classic KPP analysis seems do not admit a direct extension. On
the contrary, in several cases,  numerically, we have observed a
``maximal" speed $\l_{\rm max}>0$ such that
 \be
 \label{SpeedMax}
 \mbox{TW profiles $f(y;\l)$ exist for all $0< \l < \l_{\rm max}$, and nonexistent
 for $\l>\l_{\rm max}$},
 \ee
 and rather sharply estimated $\l_{\rm max}=\l_{\rm max}(m)>0$ for $m=2,3,4$,
 i.e., up to the eighth-order parabolic equation (see \ef{m1}
 below). However, we must admit that such a numerical
 identification of the $\l_{\rm max}$ is difficult and is not
 always reliable, since, for higher-order ODEs, the influence
of artificial boundary conditions at the end points of a
sufficiently large interval of integration cannot be completely
eliminated, especially if the solutions are highly oscillatory
therein. Sometimes, we also cannot guarantee nonexistence of TW
profiles $f(y)$ for $\l$ much larger than $\l_{\rm max}$, when the
{\tt bvp4c} solver of the {\tt MatLab} often produces a very slow
convergence.

\ssk

 In Section \ref{S3}, similarly, we study the {\em
semilinear tri-harmonic equation} (SHE-6):
 \be
 \label{E6}
 u_t= u_{xxxxxx} +u(1-u) \inB \re \times \re_+, \quad u(x,0)=H(-x)
 \inB \re.
  \ee
The corresponding TW with the speed of propagation $\l$ is then
governed by the following sixth-order ODE:
 \be
 \label{E7}
u_*(x,t)=f(y), \quad  y=x - \l t \LongA -\l f'=f^{(6)}+f(1-f),
\,\,\, \mbox{with (\ref{BC1})}.
 \ee

 \ssk

In Section \ref{S810}, we present some numerical evidence on
existence of various TWs and there properties for $2m$th-order
semilinear heat (parabolic) equations  (SHE-$2m$), such as
\be
\label{m1}
 u_t= (-1)^{m+1} D_x^{2m} u + u(1-u), \quad \mbox{with the ODEs}
 \ee
 \be
 \label{m2}
u_*(x,t)=f(y), \quad  y=x - \l t \LongA -\l
f'=(-1)^{m+1}f^{(2m)}+f(1-f)
 \ee
 (plus
(\ref{BC1})),
 and treat higher-order cases for $m=4$ and $m=5$, i.e., eight- and tenth-order
 heat equations. A related approach based on construction of a
 majorizing operator for $2m$th-order parabolic equations for  $m \ge 2$ as in
 \ef{m1} is discussed in Appendix B.

\ssk

 Finally, in Section \ref{SDiscr}, we explain how $\log t$-front
 shift occurs in those higher-order KPP-problems. Section
 \ref{S6omega} is devoted to the omega-limits of PDEs involved.



\subsection{An extension:
 towards quasilinear parabolic KPP-problems}
 The present paper
deals with higher-order {\em semilinear parabolic} equations of
reaction-diffusion type, while, in the second part of this
research \cite{GKPPII}, we extend some of the above results to the
{\em quasilinear} KPP--$4n$ problem for
\be
 \label{E4n}
 u_t= -(|u|^n u)_{xxxx} +u(1-u) \inB \re \times \re_+,
  \ee
  where $n>0$ is a parameter, as well as to some other parabolic
  equations.

\subsection{An extension:
to other PDE types and settings}

 In \cite{GKPPIII}, we will deal with
higher-order  {\em dispersion} and {\em hyperbolic} equations such
as
 \be
 \label{m3}
 u_{t}= -D_x^{11} u +u(1-u) \andA  u_{tt}= -u_{xxxx} +u(1-u),
 \quad \mbox{etc.}
  \ee
As for more ``exotic" PDE models, as a formal but quite
illustrative examples, we consider, in \cite{GKPPIII},
higher-order dispersion equations and end up with the following
one:
 \be
 \label{m4}
 D_t^9 u= -D_x^{11}u + u(1-u),
  \ee
  with eleventh-order ODE for the TW profiles
  \be
  \label{m41}
  -\l^9 f^{(9)}= -f^{(11)} + f(1-f) \inB \re \quad (\mbox{plus\,\,
  (\ref{BC1})}).
  \ee

 Overall, using the KPP-setting, we refer to such problems
 as to KPP--$(k,l)$, where $k$ stands for the order of the
 differential operator in $x$ and $l$ for the order of the
 derivative in $t$. Therefore, equation \ef{m4} represents the
 KPP--(11,9) problem.


\subsection{Summary of Main General Questions}

 Thus, for several  KPP--$(k,l)$ problems, with $k \ge 3$ and $l \ge 1$, the  main
 questions
  to study, here  and in \cite{GKPPII, GKPPIII}, are mainly:

\ssk

 {\bf (I)} {\sc The problem of TW existence:} existence of travelling waves
 via any analytical/numerical methods and estimating the
 $\l$-interval of existence:
  \be
  \label{LL1}
  \lambda= \{\l \in \re: \quad \mbox{there exists a TW profile
  $f(y;\l)$}\}.
  \ee
  Firstly, we prove the {\em positivity}, which is easy for all the parabolic models such as \ef{m1}:
   \be
   \label{l>0}
   \l \in \Lambda \LongA \l>0.
    \ee

\ssk

{\bf (II)} {\sc The problem of a ``maximal" speed:} whether
 \be
 \label{m42}
  \Lambda \quad \mbox{is bounded from above}?
   \ee
 In particular, two main questions arise:
  \be
  \label{m44}
\mbox{is there a  maximal speed} \quad \l_{\rm max} = \sup
\Lambda>0?
 \ee
 Those questions are ``remnants" of the KPP setting. It
 turns out that, unlike in \cite{KPP}, for several types of
 nonlinear PDEs, we have found that
  \be
  \label{max77}
  \Lambda=(0,\l_{\rm max}), \quad \mbox{so that}
  \quad \l_{\rm max}= \sup \Lambda, \,\,\, \l_{\rm max} \not \in
  \Lambda, \,\,\, \mbox{and}\,\,\,\, \l_{\rm min}=0 \in \Lambda \not \in
  \Lambda
  \ee
 (though analytical/numerical proof of nonexistence for {\sc all} $\l>\l_{\rm
 max}$ requires further study, as some other aspects of these
 problems).

\ssk

 {\bf (III)} {\sc The $\log t$-shift problem:}
 using general and natural instabilities of such TWs $f(y;\l_0)$, with a $\l_0 \in \Lambda$,
 to show and derive  the
 $\log t$-shifting of the moving front in the problem \ef{1.1},
 \ef{1.H},  connected with a kind of an ``(affine) centre subspace
 behaviour" for the rescaled equation.

 Thus, we believe that the $\log t$-shifting phenomenon is
quite a generic property of many semilinear or quasilinear
\cite{GKPPII, GKPPIII} KPP-problems of different types.

\ssk

{\bf (IV)} {\sc  The omega-limit set $\o(H)$:}
 the last and difficult question, which is not
being properly understood, is the actual behaviour, as $t \to
+\iy$, of the solution $u(x,t)$ of the PDE KPP-problems with data
$H(-x)$, or $\o(u_0)$, where $u_0$ is a bounded ``step-like"
function. Namely,
 defining, in a
natural sense, the $\o$-limit set $\o(H)$ of the properly shifter
orbit \ef{TW3}, to discuss whether or not
 \be
 \label{m9N}
 \o(H) \subseteq \{f(\cdot;\l), \quad \l \in  \Lambda\} \quad \mbox{and/or} \quad
 \o(H)=\{f(\cdot;\l_0), \quad \l_0 \in \Lambda\}.
  \ee
  The latter means that the rescaled (shifted) PDE orbit converges
to a single TW profile with a special ``minimal" speed $\l_0$
(though, for such higher-order parabolic flows with no Maximum
Principle and other order-preserving features, a slow evolution
within a connected subset of different TWs in the former statement
should be also taken into account).

Overall, the latter in \ef{m9N} means that, denoting by $E_{\rm S}
(\l_0)$ the stable subset (a manifold of stability: all reasonable
bounded ``step-like" data $u_0(x)$, for which there exists, after
proper shifting, uniform convergence to $f(\cdot;\l_0$) of some TW
profile $f(\cdot;\l_0)$ with a $\l_0 \in \Lambda$), we want to
know whether
 \be
 \label{m10N}
 H(-x) \in E_{\rm S}(\l_0)?
 \ee
Numerically, this assumes a full-scale numerical experiments in
general PDE settings, which is not done here\footnote{The author
is not an expert  in that at all and is not fond of such a
research, which could lead to a conflict with his current
mathematical interests.}, but it would be rather desirable to be
performed by more professional and, possibly, more applied
mathematicians in this PDE area.

\section{The basic higher-order KPP--4 problem}
\label{S2}

Thus, consider the KPP--(4,1) (or simply KPP--4, that cannot
confuse in the parabolic case) problem \ef{E4}. We begin with its
ODE counterpart \ef{E5}, \ef{BC1}.

\subsection{$\l$ is always positive}

The first simple result is the same as in the KPP--2.

\begin{proposition}
\label{Pr.Lambda1}
\be
\label{Lam1}
 \mbox{If there exists a solution $f(y) \not \equiv 0$
of $(\ref{E5})$, $(\ref{BC1})$, then $\l>0$.}
 \ee
 \end{proposition}

 \noi{\em Proof.} Multiplying the ODE \ef{E5} by $f'$,
 integrating by parts over $\re$,  and using the boundary conditions \ef{BC1}
 yield
  \be
  \label{Lam2}
   \begin{aligned}
   \tex{
  -\l \int(f')^2}= &-
 \tex{
   \int f^{(4)}f'+ \int f(1-f) f'= \int f'' f'''
  + \big[ \frac{f^2}2- \frac {f^3}3 \big]_{-\iy}^{+\iy}
  }
 \\
  = & \,
  \tex{
  \frac 12
  \,[(f'')^2]_{-\iy}^{+\iy}- \frac 16= -\frac 16<0. \qed
  }
  \end{aligned}
  \ee

\subsection{TW profiles do exist: numerical evidence}

It is convenient to present next numerical results, which directly
show the global structure of such TW profiles to be, at least
partially, justified analytically.

To this end, we used the {\tt bvp4c} solver of the {\tt MatLab}
with sufficient accuracy and both tolerances at least of the order
 \be
 \label{t1}
 \mbox{Tols}=10^{-5}- 10^{-7}.
  \ee
It is important to note that, as the initial data for further
iterations, we always took either the Heaviside function
 \be
 \label{Heav1}
 \mbox{initial data for numerical iteration are often} \quad
 H(-y),
  \ee
i.e.,  as in \ef{1.H}, or its slightly smoother version for a
better convergence. This once more had to help us to converge to a
proper ``minimal" profile (indeed, there are many other TW
profiles), though, of course, this was not guaranteed {\em a
priori}. We keep this rule for all other KPP--$(k,l)$ problems of
in \cite{GKPPII, GKPPIII}, including the hyperbolic ones \ef{m3}
and more higher-order ones such as \ef{m4}, where initial
velocity, acceleration, etc. were taken zero.

We show the profiles on sufficiently small $y$-intervals, such as
$[-100,200]$, though, to get right singular boundary conditions
\ef{BC1}, we fixed much larger intervals $[-600,600]$, up to
$[-100,5000]$, when necessary, and even bigger ones to guarantee a
proper convergence and revealing TW profiles of the Cauchy
problem.

\ssk

Figure \ref{F0} shows TW profiles for sufficiently small $\l \in
[0.2,0.5]$. For smaller $\l=0.1$, see below; the results for the
smallest $\l=0.01$ achieved numerically are shown in Figure
\ref{F0small}, where $f(y)$ gets very oscillatory for $y \gg 1$.


 \begin{figure}
\centering
\includegraphics[scale=0.85]{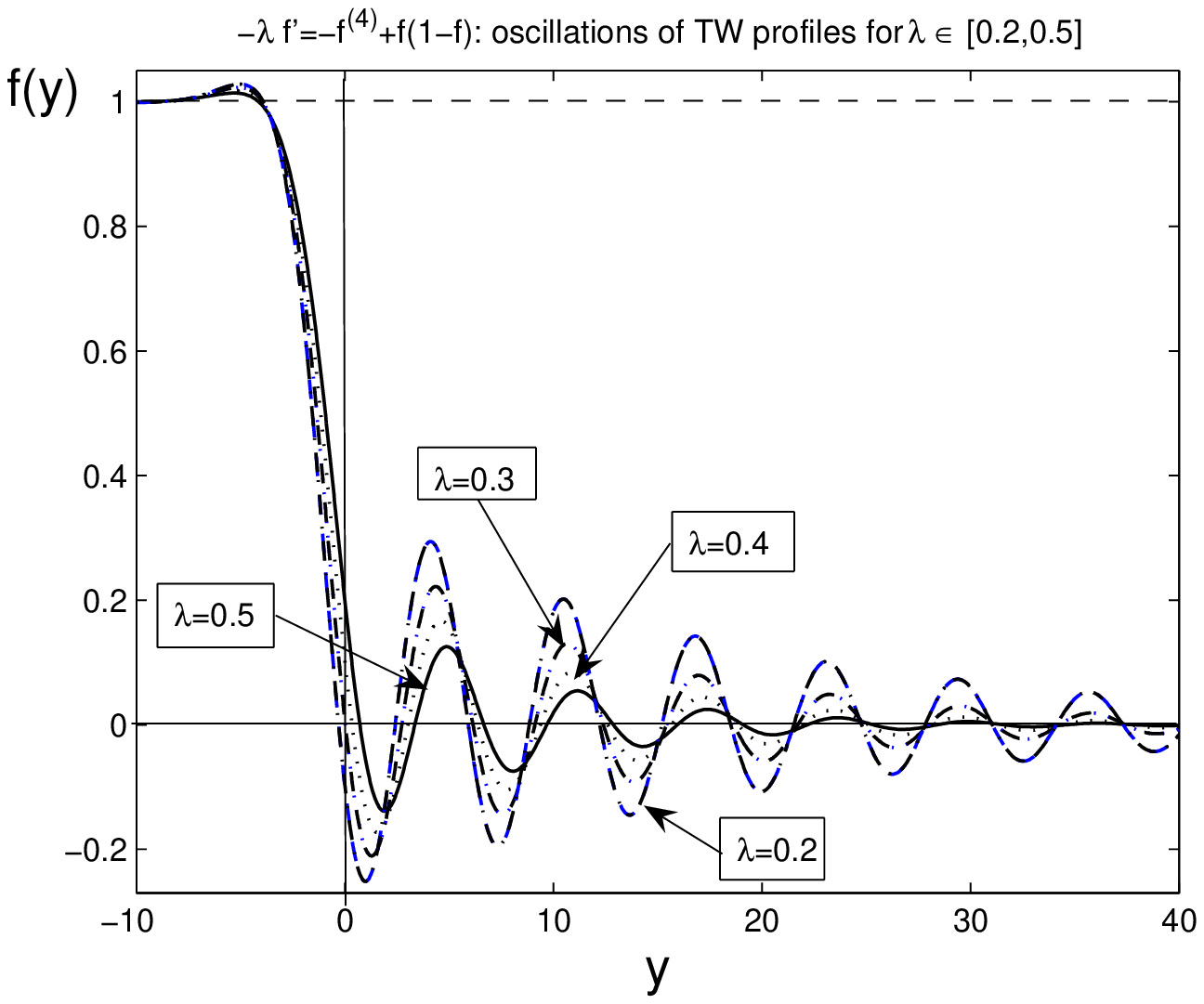}  
\vskip -.3cm
  \caption{The TW profiles $f(y)$ of
(\ref{E5}), \ef{BC1} for
 $\l=0.2$, 0.3, 0.4, and 0.5.}
 \label{F0}
\end{figure}

\begin{figure}
\centering
\includegraphics[scale=0.85]{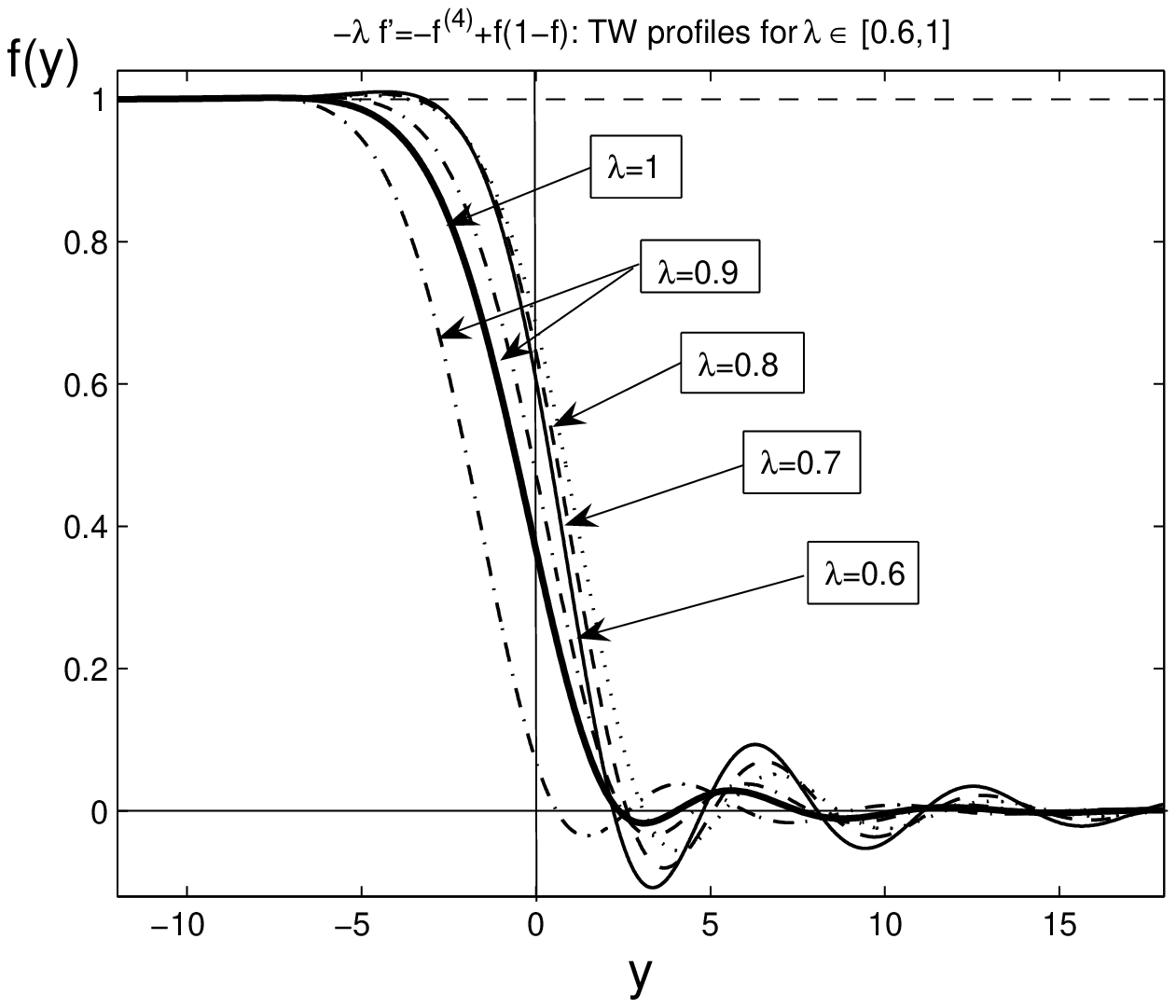}  
\vskip -.3cm
  \caption{The TW profiles $f(y)$ of
(\ref{E5}), \ef{BC1} for
 $\l=0.6$, 0.7, 0.8, 0.9, and 1.}
 \label{F1}
\end{figure}

In Figure \ref{F1}, we show TW profiles  for larger  speeds $\l
\in [0.6,1]$. All of them are clearly oscillatory as $y \to +
\iy$. Concerning the type of convergence  to 1 as $y \to - \iy$,
it  looks like it is monotone, but more close research below will
show its oscillatory character to be justified by asymptotic
expansion methods later on.



\begin{figure}
\centering \subfigure[oscillations on (0,300)]{
\includegraphics[scale=0.52]{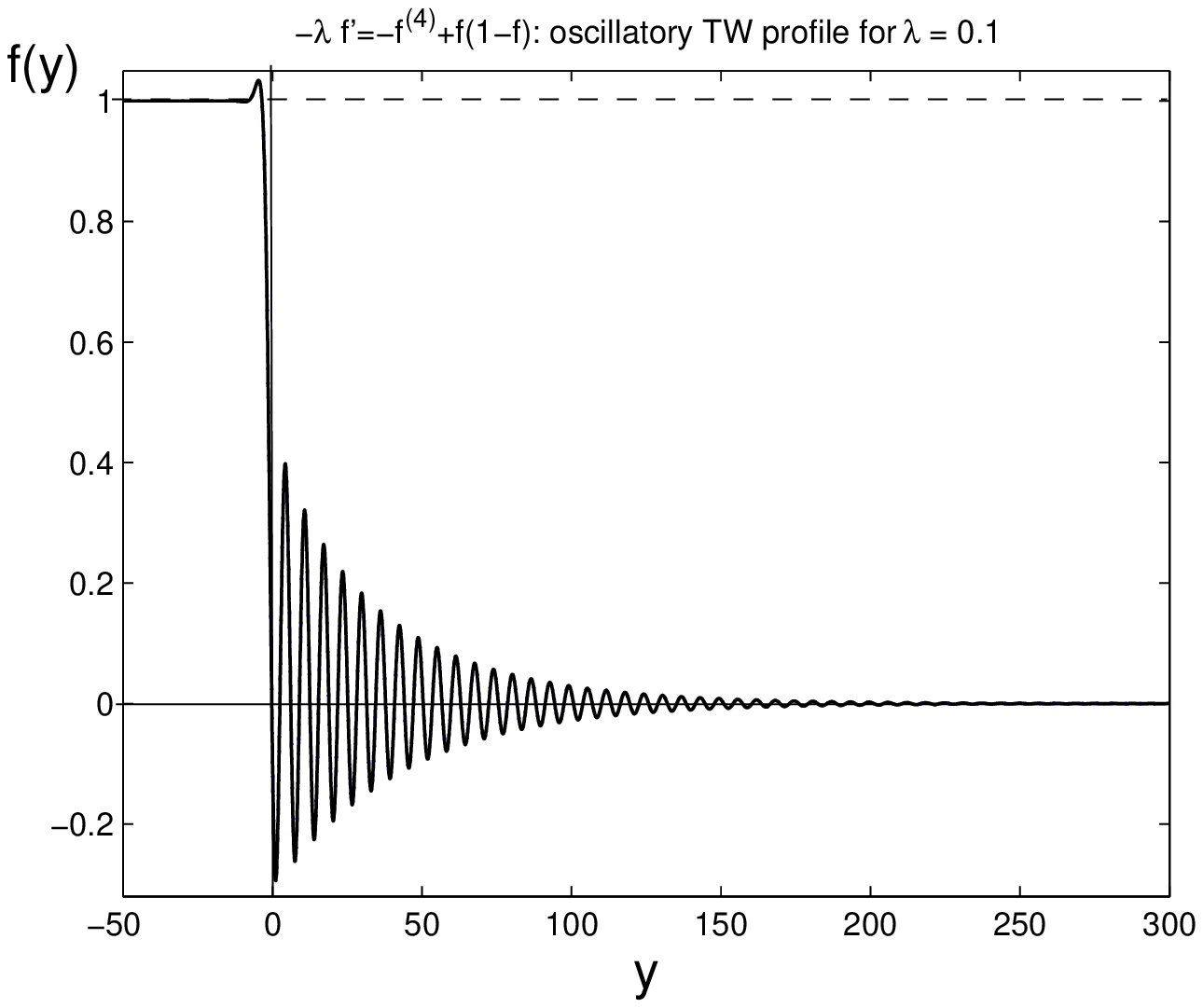}             
} \subfigure[oscillations on (0,600)]{
\includegraphics[scale=0.52]{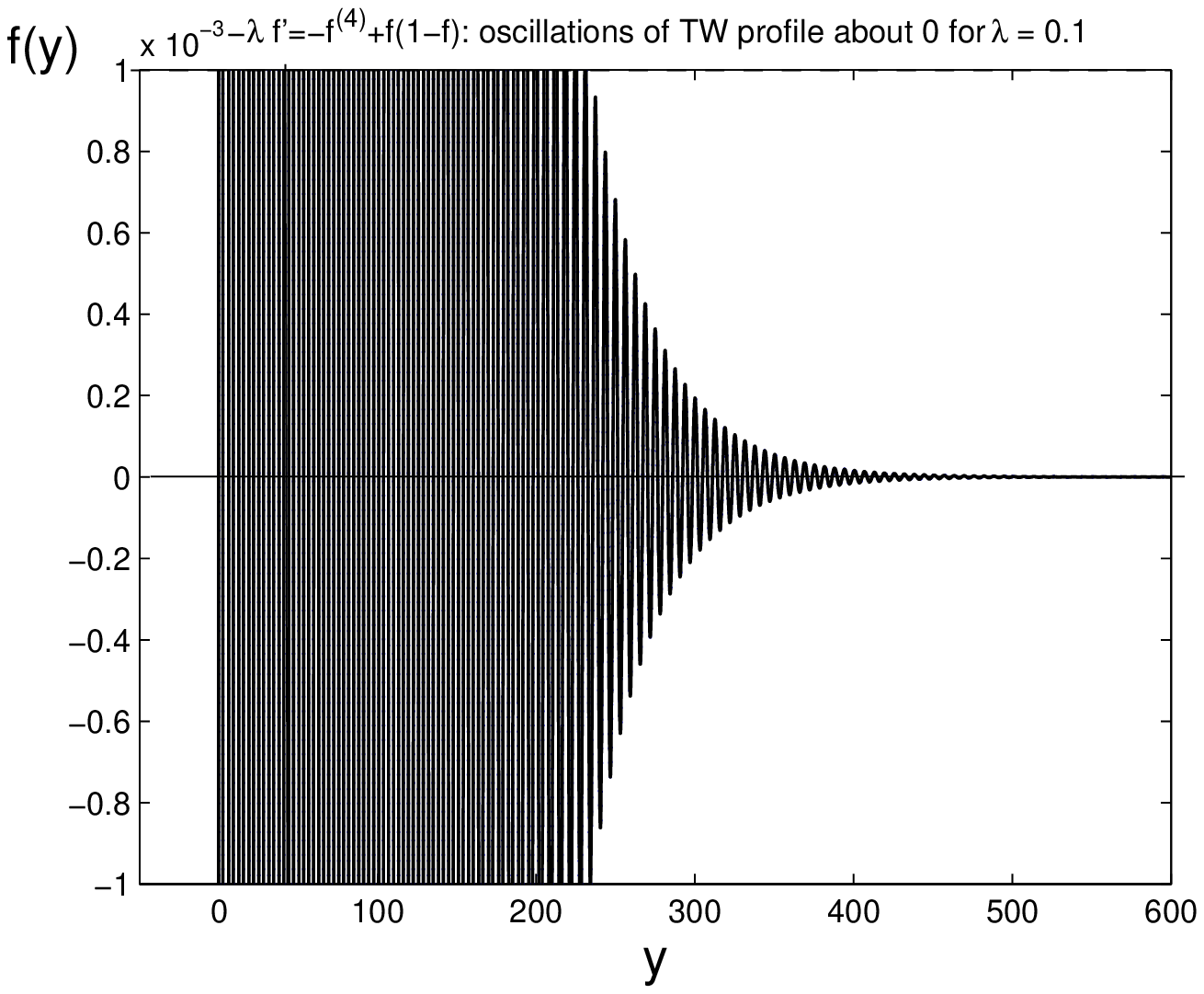}                        
}
 \vskip -.3cm
\caption{\rm\small Oscillations of the TW profile $f(y)$ for $y
\gg 1$; $\l=0.1$.}
 \label{F2}
\end{figure}

Both figures above confirm an important property:
\be
\label{t2} \mbox{the oscillatory behaviour as $y \to + \iy$
decreases as $\l\uparrow$;}
 \ee
see further comments below. Figure \ref{F2} shows the character of
those oscillations for $\l=0.1$ on intervals $[0,300]$ (a) and
$[0,600]$ (b).

\begin{figure}
\centering \subfigure[scale $\sim 10^{-2}$]{
\includegraphics[scale=0.52]{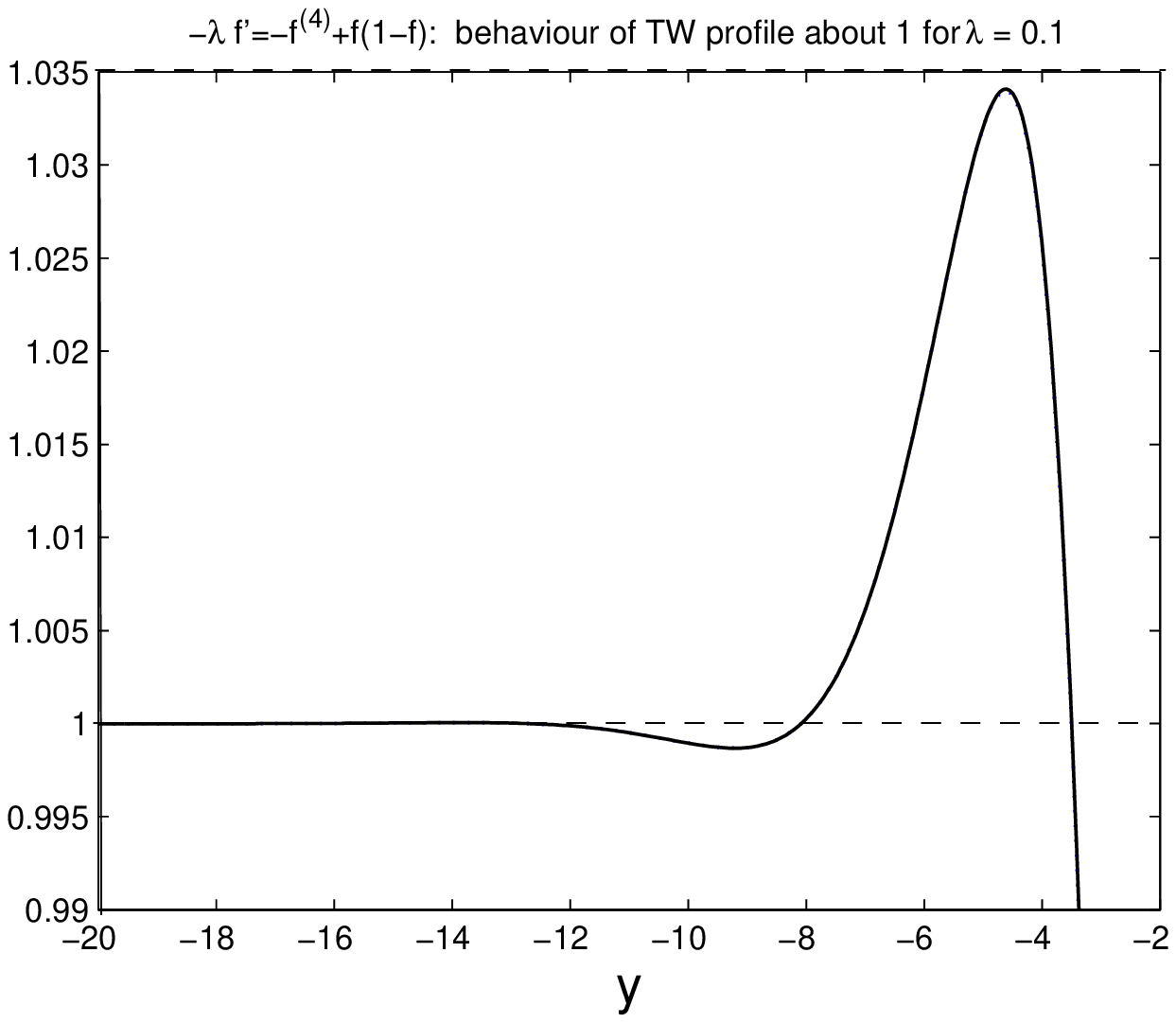}             
} \subfigure[scale $\sim 10^{-5}$]{
\includegraphics[scale=0.52]{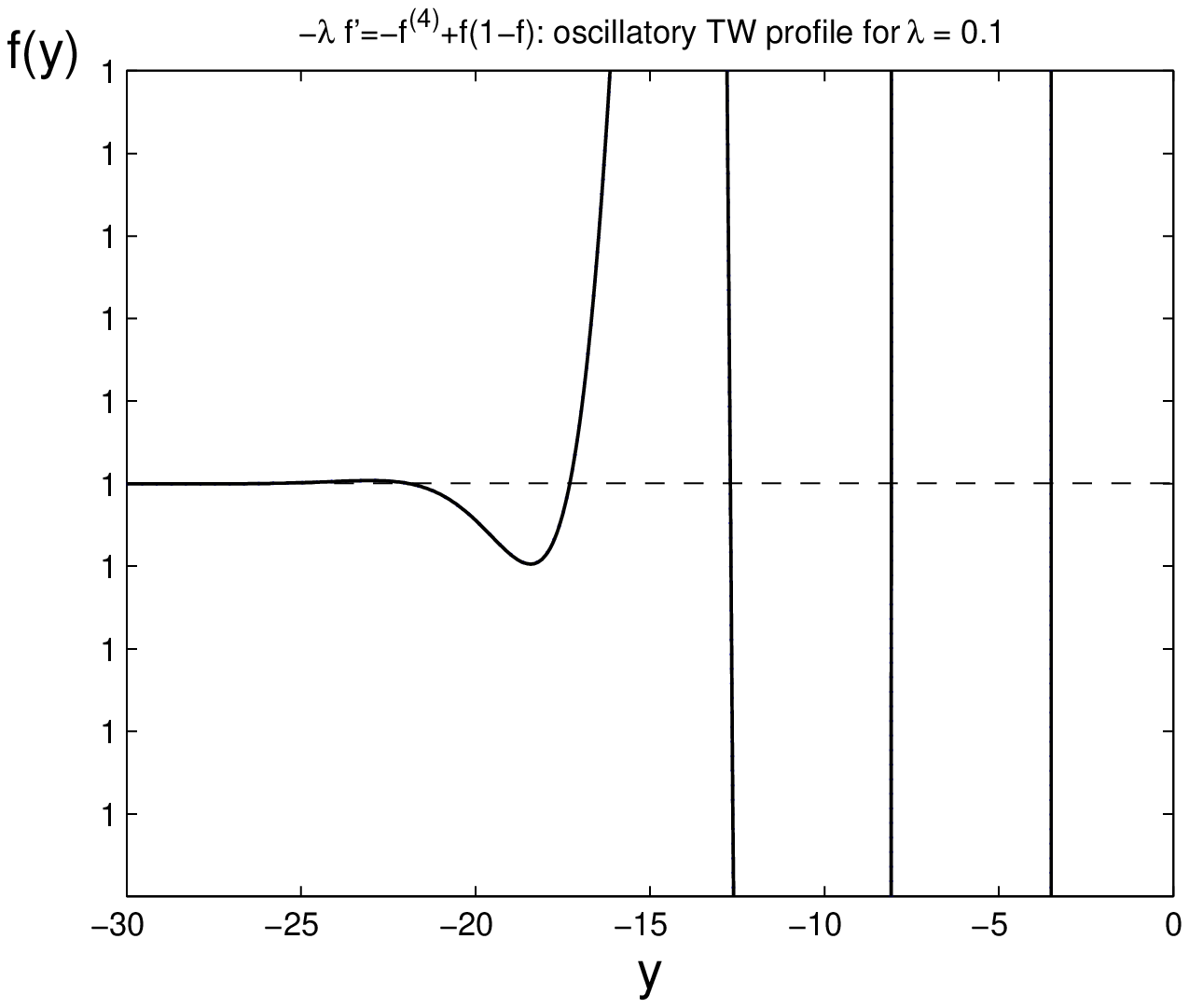}                        
}
 \vskip -.3cm
\caption{\rm\small Illustration of an oscillatory
 convergence of the TW profile $f(y)$ to 1 for $y \ll
-1$; $\l=0.1$.}
 \label{F3}
\end{figure}

 Figure \ref{F3} confirms oscillatory convergence to 1 in the
 opposite limit,
 as $y \to -\iy$. While (a) still not that convincing and reveals
 just a single intersection with the equilibrium $f=1$, the (b),
 in the scale $10^{-5}$, shows further five intersections with 1. On
 better scales, we can detect a dozen of intersections at least,
 so the convergence to 1 is clearly oscillatory.


\begin{figure}
\centering
\includegraphics[scale=0.85]{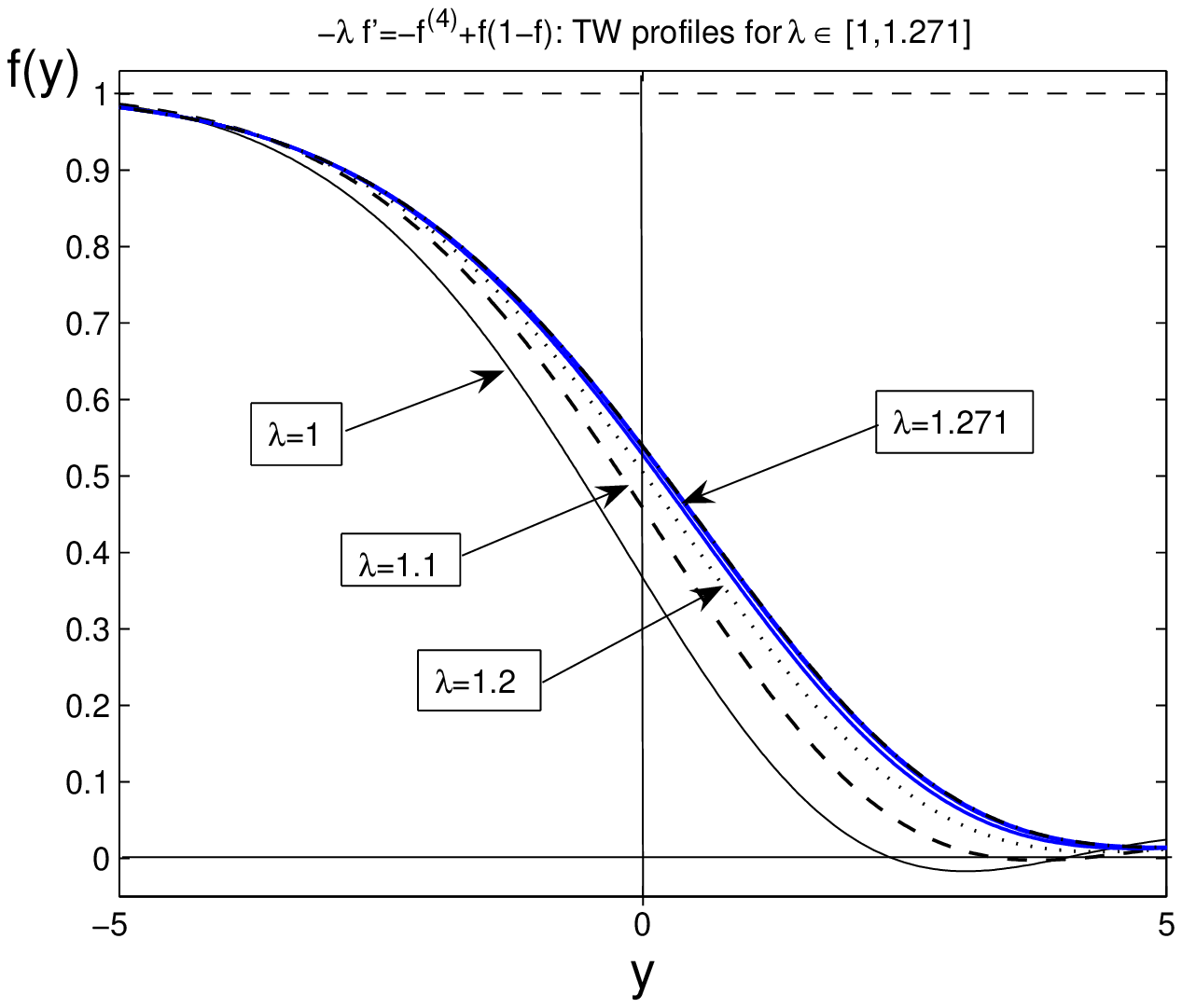}  
\vskip -.3cm
  \caption{The TW profiles $f(y)$ of
(\ref{E5}), \ef{BC1} for
 $\l \in [1, 1.271]$.}
 \label{F4}
\end{figure}

\begin{figure}
\centering
\includegraphics[scale=0.85]{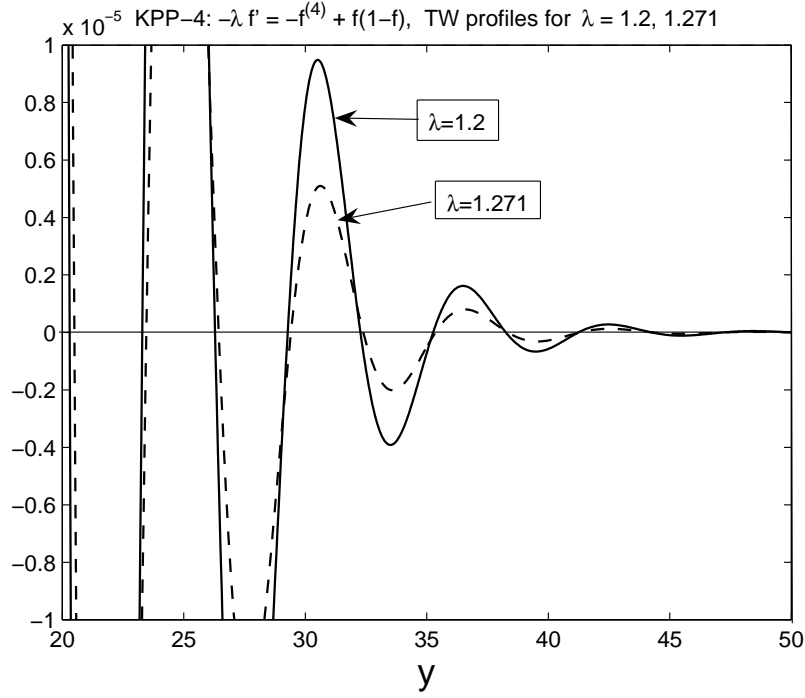}  
\vskip -.3cm
  \caption{The TW profiles $f(y)$ of
(\ref{E5}), \ef{BC1} is becoming less oscillatory  for $y \gg 1$
in the limit \ef{t4}.}
 \label{F5}
\end{figure}

Figure \ref{F4} shows  TW profiles for $\l$ larger than 1. The
phenomenon of a ``maximal" speed  \ef{SpeedMax} was clearly
observed, with the following sharp estimate:
 \be
 \label{t4}
 1.27148 \le \l_{\rm max}(2) < 1.27149.
 \ee
   In Figure \ref{F5}, we show the profiles $f(y)$ for $\l$'s close to
    $\l_{\rm max}$ in \ef{t4}.
 It turns out that, for $\l \approx \l_{\rm max}^-$, the TW
 profiles $f(y)$ remain oscillatory for $y \gg 1$, as Figure
 \ref{F5} confirms. Therefore, nonexistence of a $f(y)$ for
 $\l=\l_{\rm max}$ has nothing to do with  a possible reduction of
 the dimension (from 3D to 2D) of the asymptotic bundle as $y \to
 +\iy$, when two complex roots of the characteristic equation
 become real; cf. Theorem \ref{Th.1}(ii).

\begin{figure}
\centering \subfigure[$f(y)$]{
\includegraphics[scale=0.52]{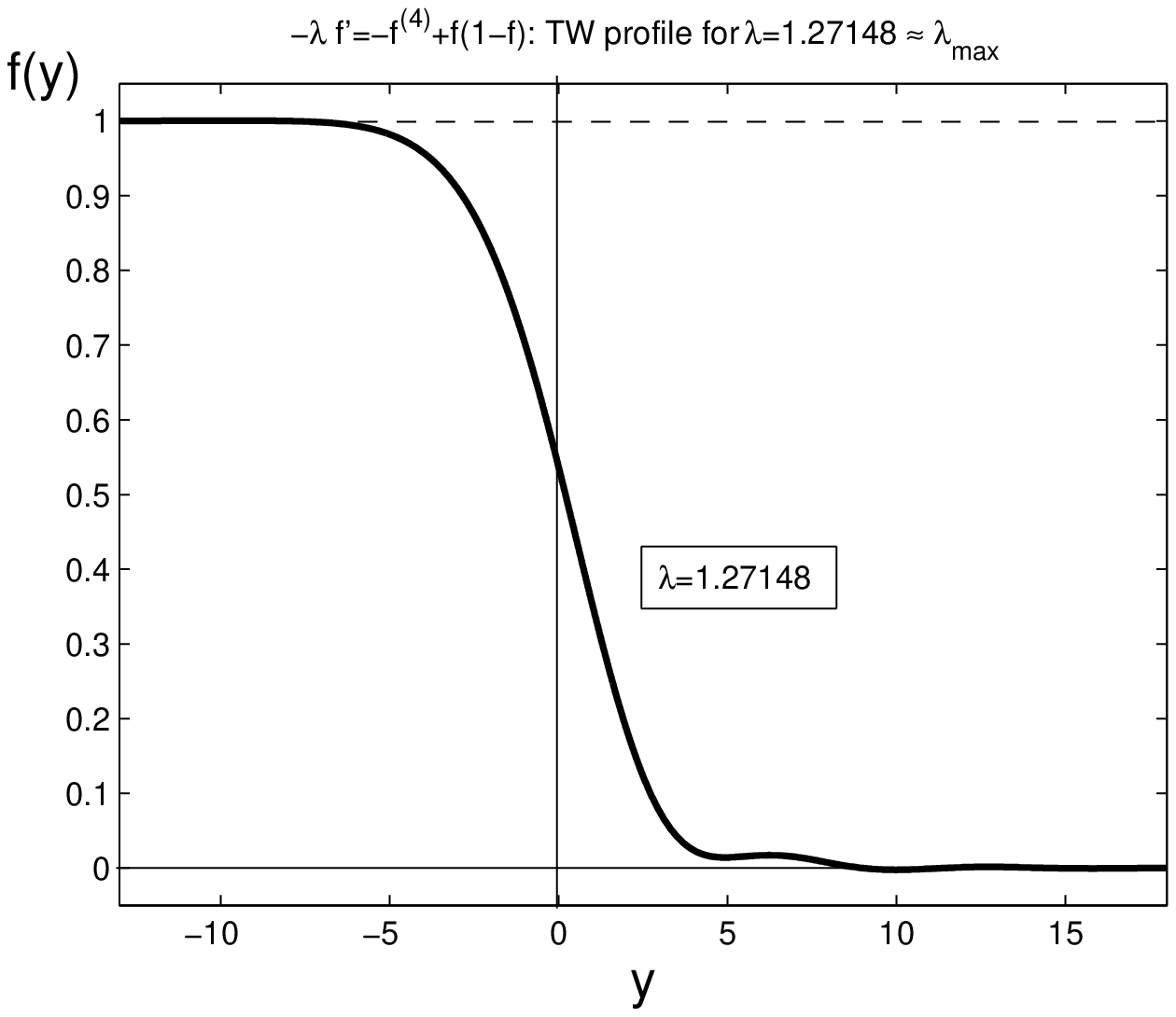}             
} \subfigure[oscillations for $y \sim 35$]{
\includegraphics[scale=0.52]{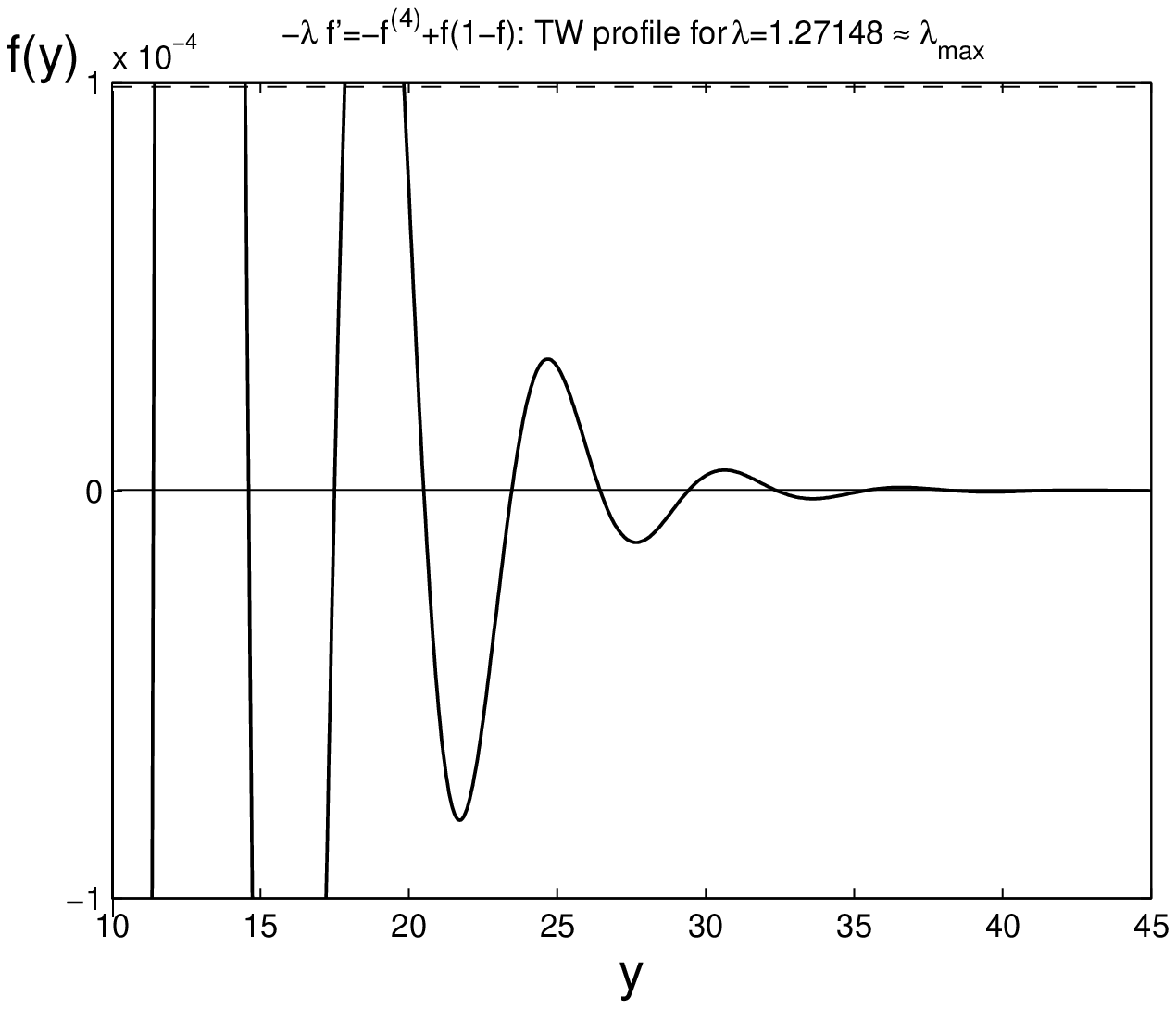}                        
}
 \vskip -.3cm
\caption{\rm\small Properties of the TW profile $f(y)$ for
$\l=1.27148 \approx \l_{\rm max}^-$.}
 \label{Fmax}
\end{figure}

In Figure \ref{Fmax}, we show the profile $f(y)$ for $\l=1.27148
\approx \l_{\rm max}^-$, (a), and its oscillatory behaviour for $y
\sim 35$, (b).




In correlation with the second half of the statement in
\ef{SpeedMax}, we observed no  TW profiles $Јf(y)$ for $\l>\l_{\rm
max}$. It should be noted that, since we always take the Heaviside
function \ef{Heav1} as the initial data for iterations via the
{\tt bvp4c}-solver, the non-convergence of the numerical scheme
could reflect {\em nonexistence} of TWs satisfying the ``minimal"
requirement,  with a ``maximal" decay at infinity, which is
available in the singular boundary conditions \ef{BC1}. In other
words, for $\l>\l_{\rm max}$, there could be other TW profiles
$f(y)$ with slower decay at infinity.

\subsection{The exponential bundle as $y \to + \iy$: linearized
analysis}

In the present semilinear ($n=0$) case, such an analysis is
reasonably easy, though, nowadays, a related stability theory
using Evans functions have been developed for more complicated
nonlinear systems; see \cite{Bert01}, with a full list of
references therein, with applications to general thin film flows.

\ssk

\begin{figure}
\centering
\includegraphics[scale=0.85]{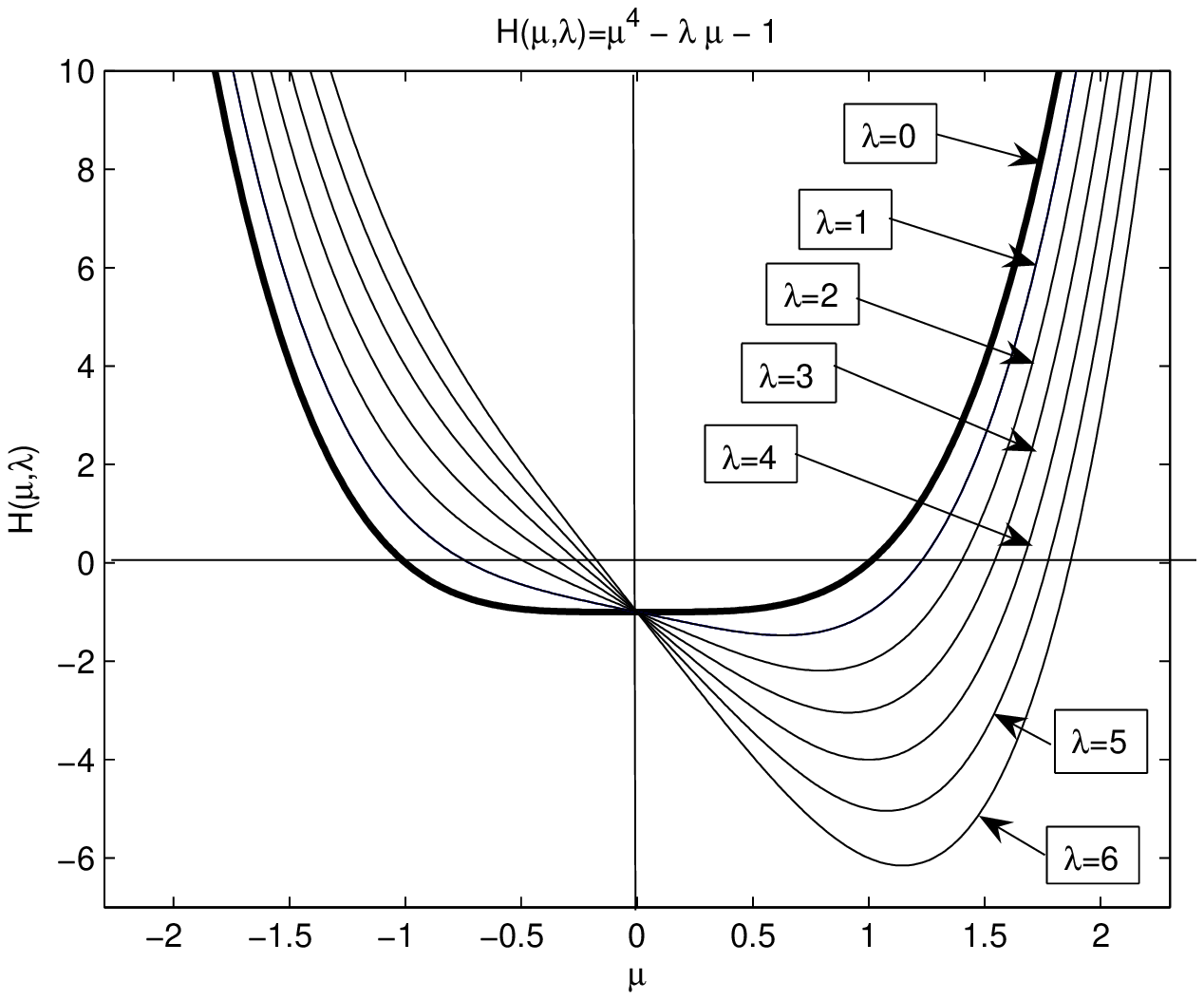}  
\vskip -.3cm
  \caption{The graphs of the function $H_+$ in \ef{t6} for various $\l$.}
 \label{F6}
\end{figure}

Consider the linearized about zero ODE in \ef{E5}:
 \be
 \label{t5}
 - \l g'=-g^{(4)} + g.
 \ee
 The corresponding characteristic polynomials and equation take  the form
  \be
  \label{t6}
  g(y)=\eee^{\mu y} \LongA H_+(\mu,\l) \equiv \mu^4 -\l \mu -1 =0.
   \ee
The graph of the function $H_+(\mu,\l)$ for various $\l$ is shown
in Figure \ref{F6}. One can see that it keeps a similar form for
all $\l >0$. It is seen that there exists a single real
$\mu=\mu_1(\l)>0$ (the unstable mode), a single $\mu_2(\l)<0$ (a
stable mode), and, for small $\l>0$, a complex stable mode with
$\mu_\pm(\l)$ having ${\rm Re} \, \mu_\pm(\l)<0$.

Indeed, for $\l=0$ in \ef{t6}, we get
 \be
 \label{t61}
 H_+(\mu,0) \equiv \mu^4-1=0 \LongA \mu_\pm(0)= \pm \ii,
  \ee
  together with the stable $\mu_1(0)=-1$ and the unstable
  $\mu_2(0)=1$.
The stable (unstable) roots for small $\l>0$ read
 \be
 \label{st11}
  \tex{
 \mu_{1}(\l)= -1- \frac \l{4 + \l}+O(\l^2) \quad \big(\mu_2(\l)= 1+ \frac \l{4 - \l}+O(\l^2)
  \big).
 }
 \ee
Similarly, the computations of complex roots of \ef{t6} yield
 \be
 \label{t7}
  \tex{
  \mu_\pm(\l)= - \frac {4\l}{\l^2+16}[1+O(\l)] \pm \ii  \big[1- \frac
  {\l^2}{\l^2+16}+O(\l^3)\big]
  + O(\l^3), \,\,\, {\rm Re} \, \mu_\pm(\l)<0.
  }
  \ee
 Together with \ef{st11}, this creates the whole 3D
  stable bundle of oscillatory linearized orbits, which we have
 seen in Figures \ref{F0} and \ref{F1}.

\begin{figure}
\centering
\includegraphics[scale=0.85]{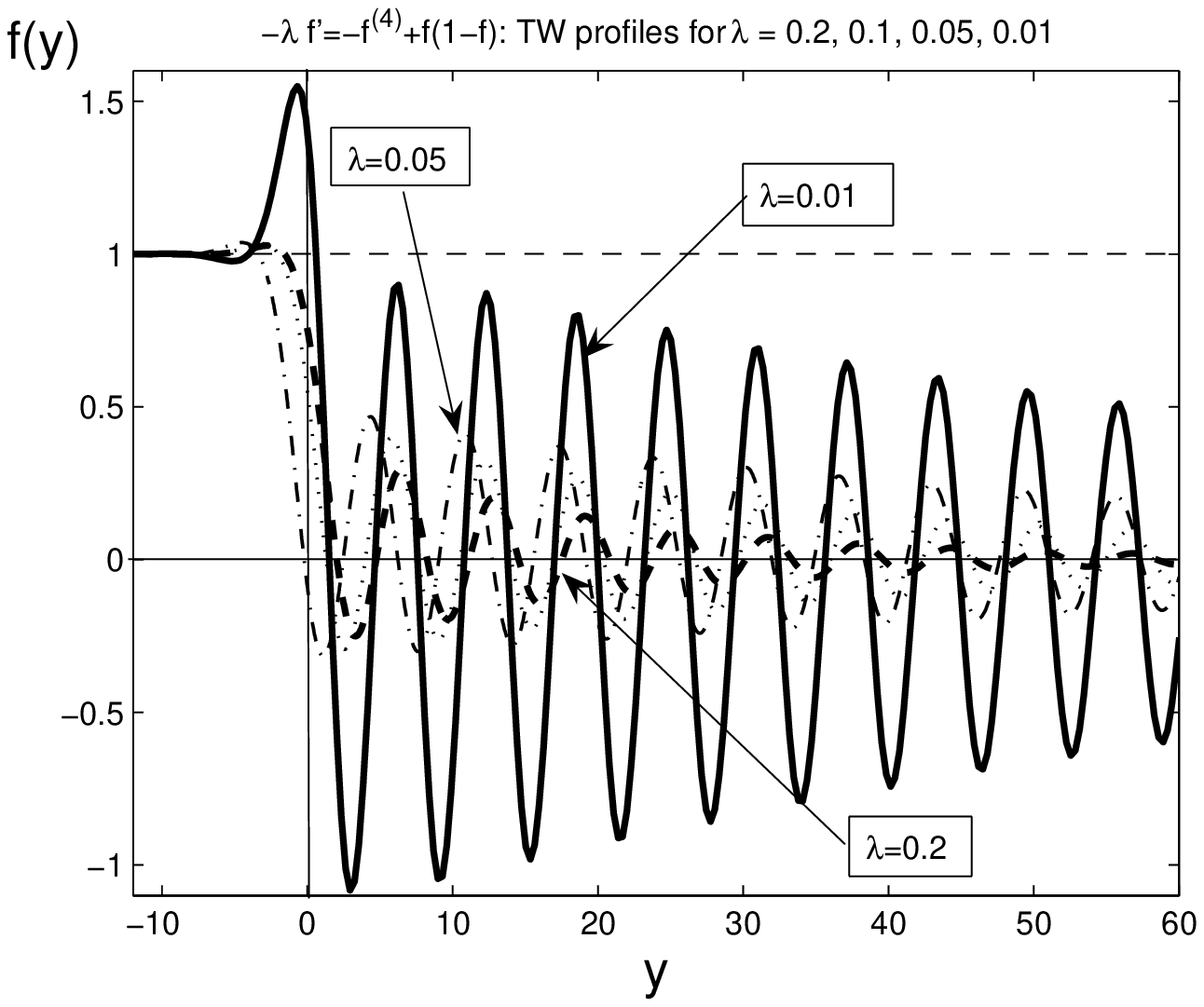}  
\vskip -.3cm
  \caption{The TW profiles $f(y)$ of
(\ref{E5}), \ef{BC1} for
 $\l=0.2$, 0.1, 0.05, and 0.01.}
 \label{F0small}
\end{figure}

Figure \ref{F0small} shows that, according to \ef{t7}, the TW
profiles get more and more oscillatory when $\l$ approaches 0,
i.e., when $\mu_\pm(\l)$ tend to pure imaginary roots in \ef{t61}.
The graph of $f(y)$ for  $\l=0.01$ in Figure \ref{F0small} was
obtained by integration over the interval $[-50,6000]$; see also
further results in the next subsection.


\ssk

 The above ``local" results can be made global by using the following
 conclusion, that differs the KPP--4 problem from the KPP--2 one,
 for which multiple characteristic zeros take place at the
 minimal speed 2; see \ef{norm2}.

 \begin{proposition}
  \label{Pr.mult}
  For any $\l \in \re$, the characteristic polynomial $\ef{t6}$ does
  not admit a root of the algebraic multiplicity $2$.
   \end{proposition}

\noi{\em Proof.} Assume the contrary: there exist two coinciding
roots $\bar \mu \in{\mathbb C}$, so that
 $$
 \mu^4- \l \mu -1 \equiv (\mu-\bar\mu)^2(\mu^2+a \mu+b), \quad
 \mbox{with some constants $a,b \in{\mathbb C}$}.
 $$
 Equating the coefficients of $\mu^3$, $\mu^2$, $\mu$, and 1 yields
 the system
  \be
  \label{sys33}
   \left\{
   \begin{matrix}
     - 2 \bar \mu+a=0, \quad\,\,\,\,\,\,\\
     \bar \mu^2-2a \bar \mu+b=0,\\
     a\bar \mu^2-2b \bar \mu=-\l,\,\,\\
     b\bar \mu^2=-1.\qquad\quad\,
    \end{matrix}
    \right.
    \ee
Obtaining $a=2 \bar \mu$ from the first equation and $b=- \frac
1{\bar \mu^2}$ from the last one and substituting into the second
one, we have
 $$
 \tex{
 3 \bar \mu^2 + \frac 1{\bar \mu^2}=0,
 }
 $$
 so that $\bar \mu$ is a complex root. Then the third equation
 reads
  $$
  \tex{
  \frac 32\, \l \bar \mu=-2,
  }
  $$
  so that then $\l \in\re$ must be complex. \quad $\qed$

  \ssk

  Thus, we claim the
 following result, which is sufficient for us for further
 applications:

\begin{proposition}
 \label{Pr.osc1}
{\rm (i)} At least for all small $\l>0$, the linearized equation
$\ef{t5}$ $($and hence the KPP--$4$ one in $\ef{E5}$ for $y \gg
1)$ admits a $3D$ stable family of oscillatory solutions as $y \to
+\iy$, and a $1D$ unstable one of exponentially divergent orbits.

{\rm (ii)} For small $\l<0$,  the stable bundle as $y \to +\iy$ is
$1D$, and the unstable one is $3D$.

 \end{proposition}

Using  Proposition \ref{Pr.mult}, one can improve the case (i) to
see that it might persist globally for $m=2$ for all $\l>0$.
However, for some larger $m$, this might be not the case globally.


\subsection{Oscillatory singular limit $\l \to 0^+$}

This has been discussed above. In fact, this is also a principle
question concerning {\em uniqueness} of $\l$-branches of
solutions. By careful numerical study of the behaviour of
$f(y;\l)$ as $\l \to 0^+$ we rule out a possibility to have a {\em
saddle-node} bifurcation  at some small $\l_{\rm s-n}>0$, at
which, via such a turning point, there appear two branches of TW
profiles for $\l>\l_{\rm s-n}$.

\begin{figure}
\centering \subfigure[on $(-20,250)$]{
\includegraphics[scale=0.52]{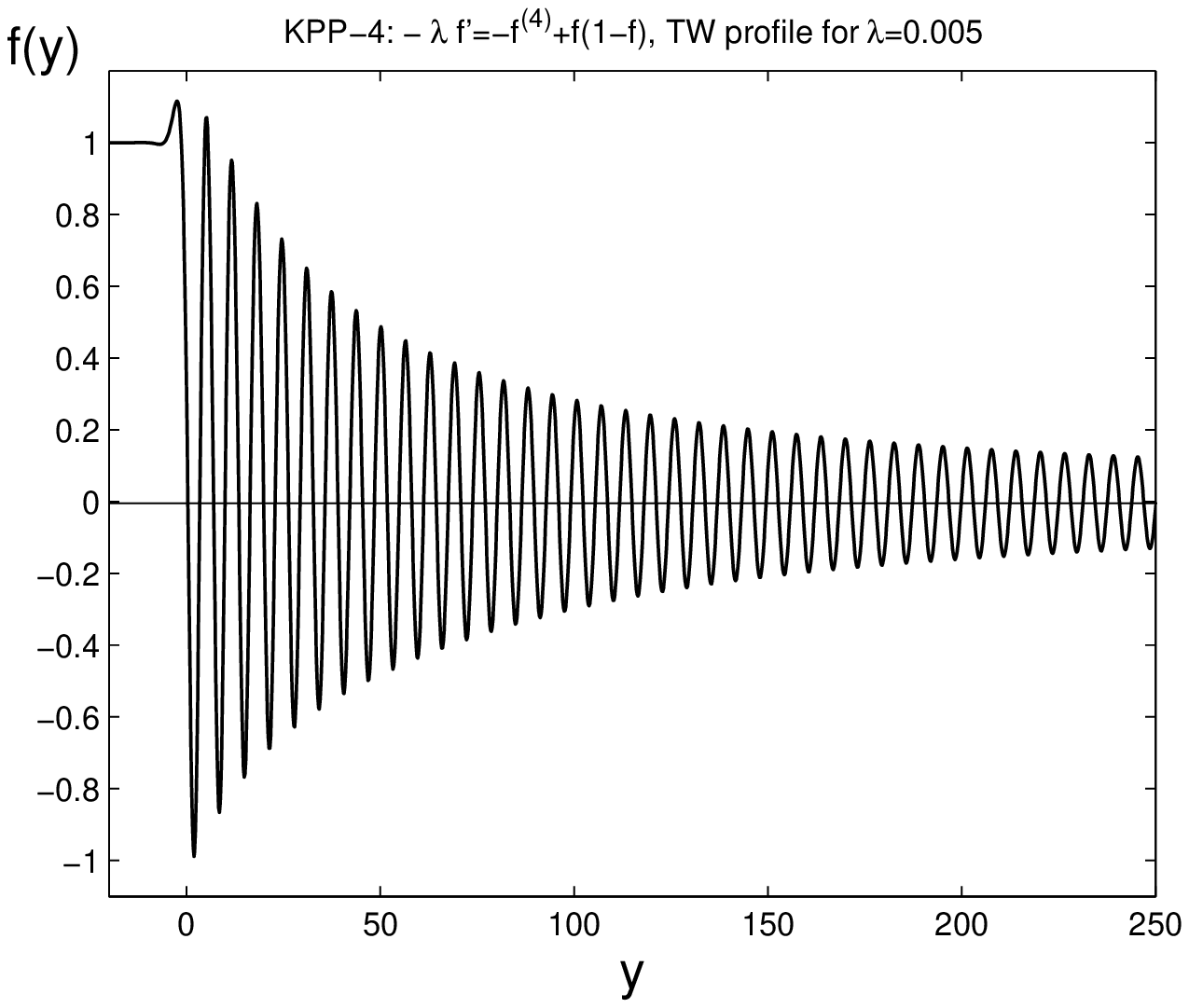}             
} \subfigure[on $(-20,1500)$]{
\includegraphics[scale=0.52]{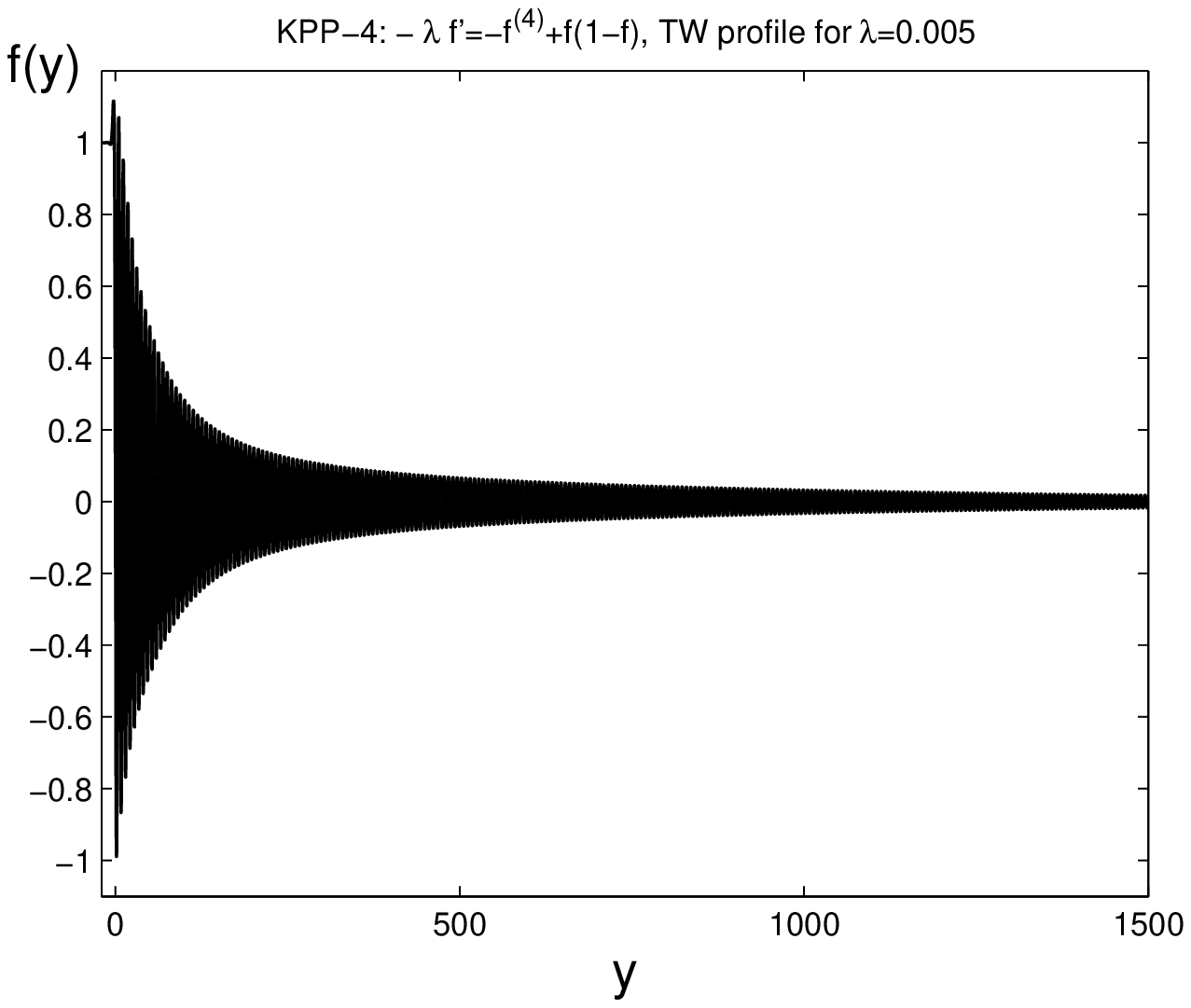}                        
}
 \vskip -.3cm
\caption{\rm\small The highly oscillatory TW profile $f(y)$ for
$\l=0.005$.}
 \label{Fl0.005}
\end{figure}

In Figure \ref{Fl0.005}, we show the TW profile $f(y)$ for
$\l=0.005$, which gets very oscillatory for $y>0$. We thus claim
that:
 \be
 \label{nosaddle}
  \begin{aligned}
 &\mbox{there is no a saddle-node bifurcation for small $\l>0$ and}\\
 & \mbox{the $\l$-branch ends up as $\l \to 0^+$ in a singular
 (oscillatory) limit.}
  \end{aligned}
  \ee
Therefore, in general, nothing contradicts so far the existence of
a single $\l$-branch of TW profiles. Of course, this well
correlates with the nonexistence of a TW at $\l=0$ (Proposition
\ref{Pr.Lambda1}) and some asymptotic computations above.

\subsection{The 2D stable bundle as $y \to -\iy$}

Then setting $f=1+g$ in \ef{E5} and linearizing yield the
following characteristic equation:
\be
\label{t9}
 - \l g'=g^{(4)} + g \andA
  g(y)=\eee^{\mu y} \LongA H_-(\mu,\l) \equiv \mu^4 -\l \mu +1 =0.
   \ee
Therefore, for $\l=0$, we have 2D stable (${\rm Re}\,(\cdot) >0$)
and unstable (${\rm Re}\,(\cdot) <0$) bundles with the roots
 \be
 \label{t10}
  \tex{
   \mu_\pm(0)= \pm \frac {1+\ii}{\sqrt 2} \andA
\bar\mu_\pm(0)= \pm \frac {1-\ii}{\sqrt 2}. }
 \ee
 By continuity, for all small $\l>0$, there exists 2D stable
 manifold of the equilibrium 1 with the roots
  \be
  \label{root11}
   \tex{
 \mbox{stable bundle}: \quad  \mu_+(\l) \andA \bar \mu_+(\l).
  }
  \ee

Unlike Proposition \ref{Pr.mult}, there exist multiple roots of
the polynomial in \ef{t9}, that follows from \ef{sys33} with the
last equation reading as $...=1$.

\begin{proposition}
\label{Pr.multminus} For the polynomial in $\ef{t9}$, there exist
two cases of double roots:
 \be
 \label{t9.1}
 \bar \mu_\pm = \pm 3^{-\frac 14} = \pm 1.3161... \forA \l_\pm= \pm 4 \cdot 3^{-\frac 34}
 =\pm 9.1180...\, .
 \ee
 \end{proposition}

 By \ef{t4}, the positive value $\l_+$ in \ef{t9.1} (at
 which the bundles may exchange their dimensions) plays no role.
Both bundles via \ef{t10} persists for small $\l>0$ ($\l<0$), so:

\begin{proposition}
 \label{Pr.osc2}
 At least for all small $|\l|>0$, the linearized equation
 $\ef{t9}$
 $($and hence the KPP--$4$ one in $\ef{E5}$ for $f \approx 1$ for $y
\ll -1)$ admits a $2D$ stable family of oscillatory solutions as
$y \to -\iy$, and a $2D$ unstable one of exponentially divergent
orbits.


 \end{proposition}

By Proposition \ref{Pr.multminus}, this result is global in
$\l>0$, but this does not help to explain why at $\l_{\rm max}$
(``minimal") TW profiles cease to exist.

\subsection {Local blow-up to $-\iy$}

In order to verify the global continuation properties of stable
bundles, one needs to check whether the {\em nonlinear} ODE
\ef{E5} admits blow-up and the dimension of such an unstable
manifold. To this end, we re-write it down and, as usual, neglect
the linear lower-order terms, which by standard local interior
regularity are negligible for $|f| \gg1$:
 \be
 \label{t12}
  f^{(4)}=-f^2 + f +\l f' = - f^2(1+o(1)) \asA f \to \iy.
   \ee
The unperturbed equation has the following exact solution:
 \be
\label{t13}
 \tex{
  f^{(4)}=-f^2 \LongA  f_0(y)=- \frac {840}{(Y_0-y)^{4}}
\to - \iy \asA y \to Y_0^+,
 }
  \ee
  where $Y_0 \in \re$ is a fixed arbitrary blow-up point. Studying the
  dimension of its stable manifold, we consider a perturbed
  solution and substitute it into the full equation \ef{t12} to get, on some linearization again:
   \be
   \label{t14}
    \tex{
   f=f_0+\e, \,\, |\e| \ll |f| \LongA (Y_0-y)^4 \e^{(4)}- 1680 \e = -840 + \frac{3360
   \l}{Y_0-y} +... \, .
   }
    \ee
One can see that the leading solution consists of balancing the
first two and the last term, that gives uniquely
 \be
 \label{t15}
  \tex{
  \e(y)= - \frac{140 \l}{69} \, \frac 1{Y_0-y}+... \, .
  }
  \ee
  Hence, in this expansion, no extra arbitrary parameter occurs.
Thus, we arrive at:

 \begin{proposition}
 \label{Pr.Blow1}
 For the ODE $\ef{E5}$, the stable manifold of blow-up solutions is
 1D depending on a single parameter being their blow-up point $Y_0 \in
 \re$.

  \end{proposition}

\subsection{Some general conclusions}
 \label{S2.5}

The KPP--4 ODE is  difficult for an analytic study in the
corresponding four-dimensional phase-space, so we do not plan here
to perform  its detailed and complete study.

However, the above local and blow-up analysis allows us to
formulate some important consequences. Recall that equation
\ef{E5} comprises analytic nonlinearities only, so by classic ODE
theory \cite{CodL}, all dependencies  on parameters in our problem
are given by analytic functions.

\begin{theorem}
 \label{Th.1}
 Fix a $\l \in \re_+$ in the ODE problem $\ef{E4}$, $\ef{E5}$:

 {\rm (i)} If the dimension of the stable bundle as $y \to +\iy$ is
 $3D$ $($as in Proposition $\ref{Pr.osc1}(i))$, and $2D$ as $y \to
 -\iy$ $($as in Proposition $\ref{Pr.osc2})$, then the number of
 solutions $f(y)$ of the ODE problem $\ef{E5}$, $\ef{BC1}$ is not more
 than finite.

 {\bf (ii)} If the stable dimension of the bundle as $y \to + \iy$
 changes to $2D$,
 or even becomes $1D$
 $($as in Proposition $\ref{Pr.osc1}(ii)$ for $\l<0)$, then, most
 probably\footnote{Actually, in a natural sense, with probability
 zero. Moreover, for $\l<0$, this means the actual nonexistence of TWs, as Proposition \ref{Pr.Lambda1} states.}, the KPP--$4$ problem does not have a solution.

 \end{theorem}

 \noi{\em Proof.} (i) Denote by $C_{\rm st}^{(1,2,3)}$ the
 real coefficients that control the 3D stable bundle at $y\to +
 \iy$, and by $D_{\rm unst}^{(1,2)}$ the coefficients of the 2D
 unstable bundle as $y \to - \iy$. Bearing in mind all three Propositions above,
  to get a global solution $f(y)$,
 one has to solve the following algebraic system:
  \be
  \label{t20}
   \left\{
    \begin{matrix}
     Y_0(C_{\rm st}^{(1,2,3)})=-\iy, \qquad\qquad\,\,\, \ssk
     \\
    D_{\rm unst}^{(l)}(C_{\rm st}^{(1,2,3)})=0, \quad l=1,2,
    \end{matrix}
    \right.
    \ee
    where the first equation means that blow-up does not occur at
    finite points and hence ``$Y_0=-\iy$" is to be understood in a
    limit sense (the solution remains uniformly bounded). Hence, we
    arrive at a system of three analytic equations with three
    unknowns $C_{\rm st}^{(1,2,3)}$. Then number of solutions is
    not more than countable, and the set of solutions is discrete with a possible accumulation point at infinity
    only.

    Moreover, by our blow-up analysis above, it follows that,
    since the stable modes as $y \to +\iy$ are linearly
    independent, blow-up always occurs, if these coefficients $C_{\rm
    st}^{(1,2,3)}$ take sufficiently large values (a stability property of blow-up for quadratic ODE
    \ef{t12}). Therefore, all the solutions of \ef{t20} belong to
    a large ``box" in $\re^3$, where only a finite number of those
    are available.

    \ssk

    (ii) The result is obvious: then the corresponding analytic
    system becomes overdetermined   (more equations than the unknowns), so it is
    inconsistent with the probability 1, meaning that its
    solvability(if any)
    can be ``accidental" only and then is destroyed by a.a. arbitrary
    perturbations
    of any of the coefficients/functions involved (no structural
    stability). $\qed$

 \subsection{Existence: shooting from $y=-\iy$}

 Actually, the system \ef{t20} corresponds to ``shooting from the
 right to the left-hand side" that was not principal for such an algebraic analysis.
  However, for a complete proof it is
 easier to shoot from $y=-\iy$. Indeed, keeping the same notations
 of ``stable" expansion coefficients, we then arrive at two
 algebraic
 equations only:
 \be
  \label{t201}
   \left\{
    \begin{matrix}
     Y_0(D_{\rm st}^{(1,2)})=+\iy, \ssk
     \\
    C_{\rm unst}(D_{\rm st}^{(1,2)})=0,
    \end{matrix}
    \right.
    \ee
 where we have used the fact that, at $y =+\iy$, there exists a 1D
 unstable manifold governed by the unique coefficient $C_{\rm
 unst}$. Thus, we obtain a system of two equations with two
 unknowns $D_{\rm st}^{(1,2)}$, and such a system admits a rather
 straightforward analysis on this plane. We refer, as a typical example, to
  \cite[\S~4]{EGW1}, where a rather similar geometric-type analysis
  applied to a fourth-order ODE arising from a semilinear Cahn--Hilliard equation.
See also \cite[\S~7]{Gl4} concerning ODEs in thin film theory.

The strategy is as follows. For a fixed $D_{\rm st}^{(1)}$, we use
the second parameter $D_{\rm st}^{(2)}$ to satisfy the first
equation in \ef{t201}. It is clear that this is always possible,
since in both limits  $D_{\rm st}^{(2)} \to \pm \iy$, the solution
$f(y;D_{\rm st}^{(1,2)})$ blows up, so that it is global for at
least one value named $\hat D_{\rm st}^{(2)}(D_{\rm st}^{(1)})$.
As the next step,  we start to vary $D_{\rm st}^{(1)}$ to get rid
of the unstable component as $y \to +\iy$, i.e., to satisfy the
second equation in \ef{t201}. In view of a  clear oscillatory
character of all the solutions as $y \to +\iy$, this always can be
done.


Thus, the existence of a proper TW profile only is not a principal
difficulty, but {\em uniqueness remains a hard open problem
still}. We do not know any mathematical reason (in particular, a
kind of a {\em monotonicity} property, usually leading to
uniqueness), for which it must be unique. But numerically, we
multiply observed convergence to the same profile from variously
perturbed initial data as a starting point of nonlinear
iterations.


\section{The  KPP--(6,1) problem}
\label{S3}

We now, more briefly than above, consider the KPP--(6 ,1) problem
\ef{E6} and its ODE counterpart \ef{E7}, \ef{BC1}. Of course,
\ef{Lam1} remains valid, so we check existence of TWs for $\l>0$
only.

\subsection{Numerical results}



 \begin{figure}
\centering
\includegraphics[scale=0.85]{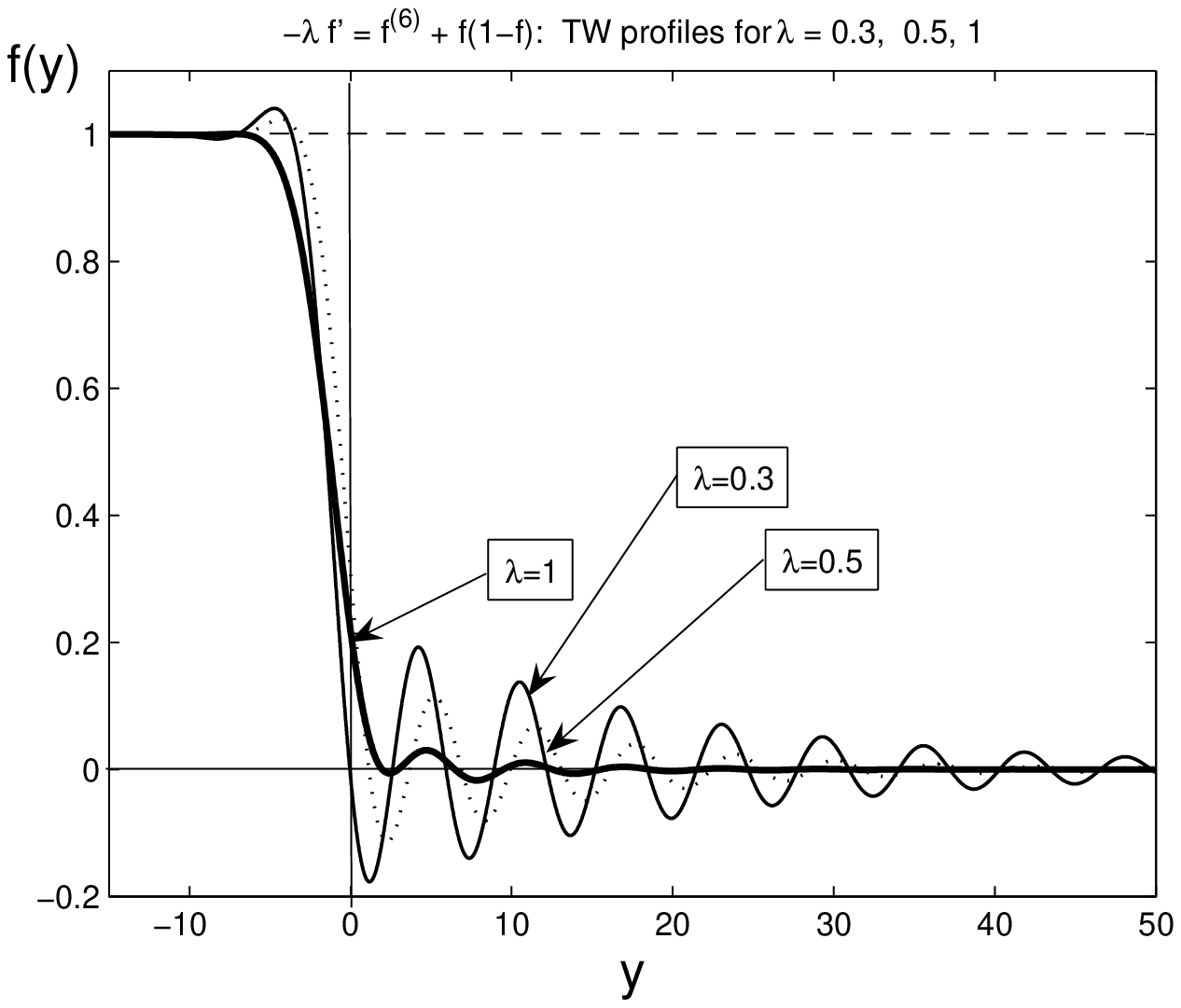}  
\vskip -.3cm
  \caption{The TW profiles $f(y)$ satisfying
(\ref{E7}), \ef{BC1} for
 $\l=0.3$, 0.5, and 1.}
 \label{F06}
\end{figure}

\begin{figure}
\centering
\includegraphics[scale=0.85]{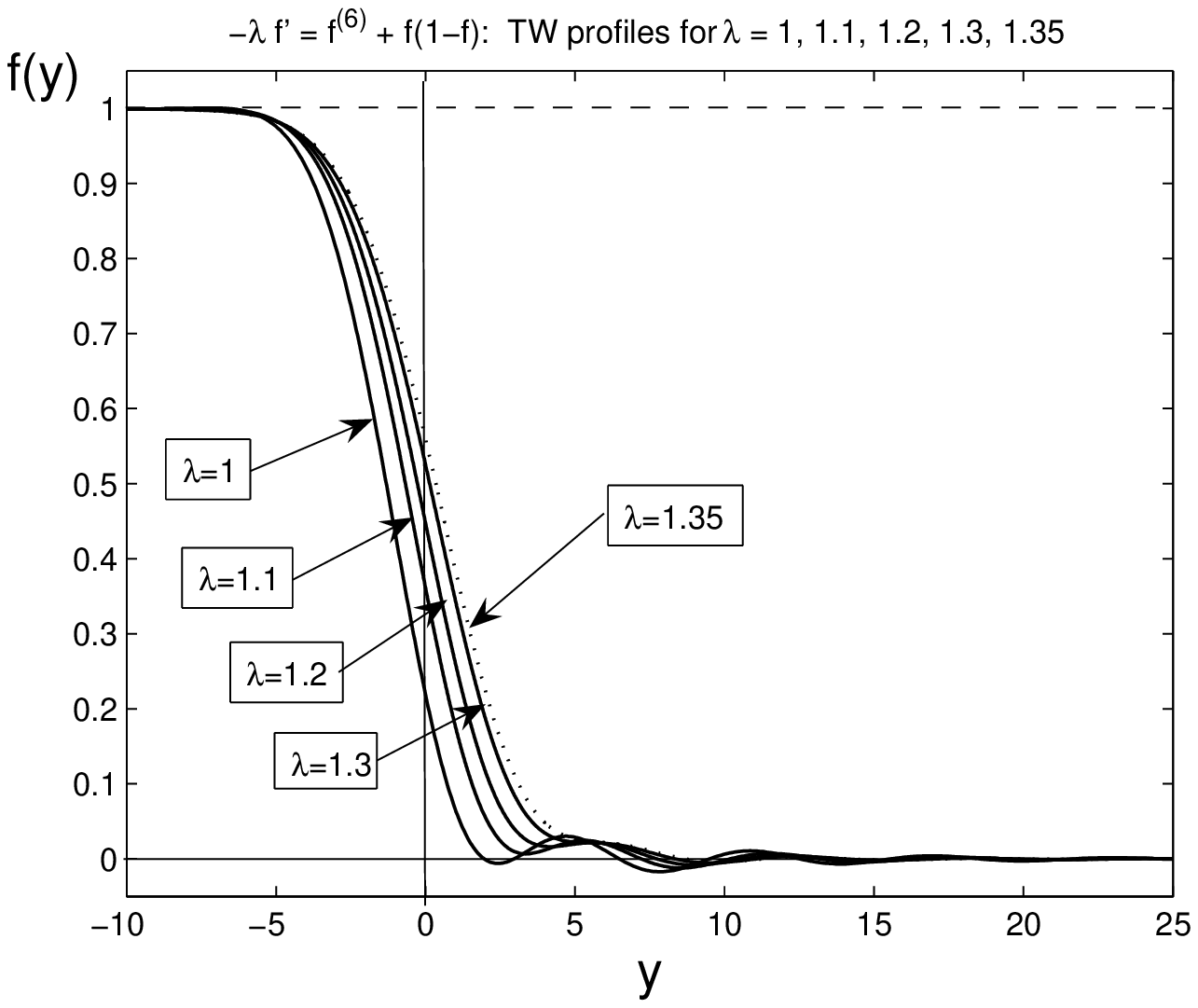}  
\vskip -.3cm
  \caption{The TW profiles $f(y)$ of
(\ref{E7}), \ef{BC1} for
 $\l=1$, 1.1, 1.2, 1.3, and 1.35.}
 \label{F16}
\end{figure}
\begin{figure}
\centering \subfigure[TW profiles]{
\includegraphics[scale=0.52]{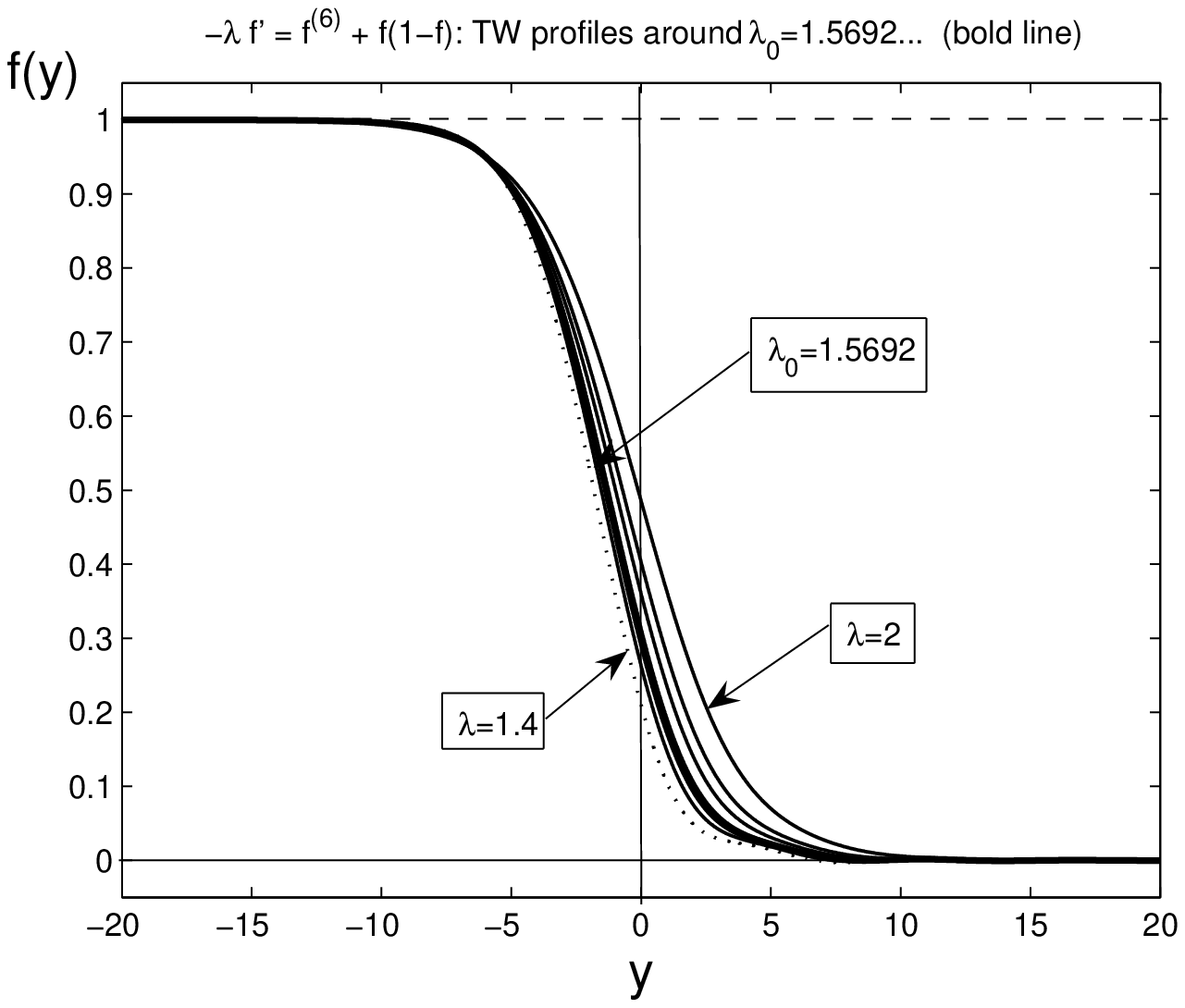}             
} \subfigure[oscillations for $y \gg 1$]{
\includegraphics[scale=0.52]{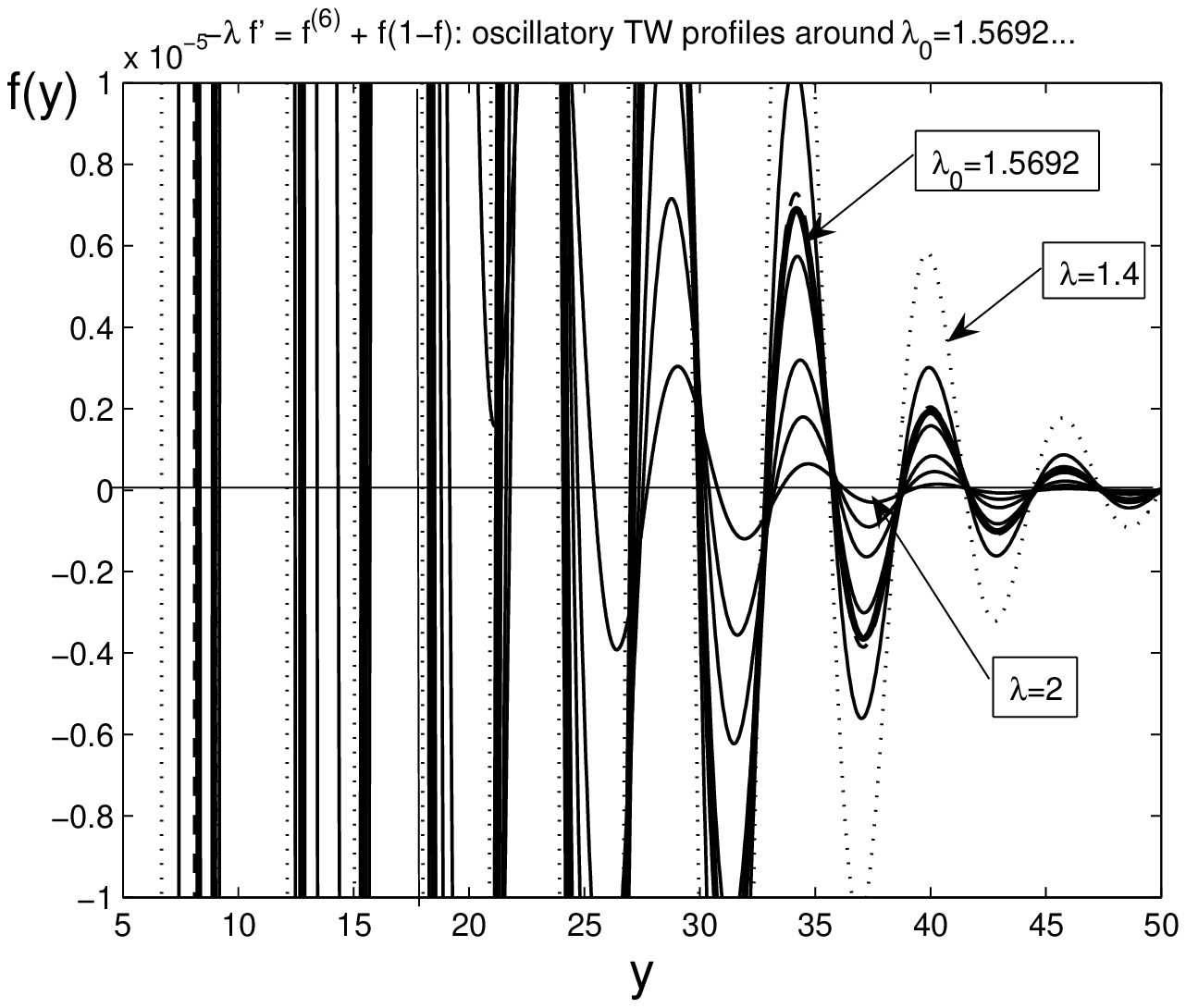}                        
}
 \vskip -.3cm
\caption{\rm\small TW profiles in the KPP--6 problem for $\l >1$.}
 \label{F36}
\end{figure}

Figure \ref{F06} shows TW profiles for $\l=0.3, \,\,0.5,$ and 1.
In Figure \ref{F16},  TW profiles are shown for larger  speeds $\l
\in [1, 1.35]$. All of them are clearly oscillatory as $y \to +
\iy$. Next, in Figure \ref{F36}, we show TW profiles for larger
$\l>1.4$.




\begin{figure}
\centering \subfigure[$f(y)$]{
\includegraphics[scale=0.52]{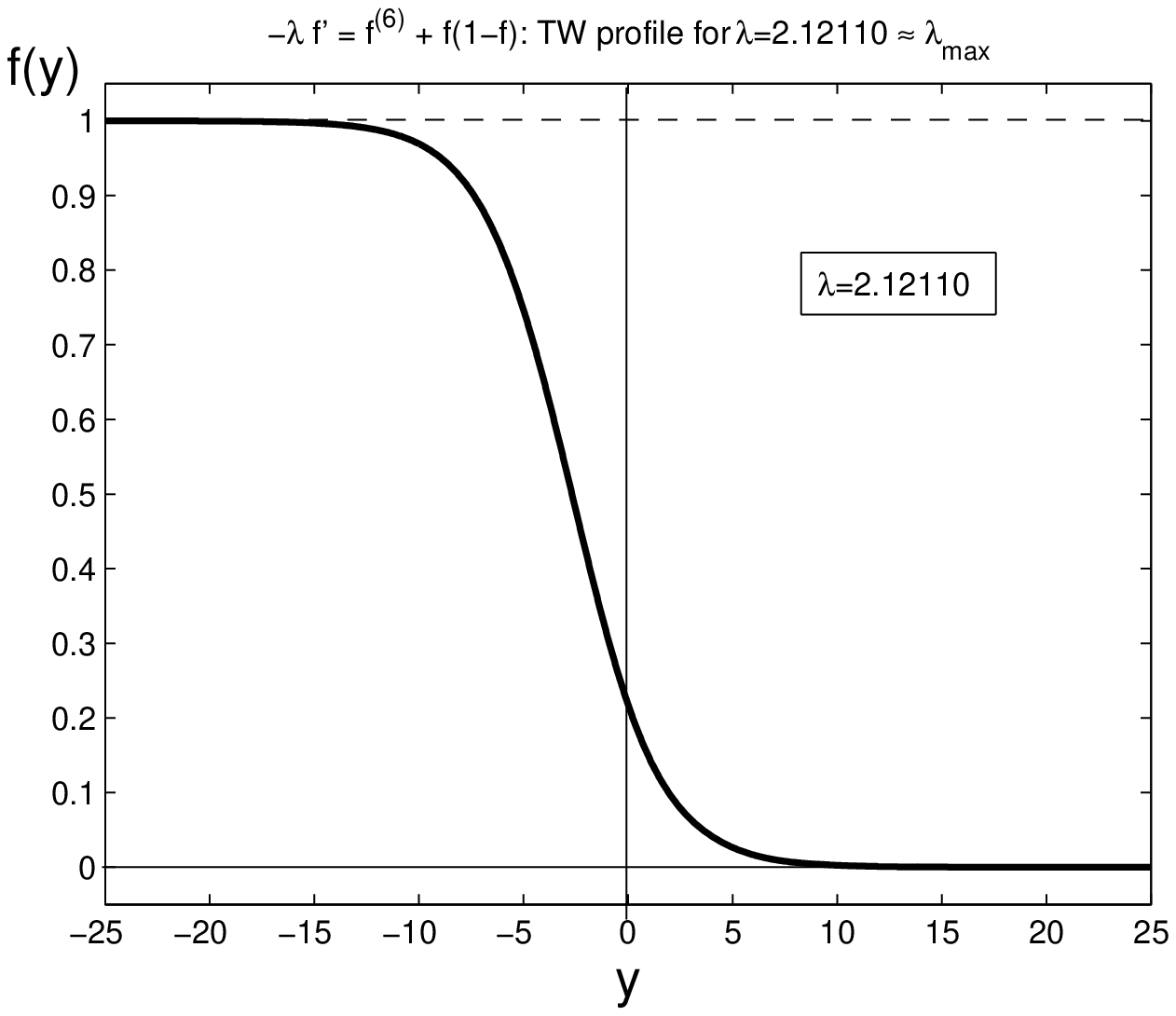}             
} \subfigure[oscillations for $y \sim 30$]{
\includegraphics[scale=0.52]{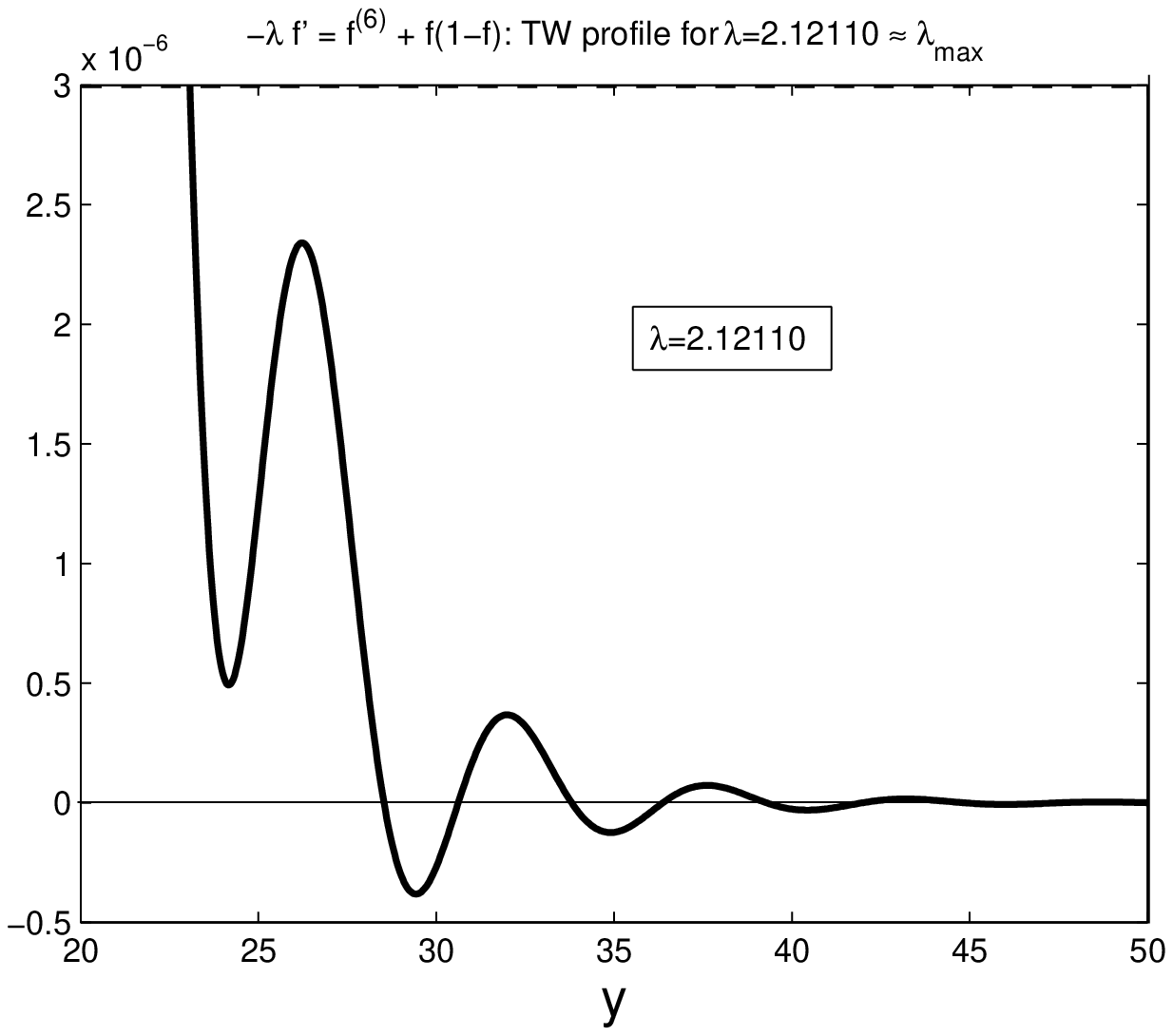}                        
}
 \vskip -.3cm
\caption{\rm\small Properties of $f(y)$ for $\l=2.12110$.}
 \label{Fmax3}
\end{figure}

 Figure \ref{F36}(a) confirms an oscillatory convergence to 1 in the
 opposite limit,
 as $y \to -\iy$.
Figure (b) explains  the same phenomenon as for the KPP--4 (cf.
\ef{t4}), but here
 \be
 \label{ll1}
m=3: \quad  2.12110 \le \l_{\rm max}(3) < 2.12111... \, ,
 \ee
 at which TW profiles cease to exist. For $\l=2.12110$, the TW
 profile and its oscillations for $y \gg 1$ are shown in Figure \ref{Fmax3}.


 \begin{figure}
\centering \subfigure[convergence to 1]{
\includegraphics[scale=0.52]{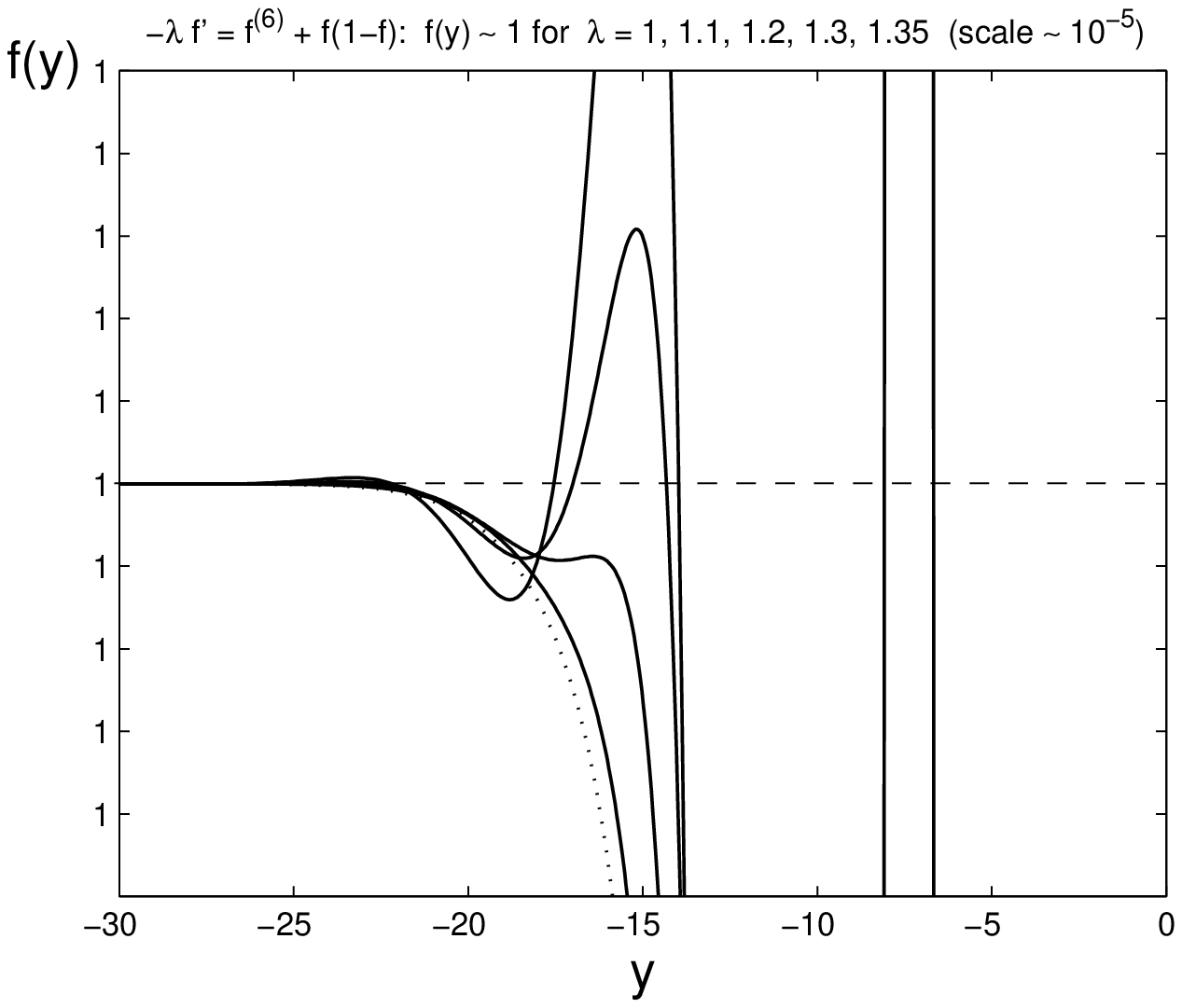}             
} \subfigure[oscillations for $\lambda_*=1.121...$]{
\includegraphics[scale=0.52]{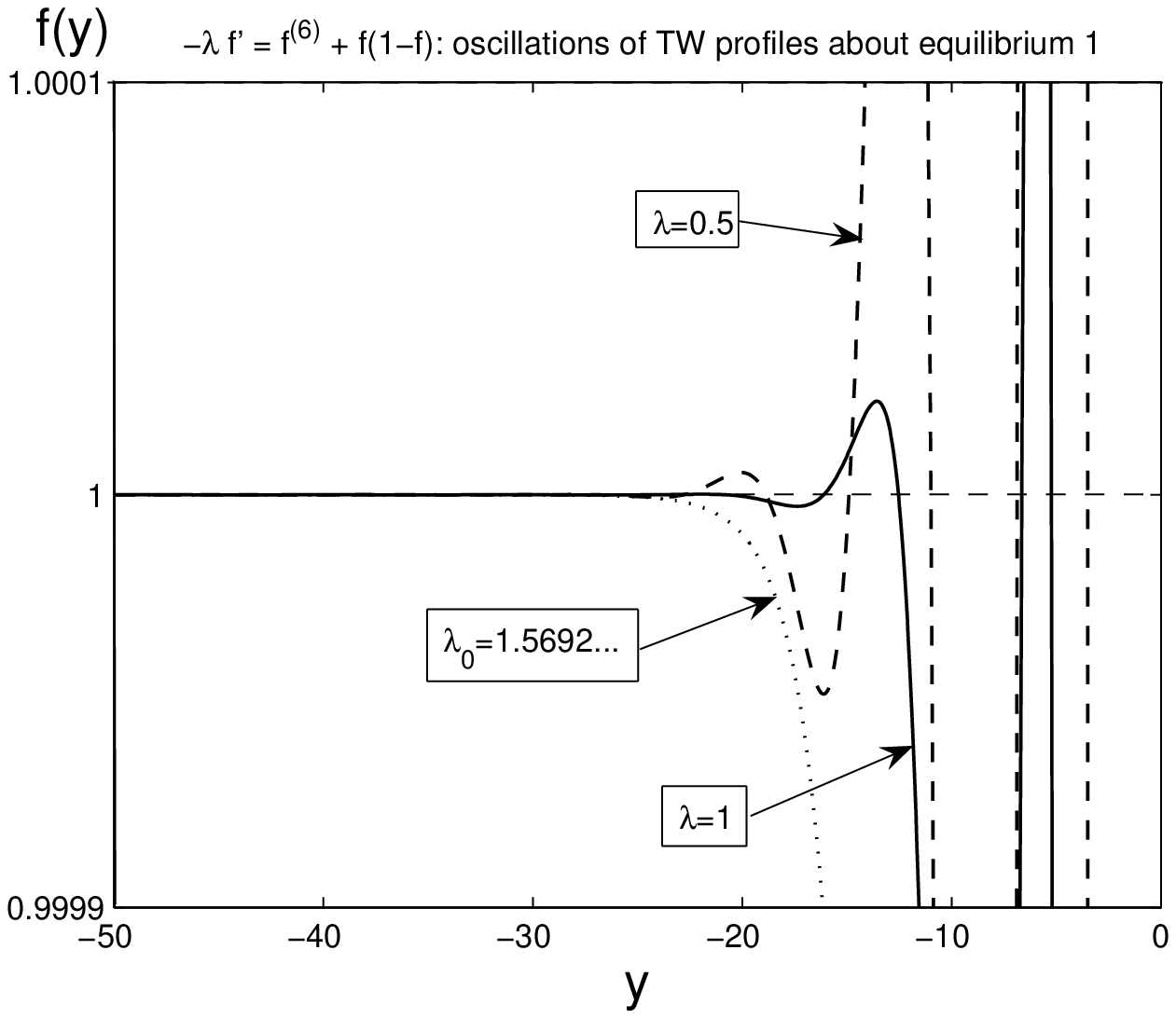}                        
}
 \vskip -.3cm
\caption{\rm\small Oscillations of the TW profile $f(y)$ as $y \to
\pm \iy$.}
 \label{F26}
\end{figure}

\subsection{Dimensional analysis of the bundles and towards
existence}

\begin{figure}
\centering
\includegraphics[scale=0.85]{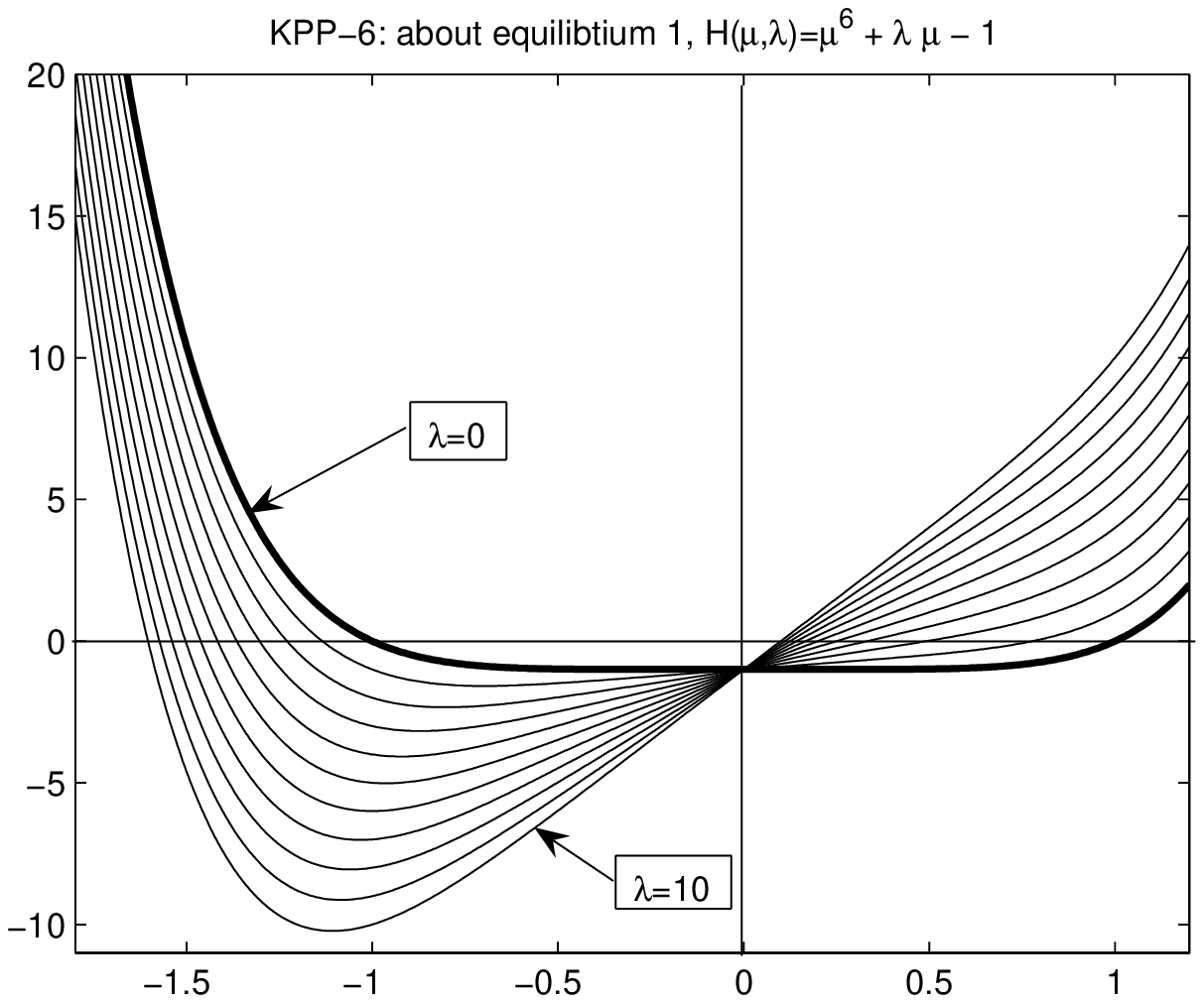}  
\vskip -.3cm
  \caption{The graphs of the function $H_+$ in \ef{t66} for various $\l=0,1,2,...,10$.}
 \label{F66}
\end{figure}

This is not essentially more difficult than in the KPP--4 problem.
Namely, consider the linearized about zero ODE in \ef{E7}:
 \be
 \label{t56}
 - \l g'=g^{(6)} + g.
 \ee
 The corresponding characteristic polynomial is now
  \be
  \label{t66}
  g(y)=\eee^{\mu y} \LongA H_+(\mu,\l) \equiv \mu^6 +\l \mu -1 =0.
   \ee
The graph of the function $H_+(\mu,\l)$ for various $\l$ is shown
in Figure \ref{F66}. It is quite similar to that for the
KPP--$(4,1)$. Moreover, Proposition \ref{Pr.mult} holds for the
polynomial in \ef{t56}.


Similarly, we study the bundles as $y \to - \iy$, where there
occurs
 \be
 \label{pol32}
H_+(\mu,\l) \equiv \mu^6 +\l \mu +1 =0.
 \ee
 It admits double roots
  \be
  \label{pol33}
  \bar \mu_\pm= \pm 5^{- \frac 16}= \pm 0.7647... \quad \mbox{at}
  \quad \l_\pm = \mp 2 \cdot 5^{-\frac 56}= \mp 0.5231... \, ,
 \ee
 which, again, have nothing to do with the maximal speed in \ef{ll1}.

Then, for small $|\l|>0$, the dimensions of stable/unstable
manifolds correspond to a proper shooting, but with more
parameters. Therefore, while we can prove that the family of TW
profiles is not more than  countable (or even finite),
the actual proving of existence is more difficult for such a
multi-dimensional shooting.

Note also, that the blow-up problem (cf. the ODE \ef{t13})
 \be
 \label{bl1}
 f^{(6)}=f^2
 \ee
 is now more difficult since \ef{bl1} admits oscillatory
 solutions, which can be studied as in \cite[\S~7]{Gl4} by
 introducing an {\em oscillatory component} represented by a
 periodic function.

 \section{Higher-order problems: KPP--(8,1) and KPP--(10,1)}
\label{S810}

 Those are two examples to illustrate a possibility of such KPP
 extensions. By \ef{Lam1}, one needs to consider $\l>0$ only.

\begin{figure}
\centering
\includegraphics[scale=0.85]{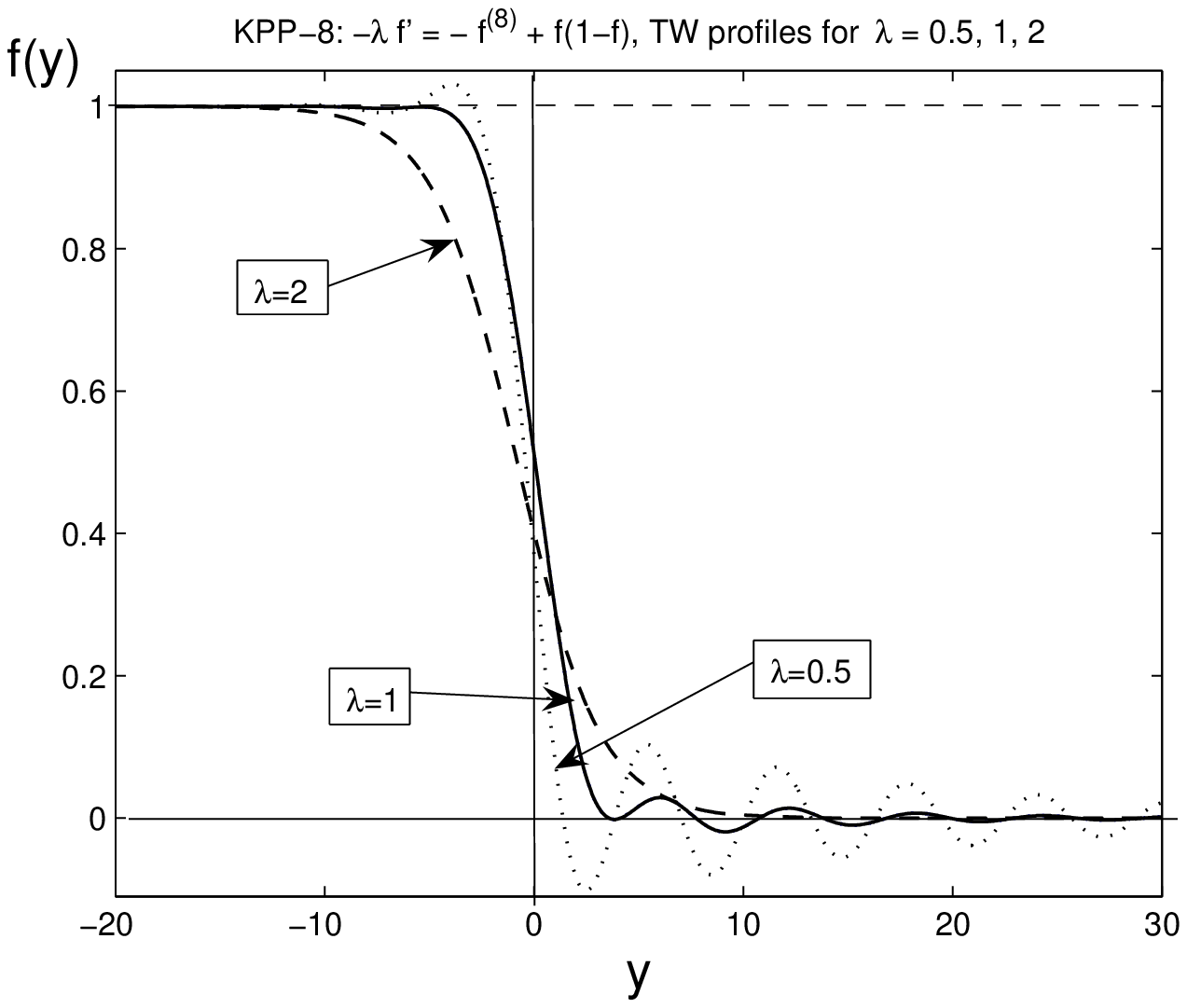}  
\vskip -.3cm
  \caption{The TW profiles $f(y)$ satisfying
(\ref{8.1}), \ef{BC1} for
 $\l=0.5$, 1, and 2.}
 \label{F81}
\end{figure}

\begin{figure}
\centering
\includegraphics[scale=0.85]{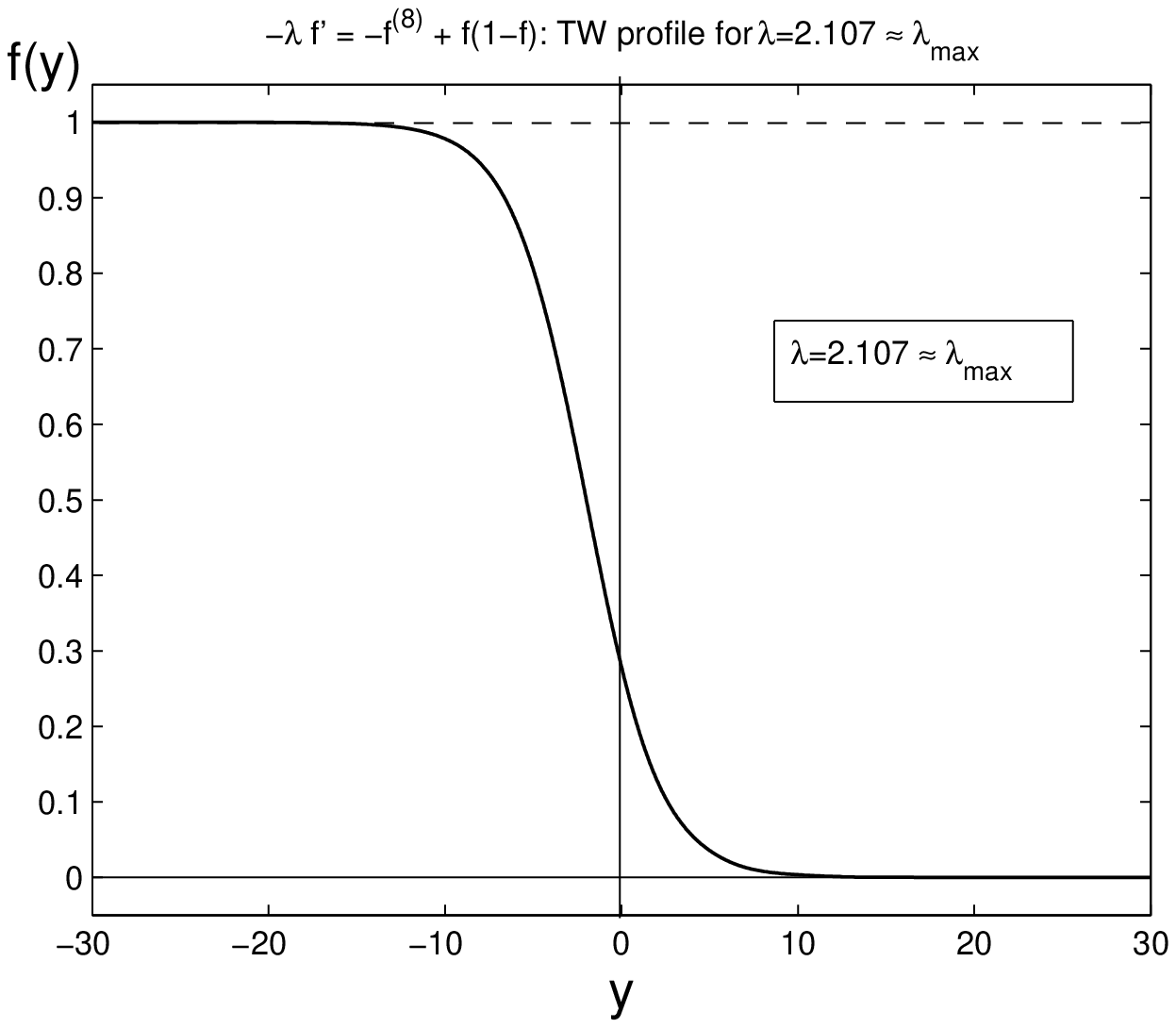}  
\vskip -.3cm
  \caption{The TW profile $f(y)$ satisfying
(\ref{8.1}), \ef{BC1} for
 $\l=2.107$.}
 \label{Fmax41}
\end{figure}

\begin{figure}
\centering
\includegraphics[scale=0.85]{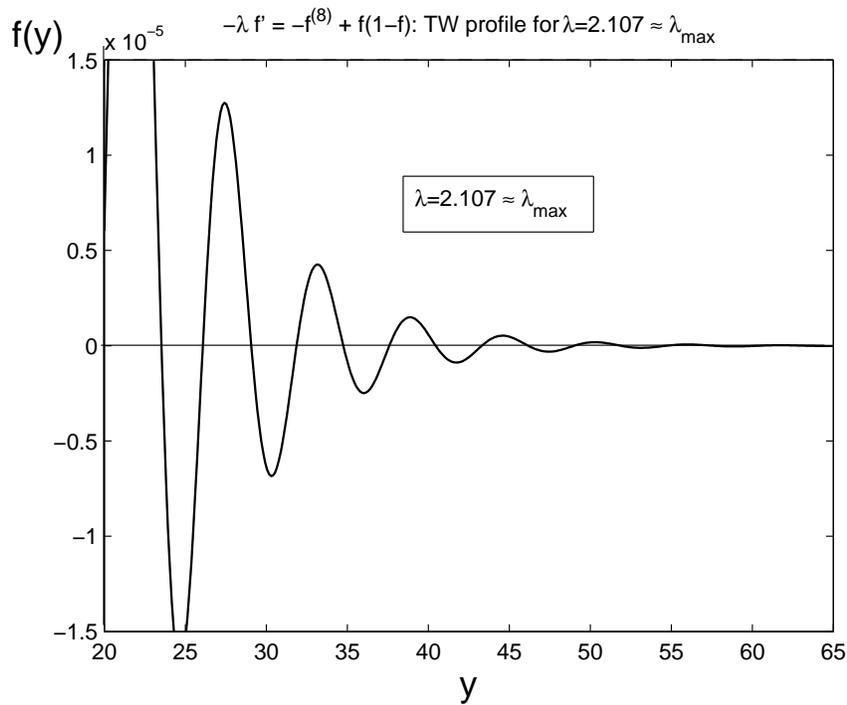}  
\vskip -.3cm
  \caption{Oscillations of the TW profile $f(y)$ from Figure \ref{Fmax41} for
 $\l=2.107$.}
 \label{Fmax42}
\end{figure}

 Thus, in the KPP--(8,1), we have the following ODE problem:
  \be
  \label{8.1}
 (\mbox{PDE})\quad u_t=-D_x^8 u+u(1-u) \LongA (\mbox{ODE}) \quad -
 \l f'=-f^{(8)}+f(1-f).
  \ee
  Figure \ref{F81} shows  TW profiles for $\l=0.5$, 1, and 2, satisfying the ODE \ef{8.1}
  with the singular boundary conditions \ef{BC1}. The critical
  existence speed value in the sense of \ef{SpeedMax} is now
   \be
   \label{m4max}
   m=4: \quad 2.107 \le \l_{\rm max}(4) < 2.108.
    \ee
Figure \ref{Fmax41} shows $f(y)$ for $\l=2.107$ from \ef{m4max},
while Figure \ref{Fmax42} explains its oscillations for $y \sim
35$.


\begin{figure}
\centering
\includegraphics[scale=0.85]{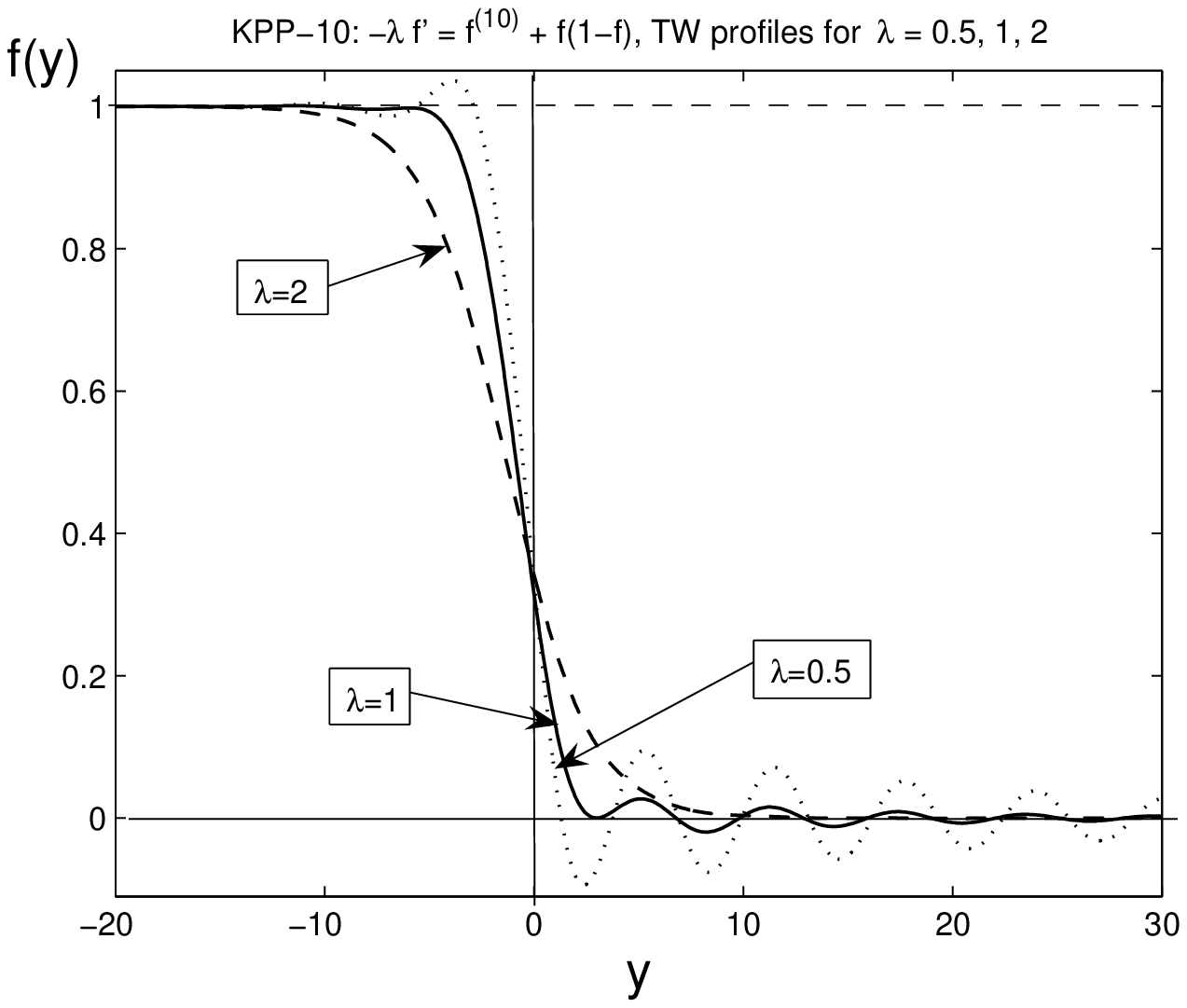}  
\vskip -.3cm
  \caption{The TW profiles $f(y)$ satisfying
(\ref{10.1}), \ef{BC1} for
 $\l=0.5$, 1, and 2.}
 \label{F101}
\end{figure}

In the KPP--(10,1), we deal with the following ODE problem:
  \be
  \label{10.1}
 (\mbox{PDE})\quad u_t=D_x^{10} u+u(1-u) \LongA (\mbox{ODE}) \quad -
 \l f'=f^{(10)}+f(1-f).
  \ee
  Figure \ref{F101} shows  TW profiles for $\l=0.5$, 1 and 2, satisfying the ODE \ef{10.1}
  with the  boundary conditions \ef{BC1}.


It is seen that the TW profiles for both problems are close to
each other. For instance, for $\l=0.5$, the difference is not more
than $2\%$. This confirms that it is possible to pass to the limit
$m \to +\iy$ in the $2m$th-order KPP--problem \ef{m1}, \ef{m2},
and construct an analogy of the KPP--problems for an
``infinite-order" PDE and ODE.

A local stability of bundles can be done for both above parabolic
models establishing for which $\l$'s the shooting problem is well
posed and  the corresponding algebraic systems are consistent.
Indeed, proving existence, to say nothing of the uniqueness, are
open problems.

\section{The origin of $\log t$-shift: centre subspace balancing}

 \label{SDiscr}

In this section, we explain the origin of the $\log t$-shifting
(retarding) from the TW for the PDE higher-order KPP--problem.

More correctly, the actual proof of such $\log t$-drift assumes
 a delicate matching of the solution behaviour on compact subsets in the TW
$y$-variable, i.e., in the {\em Inner  Region}, with a remote {\em
Outer Region} for $y \gg 1$, where the influence of the nonlinear
term $-u^2$ is negligible, and the actual behaviour is governed by
the linear bi-harmonic operator. However,
we do not present here this kind of a procedure of {\em matching}
of asymptotic expansions of those two regions, which is difficult
in the present quite general case.


\ssk

Thus, we consider a KPP-type problem for a semilinear PDE
 \be
 \label{K1}
  u_t=\AAA u+u(1-u) \inB \re \times \re_+,
  \ee
  where $\AAA$ is a proper homogeneous isotropic linear differential
  operator satisfying some extra conditions specified below.
 For instance, we may fix, as a simple example, the
  bi-Laplacian operator
   \be
   \label{Lapl22}
    \AAA=-D_x^4 \quad (\mbox{or $=D_x^6$, or $(-1)^{m+1}D_x^{2m}$, that is the same}).
   \ee

\subsection{Behaviour for small $t>0$: first appearance of $\log
t$}

Given the Heaviside data \ef{1.H}, to study \ef{K1}, \ef{Lapl22}
(with $m=2$, but any $m \ge 2$ fits quite similarly) for small
$t>0$, one needs to use the rescaled independent variables
 \be
 \label{var1}
 \tex{
  z= \frac x{t^{1/4}} \andA \t= \log t \to -\iy \asA t \to 0^+.
  }
  \ee
  Then the rescaled equation takes the form
   \be
   \label{var2}
    \tex{
    v_\t= \big(\BB- \frac 14 \, I\big) + \eee^\t \, v(1-v) \whereA \BB=- D_z^4+ \frac 14\, z
    D_z+ \frac 14 \, I.
     }
     \ee
     Denoting by $V(z)$ the unique stationary solution of the
     problem
 \be
 \label{var3}
  \tex{
 \big(\BB- \frac 14 \, I\big) V=0 \inB \re, \quad V(-\iy)=1, \quad V(+\iy)=0,
 }
  \ee
  and performing the standard linearization in \ef{var2} yield
   \be
   \label{var4}
    \tex{
v=V+w \LongA w_\t= \big(\BB- \frac 14\, I\big)+ \eee^\t\,
V(1-V)+...\,, }
 \ee
  provided that $|w| \ll 1$. Finally, in such an inner region, we
  obtain
   \be
    \label{var5}
    \tex{
    w(z,\t)=\eee^\t \var(z)+...
     \LongA \big(\BB- \frac 14\, I\big)\var= h(z) \equiv
     -V(z)(1-V(z)).
     }
     \ee
According to the linearized analysis in Section \ref{S2}, $h(z)$
has an exponential decay at infinity, so, keeping the leading
terms, we may assume that
 \be
 \label{var21}
 h(z)=B_- \eee^{b_- z}+... \asA z \to - \iy \andA h(z)=B_+
 \eee^{b_+ z}+... \asA z \to + \iy,
  \ee
  where $B_\pm, \, b_\pm \not =0$ are some constants and ${\rm Re}\,b_->0$
  and ${\rm Re}\, b_+<0$.

 Since $\BB$ is known to possess a discrete spectrum in a weighted
 $L^2_\rho$-space \cite{Eg4}, where
  \be
  \label{var6}
   \rho(z)=\eee^{a_0 |z|^{4/3}}, \quad a_0>0 \,\,\, \mbox{is a
   small constant},
    \ee
     and a complete set of eigenfunctions $\{\psi_k\}$,
the problem for $\var$ in \ef{var5} admits the unique solution
given by the eigenfunction expansion
 \be
 \label{var7}
  \tex{
  \var= \sum c_k \psi_k \whereA c_k= \langle h, \psi_k^*\rangle,
   }
   \ee
   $\{\psi_k^*\}$ being the eigenfunction set of generalized
   Hermite polynomials as eigenfunctions of the adjoint operator
    \be
    \label{var8}
     \tex{
     \BB^*=-D_z^4 - \frac 14 \, z D_z.
     }
     \ee
In \ef{var7}, $\langle \cdot, \cdot \rangle$ is the dual
$L^2$-metric, in which \ef{var8} represents the adjoint to $\BB$.
By \ef{var21}, all such integrals converge. Moreover, one can
check that
 \be
 \label{var9}
  \tex{
  \var(z)= \frac {4 B_\pm}{b_\pm z} \, \eee^{b_\pm z}+... \asA z \to \pm
  \iy.
  }
  \ee

Thus, we obtain the following expansion of the solution $u(x,t)$
for the Heaviside initial function:
 \be
 \label{var91}
  \tex{
  u(x,t)=V\big(\frac x{t^1/4}\big)+ t \, \var\big(\frac
  x{t^1/4}\big)+...\,, \quad \mbox{in the domain, where}\,\,\,
  \big|t \, \var\big(\frac
  x{t^1/4}\big)\big| \ll 1.
  }
  \ee
Hence, the last linearization condition defines the domain of its
applicability: by \ef{var91},
 \be
 \label{var10}
  \tex{
  |t \, \var(z)| \ll 1 \LongA |t \eee^{b_\pm x/t^1/4}| \ll1 \LongA
  |x| \ll t^{\frac 14} \, |\log t|.
  }
  \ee
  Here, we first observe the appearance of the $\log t$-``shift" (now, a factor) in the rescaled variable
  coming from the inner expansion region. However, a further extension of such
  an expansion towards a proper TW-frame (to see how this $\log t$-multiplier is actually transformed into the $log t$-shift)
   is a difficult procedure, which we cannot treat here.
  Instead, below, we will catch this $\log t$-shift in a simpler way,
  starting now directly from the TW-frame.

\subsection{TW-frame: linearization and rescaled equation}

We now perform a perturbation analysis in the TW-frame. Then,
   we need to assume that the corresponding ODE problem
   \be
   \label{K2}
   -\l_0 f'= \AAA f+f(1-f),
   \ee
   with the conditions \ef{BC1}, admits a unique solution $f$.
   Then, attaching the solution $u(x,t)$ to the
    front moving and setting, for convenience, $x_f(t) \equiv \l_0t -g(t)$,
    the PDE reads
    \be
    \label{K21}
    u(x,t)=v(y,t), \quad y=x- \l_0 t+g(t) \LongA v_t= \AAA v  + \l_0 v_y + v(1-v)
    - g'(t) v_y.
     \ee
   We next linearize \ef{K21} by setting
    \be
    \label{K3}
    v(y,t)=f(y)+w(y,t),
     \ee
     that yields the following perturbed equation:
      \be
      \label{gg2}
      w_t= \BB w -g'(t) f' - g'(t) w_y -w^2, \,\,\,\mbox{where}
      \,\,\,
      \BB w= \AAA w + (1-2f)w
       + \l_0 w_y.
       \ee
 Assuming that, in this $g(t)$-moving frame, there exists the
 convergence as in \ef{TW3}, so that $w(t) \to 0$ as $t \to +\iy$,
 one can see that the leading non-autonomous perturbation
in \ef{gg2} is the second term on the right-hand side. However, as
we show, the last two terms, though negligible, will define a
proper $\log t$-shift of the front.

\subsection{Where $\log t$-shift comes from}
\label{Sect.5.3}

First of all, an extra ``centre manifold"-like front shifting can
be connected with the known results on the instability of
travelling waves in a number of analogous fourth and sixth-order
parabolic equations; see \cite{Gao2004, Strauss2004, Li2012} and
references there in. However, of course, this does not directly
imply specifically the $log t$-shifting to be justified. In
addition, it was always shown that, typically, the corresponding
linearized operators (like $\BB$) possess {\em continuous}
spectra, the fact that essentially complicates our further
analysis.

In this connection, note also that the rescaled equation \ef{gg2}
is essentially {\em non-autonomous} in time, so we can hardly use
powerful  tools of nonlinear semigroup theory; see \cite{Lun}.
However, using  a formal asymptotic approach, we will trace out
some typical features of a {\em centre subspace} behaviour, after
an extra rescaling and balancing of non-autonomous perturbations.

\ssk

 Thus, as usual (see Introduction), we assume that $g'(t) \to 0$ as
 $t \to +\iy$ sufficiently fast, i.e., at least algebraically, so
 that
  \be
  \label{gg1}
  |g''(t)| \ll |g'(t)| \forA t \gg 1.
  \ee
Under the hypothesis \ef{gg1}, the only possible way to balance
{\em all} the terms therein (including the quadratic one $-w^2$)
for $t \gg 1$ is to assume the asymptotic separation of variables:
 \be
 \label{gg3}
 w(y,t) = g'(t) \psi(y)+ \e(t) \var(y)+... \whereA
 |\e(t)|
 \ll |g'(t)| \asA t \to + \iy.
  \ee
 Here, we omit higher-order perturbations. Substituting \ef{gg3}
  into \ef{gg3} yields
   \be
   \label{gg31}
    \begin{aligned}
    & g''(t) \psi + \e'(t) \var+...
   =  \,g'(t)(\BB \psi-f') \ssk \\
   & + \e(t)\BB
       \var -(g'(t))^2(\psi'+\psi^2)-g'(t)\e(t)(\var'+2\var \psi)-\e^2(t)
   \var^2+...\, .
 \end{aligned}
   \ee

Using \ef{gg1} and \ef{gg3} in
 balancing first the leading terms of the order
  $O(g'(t))$ yields the elliptic equation for $\psi$:
   \be
   \label{gg4}
    \tex{
 O(g'(t))\,\,\big(=O(\frac 1t)\big): \quad  \BB \psi - f'=0.
   }
   \ee
   Then balancing the rest of the terms in \ef{gg31} requires
   their asymptotic equivalence,
 \be
 \label{hh41}
  \tex{
 g''(t) \sim -(g'(t))^2 \sim \e(t), \,\,\, \mbox{i.e.,} \,\,\,
  g(t)= k \log t, \,\, \mbox{$g'(t)= \frac kt, \,\,
 g''(t)=- \frac k{t^2}$,} \,\, \e(t)= \frac 1{t^2}.
 }
 \ee
 Then, we obtain the second inhomogeneous singular
 Sturm--Liouville problem for $\var$:
  \be
  \label{SL1}
   \tex{
  O\big(\frac 1{t^2}\big): \quad
  \BB \var= k \psi + k^2(\psi'+\psi^2).
 }
   \ee

    Thus, the first simple asymptotic ODE in \ef{hh41} gives the $\log
    t$-dependence as in \ef{1.3}. Finally, we arrive at the
    following system for $\{\psi,\var\}$:
     \be
     \label{sys21}
     \left\{
     \begin{aligned}
& \BB \psi = f', \\
   &
   \BB \var= k \psi + k^2(\psi'+\psi^2).
    \end{aligned}
     \right.
 \ee
  Solving this system, with typical boundary conditions as in
  \ef{BC1}, allows then to continue the expansion of the solutions of
  \ef{gg2}
  close to an  ``affine ({\rm i.e., shifted via $f'$ on the RHS}) centre subspace" of $\BB$
  governed by the spectral pair obtained by translation in
  \ef{K2}:
   \be
   \label{pair1}
   \hat \l_0=0 \andA \hat \psi_0(y)= f'(y).
    \ee
 The asymptotic expansion for $t \gg 1 $ then takes the form
  \be
  \label{as22}
  \tex{
  w(y,t)= \frac kt \, \psi(y) + \frac 1{t^2} \,\var(y) + ...\, ,
  }
  \ee
  which can be easily extended by introducing further terms, with
  similar inhomogeneous Sturm--Liouville problems for the
  expansion coefficients.

  \ssk

{\bf Final ``spectral" remark: on $k$-values.} Thus, $\BB$ does
not have a discrete spectrum (see references at the beginning of
Section \ref{Sect.5.3}), so one cannot get a simple algebraic
equation for $k$ by demanding the standard orthogonality of the
right-hand side in the second equation in \ef{sys21} to the
adjoint eigenvector $\hat \psi^*_0$ of $\BB^*$ in the $L^2$-metric
(in which the adjoint operator $\BB^*$ is obtained), like
 \be
 \label{as23}
k: \quad   \langle k \psi + k^2(\psi'+\psi^2), \, \hat \psi_0^*
\rangle =0.
 \ee
 Therefore, the system \ef{sys21} cannot itself  determine
    the actual value of $k$ therein. As we have mentioned, the latter requires a
    difficult matching analysis of Inner and Outer Regions, which,
    for the KPP--4 (and all other problems), remains an open
    problem.

\subsection{On some other related results and references}

It is worth mentioning that, in similar cases, it is known that,
after a suitable time-rescaling and necessary transformations, the
orbit can approach the center subspace ${\rm Span}\{\hat \psi_0\}$
locally in space (in $L^2_{\rho, {\rm loc}}$). First results for
the semilinear heat equation
 $$
 u_t=\D u -u^p \quad \inA \quad (u>0)
  $$
   treated in
such a way were proved in \cite{GmV, KP, GKS, BPT}, etc. in the
middle of the 1980s; a full list of references is given in
\cite[Ch.~2]{SGKM}; see also various chapters in \cite{AMGV}. The
first realization of the scaling idea
 in quasilinear degenerate parabolic equations, which is one of the powerful tools to
  study such a stabilization 
 was established  by Kamin in 1973, \cite{K}.
  Sometimes, it is called now  a  renormalization group method.
 Basic ideas go back to the dimensional ideas in
 nonlinear problems and to the notion of self-similarity of the  second kind
  introduced by Ya.B.~Zel'dovich at the beginning of the 1950s. See Barenblatt's  book
  \cite{B},
   where several such ideas were discussed first. Of course, it
   should be noted that rigorous results obtained in the above
   papers are always related to semilinear and quasilinear {\em second-order}
parabolic equations and are heavily based on the Maximum
Principle. Therefore, any justification of such approaches to bi-
and poly-harmonic equations is expected to be very difficult.

    \ssk

    One can see that the above elementary conclusion well
    corresponds to an ``(affine) centre subspace analysis" of the
    non-autonomous PDE \ef{gg2}, and then $\t = \log t$ naturally
    becomes the corresponding ``slow" time variable; see various
    examples in \cite{GV1, GV2, AMGV} of such a slow motion along centre subspaces in
    nonlinear parabolic problems with global and blow-up
    solutions.
    In the latter case, the slow time variable is
     $$
     \t = - \ln
    (T-t) \to +\iy \quad \mbox{as $t$ approaches finite blow-up time $T^-$}.
    $$


 For the semilinear higher-order reaction-absorption equations
such as
 \be
 \label{utp}
 u_t=-u_{xxxx}-|u|^{p-1}u \inB \re \times \re_+ \withA p>1,
 \ee
existence of $\log t$-perturbed global asymptotics was established
in \cite{GalCr}. For finite-time extinction, with $-1<p<1$ in
\ef{utp}, this was done in \cite{Galp1}. For the corresponding
blow-up problem with the combustion source
 $$
 u_t=-u_{xxxx}+|u|^{p-1}u \inB \re \times \re_+ \withA p>1,
 $$
centre manifold-like $\log(T-t)$-dependent blow-up singularities
were constructed in \cite{Gal2m}. We must admit that a
justification of such $\log t$-corrections in any of KPP--$2m$
problems with $m \ge 2$ is  more difficult and has been obtained
for global orbits as $t \to +\iy$. For blow-up orbits, with
$|\log(T-t)|$, this remains an open problem.

\section{When the $\o$-limit set consists of  TW profiles}
 \label{S6omega}

Consider the $2m$th-order KPP problem for the PDE \ef{m1}, with
smooth bounded ``step-like"  initial data,
 \be
 \label{ub1}
 u(x,0)=u_0(x) \inB \re: \quad u_0(x) \to
 0, \,\,\, x \to +\iy \andA  u_0(x) \to 1, \,\,\, x \to -\iy,
  \ee
  where the convergence is assumed to be
sufficiently exponentially fast.

  \begin{theorem}
  \label{Th.Lam}
  Let \ef{m1}, \ef{ub1} admit a global uniformly bounded
  solution $u(x,t)$ having a ``step-like" form for all $t \gg 1$,
  i.e.,
 \be
 \label{ub2}
 u(x,t) \to 0 \asA x \to +\iy, \quad u(x,t) \to 1 \asA x \to -\iy
 \ee
 sufficiently exponentially fast uniformly in $t \gg1 $. Let, for
 all $t \gg 1$, the front location $x_f(t)$ can be uniquely and
 smoothly $($say, analytically$)$ defined from the equation
  \be
  \label{un3}
   \tex{
  u(x_f(t),t)= \frac 12 \forA t \gg 1,
 }
 \ee
with the following asymptotic representation:
  \be
  \label{un4}
  x_f(t) = \l_0 t - g(t) \whereA g'(t) \to 0 \asA t \to +\iy,
   \ee
   where $\l_0>0$ is a constant.
   Then
   $\l_0 \in \Lambda$, i.e.,   there is a  TW profile (maybe, non-unique) $f(\cdot;\l_0)$
    satisfying the ODE $\ef{m2}$, with $\ef{BC1}$,  and
  the omega-limit set $\o(u_0)$,
   defined via uniform convergence in the TW frame $y=x-\l_0
   t+g(t)$,
    is contained in this connected closed family of TW profiles:
     \be
     \label{un5}
     \o(u_0) \subseteq\{f(\cdot;\l_0)\}.
      \ee
 \end{theorem}

 \noi{\em Proof.} (i). Given a sequence $\{t_k\} \to +\iy$, passing to the
 limit(in a weak sense, implying, by parabolic regularity,
 stronger convergence)
in the rescaled perturbed equation \ef{K21} as $\{t_k+s\} \to
+\iy$, one concludes that
 \be
 \label{un6}
 v(\cdot,t_k+s) \to w(\cdot,s) \quad \mbox{in, say,} \quad L^\iy_{\rm
 loc}(\re; \re_+) \quad (\mbox{or in $C_{\rm loc}(...)$, etc.}),
  \ee
 where, by the assumption in \ef{un4}, $w(y,s)$ solves (in the classic sense) the
 corresponding unperturbed equation
  \be
  \label{un61}
  w_s=\AAA w +\l_0 w_y + w(1-w) \inB \re \times \re_+,
   \ee
   with data $w_0 \in \o(u_0)$.
   Moreover, by the definition of the front tracking \ef{un3},
   \be
   \label{un7}
    \tex{
   w(0,s) \equiv \frac 12 \LongA D_s^k w(0,s) =0 \quad \mbox{for
   any} \quad k=1,2,3,...\, .
   }
   \ee
 Since the solutions of \ef{un6} are analytic in both $y$ and $s$ for $s>0$ \cite{Fr58}
 (see also related more recent results in
  \cite{Lun}),
 \ef{un7} implies that $w(y,s)$ is independent of $s$, so, under
 the above hypothesis, $w(y,s) \equiv f(y;\l_0)$.
 \quad $\qed$

\ssk

\noi{\bf Remark: Blow-up in the parabolic problem is possible.}
The conditions on data $u_0(x)$ at the beginning of Theorem
\ref{Th.Lam} are essential. Indeed, if $-u_0^2(x) \ll u_0(x) \le
-1$ on some interval, then any such solution blows-up in finite
time as $t \to T^- \in \re_+$:
 \be
 \label{BL1}
 \tex{
 \lim_{t \to T^-} \,\, \sup_{x \in \re} |u(x,t)|=+ \iy.
 }
 \ee
 Different proofs can be found in \cite{EGKP1, GPInd}.
 Then blow-up can be self-similar \cite{BGW1} and then blow-up {\em is
 not uniform}, i.e.,
  \be
  \label{BL2}
\tex{
 \lim_{t \to T^-} \,\, \inf_{x \in \re} u(x,t)=- \iy \andA
 }
 \tex{
 \lim_{t \to T^-} \,\, \sup_{x \in \re} u(x,t)=+ \iy.
 }
 \ee
As an alternative, blow-up may create centre and stable manifold
patterns, \cite{Gal2m}.

\ssk

\noi{\bf Remark 1: Is the $\o$-limit set a point? -- No answer.}
For various smooth (or analytic) gradient systems, under some
extra assumptions, it is known that the omega limit sets
 consists of a single point. We refer to, e.g., as a first such result,
to Hale--Raugel's approach that applies to guarantee convergence
of the orbits of gradient systems; see \cite{Hale92}, where a
survey and further references are given. This approach essentially
relies on spectral properties of the linearized operator (main
Hale--Raugel's hypothesis is that the equilibrium has multiplicity
at most one, or $k>1$ under special hypothesis) and uses
properties of stable, unstable, and center manifolds, that are
difficult to justify for some less smooth equations. We also refer
to  \cite{Busca02}, where a similar approach to stabilization is
used and other references can be found.

An alternative application is  the \L{ojasiewicz}--Simon approach
to parabolic equations. Namely, the \L{ojasiewicz}--Simon
inequality is an effective tool of studying of the stabilization
phenomena in various evolution problems. In particular, it
completely settles the case of analytic nonlinearities; see
references in \cite{Chill}; see also \cite{Fer, Haraux, Har2011,
Ver2011, Yas2011} and references therein.
Another approach \cite{GPSGrad} is based
on Zelenyak's ideas (1968) \cite{Zel} from  one dimension that are
mainly connected with Lyapunov functions only for gradient
dynamical systems; see also a short survey and references in
\cite{GPSGrad}.

Unfortunately, \ef{un61} is not a gradient system in the usual and
necessary sense (and \ef{K21} is a not a ``perturbed" one, for
which there exists a Lyapunov-like functional that is ``almost"
monotone on evolution orbits). In fact, the bi-harmonic equation
\ef{E4} (as well as any of \ef{m1}) admits a ``pseudo-Lyapunov
function", obtained, as usual, by multiplication of the equation
 by $u_t$ in $L^2$:
  \be
 \label{Z1}
  \tex{
  \frac{\mathrm d}{{\mathrm d}t}\, L[u](t) \equiv
   \frac{\mathrm d}{{\mathrm d}t} \, \big[\frac 12 \int(u_{xx})^2+
   \int\big(c(x)- \frac 12 \, u^2 + \frac 13\, u^3\big)=\int
   (u_t)^2 \le 0,
   }
   \ee
 where $c(x)$ is any smooth function with an exponential convergence  as
 $x \to \iy$ satisfying, for the sake of convergence of the second
 integral in \ef{Z1},
  \be
  \label{Z2}
   \tex{
  c(-\iy)= \frac 16 \andA c(+\iy)=0.
  }
  \ee
Unfortunately, $L[u](t)$ is not bounded below, so cannot be used
in the Lyapunov classical analysis. This is also easily seen
 in the TW frame, where, for $v=v(y,s)$, \ef{Z1} reads
  \be
  \label{Z3}
  \tex{
\frac{\mathrm d}{{\mathrm d}t} \, \big[\frac 12 \int(v_{xx})^2+
   \int\big(c(y+\l_0 t-g(t))- \frac 12 \, u^2 + \frac 13\, u^3\big)=\int
   (v_t-\l_0 v_y-g'(t)v_y)^2 \le 0.
   }
   \ee
Observe that, here, $c(y+\l_0 t-g(t)) \to 0$ as $t \to +\iy$
uniformly on compact subsets in $y$.
 Since such a Lyapunov function has almost nothing to do with
TW profiles (i.e., equilibria), we do not think that these
inequalities can be of any help in the stabilization problem.

On the other hand, for generalized gradient KPP-problems in
Appendix A, where the ODEs for $f$ are variational, these
approaches can apply; see comments therein.



\ssk

\ssk

\noindent{\bf Acknowledgement.} The author would like to thank
D.~Williams, as well as other active participants of seminars of
Reaction-Diffusion and Probability Groups,  Department of
Mathematical Sciences, University of Bath, who stimulated a number
of discussions of PDE aspects of the KPP--2 problem in 1995, which
led the author to initiate such a research, \cite{Gal1995}.


\section*{\sc Appendix A. Some modified variational KPP-problems}
 \label{SVar}

\subsection{A more general semilinear bi-harmonic equation}

Consider the following bi-harmonic equation, with parameters
$\a,\b,\g \in \re$:
 \be
 \label{A1}
  \begin{aligned}
  \mbox{(PDE)}: \quad &
 u_t=-u_{xxxx} + \a u_{xxx} + \b u_{xx}+\g u_x+u(1-u), \\
 \mbox{(ODE)}: \quad &
  \hat \AAA f \equiv -f^{(4)}+\a f'''+\b f''+ (\l+\g) f'=-f(1-f).
 \end{aligned}
  \ee

\begin{proposition}
\label{Pr.Var1}
 For a fixed $\l \in\re$,
the ODE in \ef{A1} is $($formally$)$ variational in the metric of
the weighted space $L^2_\rho(\re)$, where
 \be
 \label{A2}
  \tex{
 \rho=\eee^{- \frac {\a y}2} \quad \mbox{and $\a \in \re$ solves
 the cubic equation} \quad \frac \a 2 \big( \frac {\a^2}4+\b
 \big)=-(\l+\g).
 }
  \ee
  \end{proposition}

\noi{\em Proof.} Using the general form of symmetric
(self-adjoint) ordinary differential operators \cite{Nai1}, we
have to find two coefficients $p_1(y)$, $p_2(y)$, and the weight
$\rho(y)$ such that
 \be
 \label{A3}
  \tex{
  \hat \AAA f \equiv -\frac 1\rho\, [(p_1 f'')''+(p_2f')'] \equiv
 - \frac 1 \rho \,[p_1f^{(4)}+2 p_1' f'''+(p_1''+p_2)f''+p_2'f'].
  }
  \ee
This yields the following system of linear differential equations:
 \be
 \label{A4}
  \left\{
  \begin{aligned}
  & \tex{\frac{p_1}\rho=1},\\
 & \tex{\frac {2p_1'}\rho=-\a},\\
 & \tex{\frac{p_1''+p_2}\rho=-\b},\\
 & \tex{ \frac{p_2'}\rho=-(\l+\g)}.
 \end{aligned}
 \right.
 \ee
Hence, $p_1=\rho$, so, the second equation yields the positive
weight in \ef{A2}:
 \be
 \label{A5}
  \tex{
   \frac {2 \rho'}\rho=-\a \LongA \rho(y)=\eee^{-\frac{\a y}2}>0.
   }
   \ee
The third equation then defines
 \be
 \label{A6}
  \tex{
 p_2=-p_1''-\b \rho=-\rho''- \b \rho=-\big(\frac{\a^2}4+\b\big)\,\eee^{-\frac{\a
 y}2}.
 }
 \ee
Finally, the last equation in \ef{A4} yields the necessary
condition on $\a$:
 \be
 \label{A7}
  \tex{
  p_2'=-(\l+\g)\rho \LongA p_2(y)= \frac 2 \a\,(\l+\g) \, \eee^{-\frac{\a
  y}2},
   }
   \ee
   i.e.,  comparing  with $p_2$ in \ef{A6} gives the cubic
   algebraic equation in \ef{A2}. \quad $\qed$

   \ssk

 Thus, expressing all the coefficients in \ef{A3} in terms of
 the weight $\rho$ given in \ef{A2} yields the following problem
 with a symmetric operator:
  \be
  \label{A8}
   \tex{
  \hat \AAA f \equiv -\frac 1\rho\,\big[ (\rho f'')''+
   \frac{\l+\g}\a\,(\rho f')'\big]=-f(1-f),
    }
    \ee
    admitting the (formal) potential:
     \be
     \label{A9}
      \tex{
      \hat \Phi(f)=- \frac 12\, \int \rho(f'')^2+ \frac{\l+\g}{2 \a}
 \int \rho(f')^2 + \int \rho \big( \frac 12\, f^2- \frac 13
 f^3\big) \inB H^2_\rho \cap L^3_\rho.
  }
  \ee

  On the other hand, for the corresponding rescaled limit ($t=+\iy$) parabolic
  equation in the TW-frame \ef{un61}, where
   \be
   \label{A12N}
    \tex{
   \hat \AAA = \AAA +
  \l_0 I \LongA \hat \Phi[f]=
 \Phi[f]+ \frac {\l_0}{2 \a}\, \int \rho(f')^2,
 }
 \ee
   the potential structure of the operator means existence
  of a good Lyapunov function:
   \be
   \label{A11N}
    \tex{
     \frac{\mathrm d}{{\mathrm d}s}\, \Phi[w](s) = - \int \rho
     (w_s)^2 \le 0.
      }
      \ee
As was mentioned in Section \ref{S6omega}, this allows a further
identification of the $\o(u_0)$.

For the former perturbed equation in \ef{K21}, similar
computations yield
 \be
   \label{A11N1}
    \tex{
     \frac{\mathrm d}{{\mathrm d}t}\, \Phi[v](t) = - \int \rho
     (v_t)^2 +\l_0 g'(t)\, \int \rho v_y v_t.
      }
      \ee
Therefore, a proper control of the second term on the right-hand
side is necessary. It perfectly converges and is $\sim O(\frac
1{t^2})$ as $t \to +\iy$ (as the first term), if we trust the
formal expansion \ef{as22}. In this case, $\Phi[v](t)$ is
``almost" Lyapunov function that allows to pass to the limit and
to identify $\o(u_0)$ as a point. However, a full rigorous proof
of \ef{as22} remains an open problem.

  \subsection{Application to a particular potential model}

  We set $\b=\g=0$ in \ef{A1},
   so that, for a fixed $\l>0$, we have from
  \ef{A2}
   \be
   \label{A10}
   \a_0(\l)=-2 \l^{\frac 13}<0 \LongA
    \rho=\eee^{\l^{1/3} y} \to 0 \asA y
   \to -\iy.
   \ee
The model now reads
 \be
 \label{A1N}
 u_t=-u_{xxxx}+ \a_0(\l) u_{xxx}+u(1-u) \LongA - \l f'=-f^{(4)}+
 \a_0(\l)
 f'''+f(1-f).
  \ee

Then, \ef{A10} means that the metric is now suitable for the ODE
in \ef{A1N} with the second condition in \ef{BC1} saying that
$f(y) \to 1$ as $y \to - \iy$.

Concerning the first condition, $f(y) \to 0$ as $y \to +\iy$, the
suitability of the $L^2_\rho$-metric (the condition $f \in
L^2_\rho(\re_+)$) will depend on the fact whether the necessary
asymptotic bundle (to be matched with that about the equilibrium 1
as $y \to -\iy$, as in Section \ref{S2.5}) governed by the
linearized operator,
 \be
 \label{A11}
 g^{(4)}-\a_0 g'''-\l g'-g=0, \quad g(y)=\eee^{\mu y} \LongA H_+(\mu,\l) \equiv\mu^4+2
 \l^{\frac 13} \mu^3-\l \mu-1=0,
  \ee
  belongs to $L^2_\rho(\re_+)$.
The characteristic equation in \ef{A11} is reduced to a
bi-quadratic one:
 \be
 \label{A12}
  \tex{
 \mu=- \frac 12\, \l^{\frac 13}+\nu \LongA \nu^4- \frac 32\,
 \l^{\frac 23} \nu^2+ \frac 5{16}\, \l^{\frac 13}-1=0
 }
 \ee
 \be
 \label{A13}
  \tex{
 \LongA \nu^2_\pm = \frac 34\, \l^{\frac 23} \pm \sqrt{ \frac 14\,
 \l^{\frac 43}+1}.
 }
  \ee
  Then $\nu_+^2>0$, while $\nu_-^2$ yields
 the following ``critical" value of the TW speed:
  \be
  \label{A14}
   \tex{
   \nu_-=0 \LongA \l_*= \big( \frac {16}5 \big)^{\frac
   34}=2.3926...\, .
   }
   \ee
   Indeed, at $\l=\l_*$, the characteristic equation in \ef{A11}
   has a double root $\mu_0=- \frac 12\, \l_*^{\frac 13}$, and,
    \be
    \label{A14N}
    \eee^{\mu_0 y}, \,\,\, y \,\eee^{\mu_0 y} \not \in
    L^2_\rho(\re_+),
     \ee
so that, by this functional variational setting, a 2D bundle of
exponentially decaying orbits as $y \to +\iy$ are excluded.
Clearly, this diminishes any possibility to have a proper critical
point of the corresponding functional \ef{A9}, as a TW profile.
Indeed, the explicit eigenvalues \ef{A13} allow one to perform
further asymptotic analysis, which, together with a similar one as
$y \to -\iy$, eventually may describe a precise range of $\l$, for
which  critical points of $\Phi$ correspond to TWs. But we will
not do that and explain our reasons later on.

Moreover, since the weight \ef{A10} is decaying as $y \to -\iy$, a
possible critical point of $\Phi(f)$ may have a behaviour therein
that has nothing to do with a TW $f(y)$. In other words, in this
case, a full and a complete classification of all the asymptotics
as $y \to -\iy$ of solutions of the ODE in \ef{A1N} is crucially
necessary.

Overall, though potential operators allow to perform a more
rigorous analysis of critical points, but an extra essential work
is inevitable to guarantee the necessary result concerning TW
profiles.

 \ssk

However, our numerical experiments show a different and larger
than \ef{A14} ``critical" value of $\l$. Namely, integrating
 the ODE in \ef{A1N} on the interval $[-100,50]$ yields a kind of
 a ``minimal" speed of existence
  \be
  \label{A15}
   3.38<\l_{\rm min} \le 3.39,
 \ee
 such that, for $\l < \l_{\rm min}$,  there are no TW profiles
 $f(y)$. However, the convergence is very slow and, actually, a
 full convergence with Tols$\sim 10^{-4}$ never happened, even
 using 16000 points or more. Therefore, we do not know how much
 the value \ef{A15} may depend on boundary conditions and the
 left-hand end point $L_-=-100$.

\begin{figure}
\centering
\includegraphics[scale=0.8]{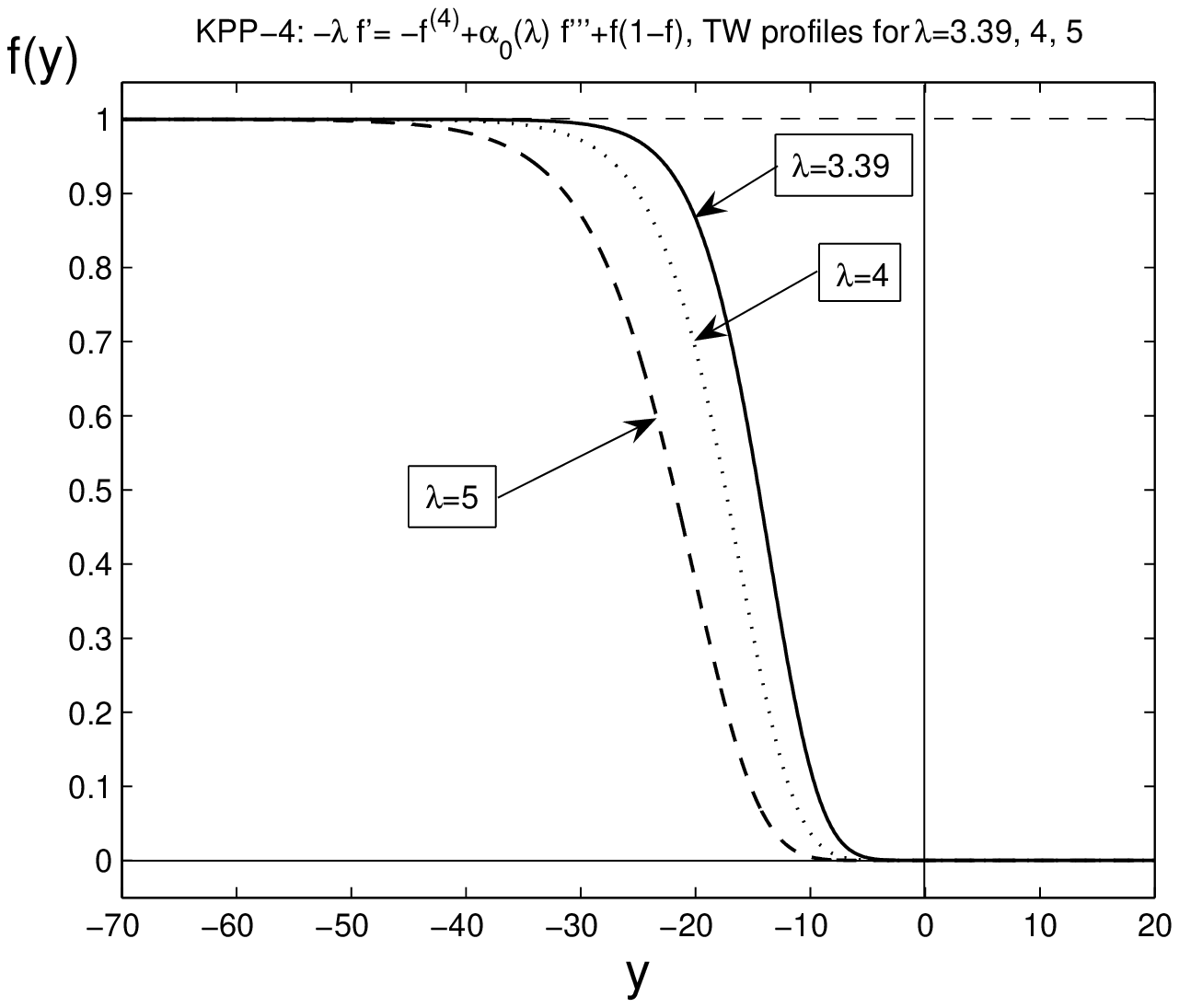}  
\vskip -.3cm
  \caption{The TW profiles $f(y)$ of
(\ref{A1N}), \ef{BC1} for
 $\l=3.39$, 4, and 5.}
 \label{Fuxxx1}
\end{figure}

In Figure \ref{Fuxxx1}, we show TW profiles satisfying \ef{A1N}
for $\a=3.39$, 4, and 5. Observe a non-oscillatory character of
solutions in a neighbourhood of the constant equilibria that well
corresponds to the distribution of the eigenvalues in \ef{A13} for
$\l>\l_*$.


\subsection{Non-variational KPP-problem: arbitrary $\alpha$}

Indeed, with a certain extra analysis (and certain doubts), the
above variational character of the problem \ef{A1N}, with
$\a=\a_0(\l)$, $\l>0$, improves its overall mathematical quality.
However, we paid so little attential to the actual possible
rigorous study of \ef{A1N}, because the potential functional
setting here has nothing to do with the existence of TW profiles.

To show this, consider the same problem but now with an {\em
arbitrary} $\a \in \re$:
 \be
 \label{ar1}
 u_t=-u_{xxxx}+\a u_{xxx}+u(1-u) \LongA -\l f'= -f^{(4)}+ \a
 f'''+f(1-f),
  \ee
  which {\em is not variational} for any $\a \not = \a_0(\l)=-2
  \l^{\frac 13}$. Note that, for $\a>0$, the positivity result
  \ef{Lam1} remains valid and the identity \ef{Lam2} now reads
   \be
   \label{Lam4}
    \tex{
    -\l \int(f')^2=-\a \int (f'')^2 -\frac 16<0.
    }
    \ee

\begin{figure}
\centering
\includegraphics[scale=0.8]{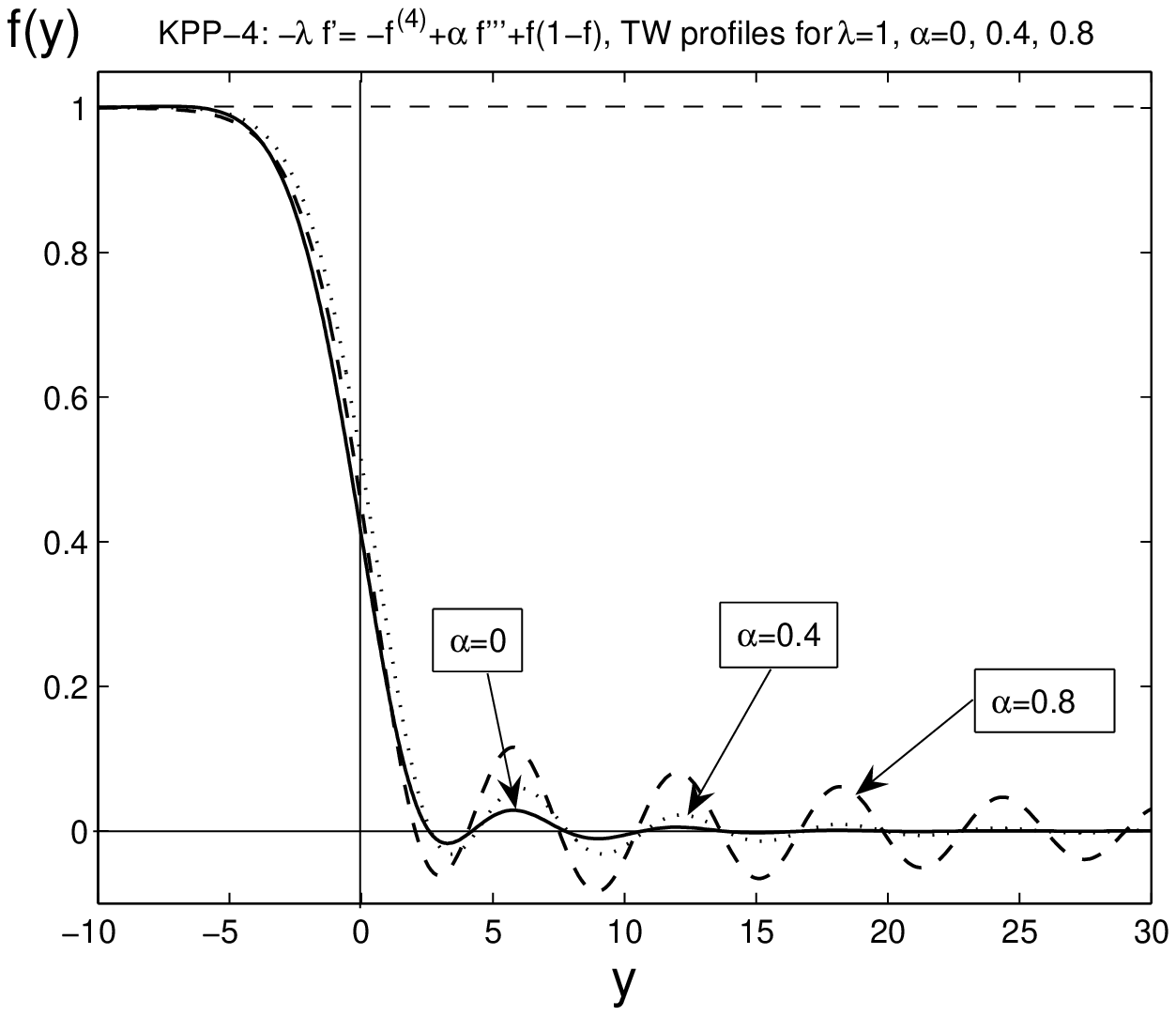}  
\vskip -.3cm
  \caption{The TW profiles $f(y)$ of
(\ref{ar1}), \ef{BC1} for $\l=1$ and
 $\a=0, \, 0.4, \,0.8$.}
 \label{Fuxxx2}
\end{figure}

\begin{figure}
\centering
\includegraphics[scale=0.8]{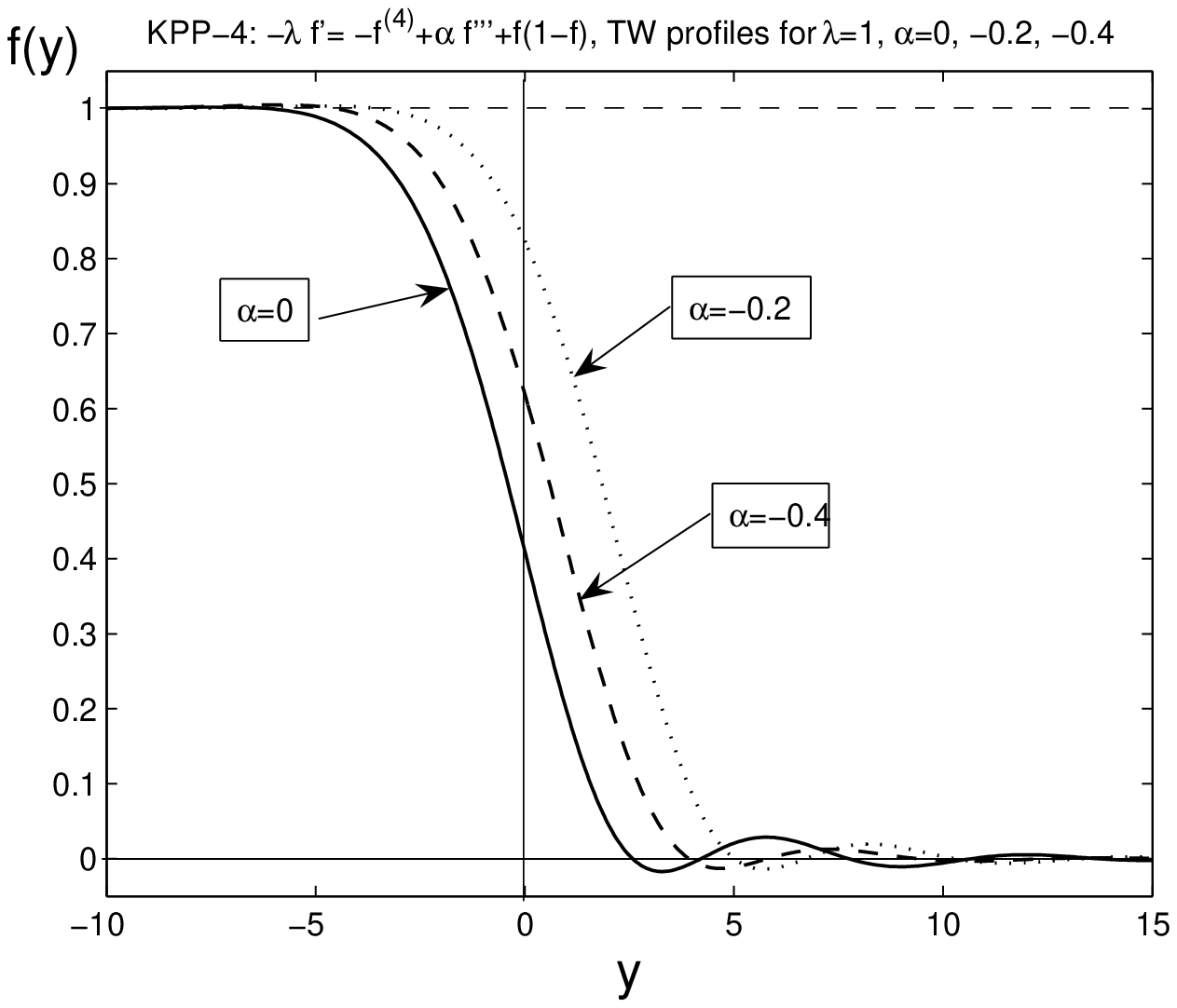}  
\vskip -.3cm
  \caption{The TW profiles $f(y)$ of
(\ref{ar1}), \ef{BC1} for $\l=1$ and
 $\a=0, \, -0.2, \,-0.4$.}
 \label{Fuxxx3}
\end{figure}

  It turns out that, numerically, with the full convergence up to
  $10^{-5}-10^{-7}$
  (unlike the case $\a=\a_0(\l)$ above), we observe existence of
  TW profiles for both values $\a>0$ (Figure \ref{Fuxxx2}, $\l=1$) and
  $\a<0$ (Figure \ref{Fuxxx3}, $\l=1$).
 In other words, there exists a local curve of solutions
 originated at $\a=0$ with existence in some neighbourhood $\a
 \approx 0$.

   For this $\l=1$, we have
  $\a_0(1)=-2$, but, for $\a<0$, we were not able to extend the
  solution for such negative $\a$'s. Namely, we
found a ``{\em minimal}" value
 \be
 \label{almin1}
 -0.413 < \a_{\rm min}(1) \le -0.412
 \quad(\mbox{integrating on $[-50,50]$}),
  \ee
such that, for $\a<\a_{\rm min}(1)$, there is no convergence at
all, i.e., $f$ seems to be nonexistent.

Therefore, here and in \cite{GKPPII, GKPPIII}, we always did not
follow any ``variational temptation" to make our analysis  more
rigorous for some particular and specially designed KPP-problems
with potential operators, and, instead, concentrated on more
general aspects of the KPP-ideology and necessary and required
related mathematics.


  \end{small}
   \end{appendix}

\begin{appendix}
\begin{small}
\setcounter{section}{2} \setcounter{equation}{0}

\section*{\sc Appendix B. A linear model: level propagation estimates from above by comparison with majorizing
order-preserving integral evolution}
 \label{SMaj}

 For the semilinear KPP--$2m$ problem \ef{m1}, \ef{1.H}, there is
 no the usual
 comparison of solutions for any $m \ge 2$. We
 consider the corresponding linear {\em poly-harmonic equation}
 \be
  \label{PPP1}
   u_t=(-1)^{m+1} D_x^{2m} u  \inB \re \times \re_+, \quad
   u_0(x)=H(-x),
    \ee
  which admits certain comparison techniques
  \cite{GPInd} by  constructing
 a so-called {\em majorizing order-preserving operator}.
This is done as follows. We first present the fundamental solution
of the poly-harmonic \ef{PPP1}:
 \be
 \label{PP2}
  \tex{
 b(x,t) = t^{-\frac 1{2m}}\, \eee^t \,F(y), \quad y=\frac x{t^{1/2m}}.
 }
\ee
  Substituting $b(x,t)$ into (\ref{PPP1}), one obtains the
symmetric (even) profile
 $F$ as a unique solution of the linear  ODE
\begin{equation}
\label{ODEf}
 \tex{
 {\mathbf B}F \equiv (-1)^{m+1} F^{(2m)} + \frac 1{2m} \, F'
y + \frac 1{2m} \,F = 0 \quad  {\rm in} \,\,\, \re, \quad
\int_{{\mathbb R}} F(y) \, {\mathrm d}y = 1.
 }
\end{equation}
 The rescaled kernel $F(y)$ is oscillatory as $y \to \iy$ and
satisfies a standard pointwise estimate \cite{EidSys}, which is
convenient to present in the  form
\begin{equation}
\label{kern}
 \tex{
 |F(y)| \le  D_0 \bar F(y), \,\, \bar F(y) = \omega_1
\eee^{-d_0|y|^{\alpha}}>0 \quad {\rm in} \,\,\, \re,
\,\,\mbox{where}\,\, \o_1= \left(\int_{\re}
\eee^{-d_0|y|^{\alpha}} \, {\mathrm d}y \right)^{-1},
 }
\end{equation}
 and $D_0$, $d_0$ are some positive
constants depending on $m$. The normalizing constant $\o_1$ in
\ef{kern}  guarantees that the positive kernel $\bar F$ satisfies
the standard normalization:
 \be
 \label{barF}
  \tex{
 \int_{\re} \bar F(y) \, {\mathrm d}y = 1,
 }
  \ee
 that will allow us to use it to create below a majorizing
 order-preserving operator. The exponent $\a$ is given by
 $$
  \tex{
 \a =  \frac {2m}{2m-1} \in (1,2) \quad \mbox{for all} \,\,\, m
 \ge 2 \quad (\mbox{$\a=2$ for $m=1$, i.e., for the Gaussian
 kernel}).
 }
 $$

 Here $D_0>1$ (or $D_0-1>0$) is called the {\em order deficiency} of
 the poly-harmonic
semigroup. One can see that $D$ can be made to be equal to 1,
$D=1$, iff $F \ge 0$, i.e., the semigroup is {\em
order-preserving}, and this happens for $m=1$ only.

Next, it is key for us to get a sharp representation of the
constant $d_0=d_0(m)$. By WKBJ-type two-scale asymptotics of the
solution $F(y)$ of \ef{ODEf}, it is proved that $d_0$ is
associated with the root of the algebraic equation
 \be
  \label{alg11}
   \tex{
(-1)^m a^{2m-1} =-  \frac 1{2m}.
 }
 \ee
 It has $2m-1$ different
roots $\{a_k\}, \,\,$ and in (\ref{kern}), we  choose
 \be
 \label{d111}
  \tex{
 - d_0 =
\max_k \, \{{\rm Re} \, \a_k <0\}.
 }
 \ee
 For the bi-harmonic flow \ef{PPP1} with $m=2$, the
 rescaled kernel estimate \ef{kern} is
\begin{equation}
  \label{F11}
 |F(y)| \le D \,\bar F(y) \equiv D \, \o_1\,{\mathrm e}^{-d_0 |y|^{4/3}} \quad \mbox{in} \quad
 \re \whereA d_0=3 \cdot 2^{-\frac{11}3} =0.2362... \,.
   \end{equation}

 In \ef{kern}, there appears $\bar
F$ as the rescaled kernel of the majorizing order-preserving
integral evolution equation,
 \be
 \label{maj66}
 \tex{
 \bar b(x,t)= t^{\frac 1{2m}} \, \bar F(y)>0, \quad y= \frac x{t^{1/2m}},
  }
  \ee
 whose solutions can be compared with
those of (\ref{PPP1}), \cite{GPInd}.
 Obviously, for the linear KPP equations (i.e., \ef{m1} without
 the quadratic term $-u^2$) and the corresponding majorizing
 equation, the fundamental solutions are
  \be
  \label{maj67}
b_*(x,t)= \eee^t \, b(x,t) \andA \bar b_*(x,t)= \eee^t \, \bar
b(x,t). \ee

 Consider this majorizing
equation and the corresponding poly-harmonic flow
 \be
 \label{maj1}
  \tex{
  \bar u(t)= \bar b_*(t)*\bar u_0, \,\,\, \bar u_0(x) \ge 0  \andA u(t)=  b_*(t)* u_0 \forA t \ge
  0,
  }
  \ee
 with some initial data $u_0$ and $\bar u_0$ from the
 corresponding weighted $L^2_\rho(\re)$ ($\rho = \eee^{- a_0|x|^\a}$, $a_0 \in
 (0,2d_0]$)
  space of correctedness.
 Subtracting these two and using \ef{kern} yield
  \be
  \label{maj2}
   \aligned
  \bar u(t)-u(t) =\, & \bar b_*(t)* \bar u_0 - b_*(t)*u_0 \ge \bar b_*(t)* \bar u_0 -
  |b_*(t)|*|u_0|
  \ssk\\
   \ge \, & \bar b_*(t)* \bar u_0 -  D_0 \bar b_*(t)*|u_0| \ge
  \eee^t \, \bar b(t)*(\bar u_0- D_0 |u_0|),
 \endaligned
  \ee
  whence the following result:

  \begin{proposition}
  \label{Pr.Comp2}
  $(${\bf Comparison}$)$
There is comparison of solutions of $\ef{maj1}$, i.e.,
 \be
 \label{maj12}
 u(x,t) \le \bar u(x,t) \inB \re \times \re_+ \,\,\, \mbox{for majorizing data}
 \,\,\,\bar u_0(x) = D_0 |u_0(x)| \ge 0 \,\, \mbox{in} \,\, \re.
 \ee

 \end{proposition}

Consider now the majorizing solution
 \be
 \label{maj99}
 \bar u(x,t)= m_0 \eee^t \, \bar b(x,t) \,\,\, \mbox{corresponding
 to Dirac's delta function} \,\,\, \bar u_0(x)=m_0 \d(x),
 \ee
 where $m_0$ is a sufficiently large constant to be estimated.
 We now treat \ef{maj99} as a solution of the boundary value problem in $\re_+
 \times \re_+$ with the symmetry conditions at the boundary point
  $x=0$ (of course, at the same time, this is a solution of the Cauchy problem):
   \be
   \label{maj90}
   u_x(0,t)=u_{xxx}(0,t)=...=D_x^{2m-1}u(0,t)=0 \forA t \ge 0.
   \ee
 We now choose $m_0>0$ so large that
  \be
  \label{maj91}
 D_0  u(0,t) \le     \bar u(0,t) \forA t \ge 0.
   \ee
 Then the same argument of majorizing comparison guarantees
that
 \be
 \label{maj92}
 u(x,t) \le \bar u(x,t) \inB \re_+ \times \re_+.
  \ee

This inequality allows to estimate the $\frac 12$-level
propagation \ef{TW2} for $u(x,t)$: indeed, $x_f(t) \le \bar
x_f(t)$, where $\bar x_f(t)$ solves
 \be
 \label{maj93}
  \tex{
  \eee^t \, t^{-\frac 1{2m}} \bar F\big( \frac {\bar x_f(t)}{t^{1/2m}}\big)=
  \frac 12 \orA
 \bar F\big( \frac {\bar x_f(t)}{t^{1/2m}}\big)= t^{\frac 1{2m}} \,
 \eee^{-t} \to 0 \asA t \to \iy.
 }
 \ee
 Using \ef{kern}, we finally obtain the estimate
  \be
  \label{maj94}
   \tex{
   x_f(t) \le \bar x_f(t)= d_0^{- \frac{2m-1}{2m}} t- \frac
   {2m-1}{(2m)^2}\,d_0^{- \frac{2m-1}{2m}} \log t(1+o(1)), \quad t \to
   \iy.
}
 \ee

 Therefore, formally (without a proof, since comparison was not
 performed for full nonlinear KPP--problems; see also comments below), we have the following
 bound for the minimal speed:
  \be
  \label{maj95}
   \l_0 \le \l_*= d_0^{- \frac{2m-1}{2m}} \andA
    \l_0 \le \l_*=3^{-\frac 34} \cdot 2^{\frac{11}4}=2.9512... \forA m=2.
     \ee
 If this minimal speed actually takes place\footnote{As we know, by \ef{t4}, this cannot happen
  for $m=2$.}, then \ef{maj94} gives
 a lower estimate of $k=k(m)$:
  \be
  \label{maj96}
   \tex{
   \l_0=\l_*: \,\,\,
  k \ge k_*=
  \frac
   {2m-1}{(2m)^2}\,d_0^{- \frac{2m-1}{2m}}, \,\,\, k \ge k_*=
 3^{\frac 14} \cdot 2^{-\frac{5}4}=0.5533... \forA m=2.
  }
  \ee

Finally, it is worth mentioning  that the estimates \ef{maj95} and
\ef{maj96} make more sense for the KPP-like problems for the
semilinear majorizing order-preserving integral,
``pseudo-differential" evolution equation\footnote{This kind of
equations we  meant  in Section \ref{SDiscr}, thought, it seems
that, for the majorizing operator in \ef{maj999}, the
corresponding  rescaled operator $\BB$ in \ef{gg2} still does not
have a discrete spectrum.} (recall again that it is not a flow,
with no a  semigroup and $t$-translational invariance available)
 \be
 \label{maj999}
  \tex{
\bar u(t)=\bar b(t)* \bar u_0 + \int\limits_0^t \bar
b(t-s)*u(s)(1-u(s)) \, {\mathrm d}s,
 }
\ee which is far away from semigroups for the original
KPP--problems for higher-order differential operators.


\ssk

 \noi{\bf Remark: on comparison in semilinear KPP--$2m$
problem.} For $m=1$, there is no problem, since $-u^2 \le 0$, so
that
 $\bar u(x,t)$ is always a super-solution.

 Let $m \ge 2$, then the integral equation via the Duhamel
 principle reads
  \be
  \label{maj97}
   \tex{
  u(t)= \eee^t b(t) * u_0 - \int\limits_0^t b(t-s)*(\eee^s u^2(s))\,
  {\mathrm d s}.
  }
 \ee
 Therefore, the final inequality similar to \ef{maj2} now
 contains an extra term:
 \be
 \label{maj98}
  \tex{
 \bar u(t)-u(t)  \ge
  \eee^t \, \bar b(t)*(\bar u_0- D_0 |u_0|) + \int\limits_0^t b(t-s)*(\eee^s u^2(s))\,
  {\mathrm d s} \equiv J_1(t) + J_2(t),
  }
    \ee
    where $J_1(t) \ge 0$ by Proposition \ref{Pr.Comp2}.
   However, since $b(x,t)$ changes sign, in general, we cannot control
   the positive sign of the last term (though $u^2 \ge 0$) for
   absolutely arbitrary solutions. Thus, any ``unconventional"
   comparison results are not available, and this is well known
   for higher-order parabolic flows.

   But we recall that we need a proper estimate above for $t \gg
   1$, where our solution is assumed to be approaching (with a
   shift) to a TW, i.e., we may assume that
    \be
    \label{maj40}
    u(x,t) \to f(x- x_f(t)), \,\,\, x_f(t)= \l_0 t - g(t), \quad \mbox{uniformly in $x \in \re$ as} \quad
    t \to +\iy.
     \ee
 Moreover, we need such an upper estimate on a {\em level set}
 $\{u(x,t)= \frac 12\}$ for $t \gg 1$. Therefore, the
 second term on the right-hand side of \ef{maj98} can be
 arbitrarily sharply  approximated by
  \be
  \label{maj54}
   \tex{
  J_2(t) \sim  \int\limits_0^t b(t-s)*(\eee^s f^2(\cdot))\,
  {\mathrm d s} \asA t \to +\iy,
  }
   \ee
   and this approximation gets sharper as $t $ increases.
   Bearing in mind rather smooth TW profiles $f(y)$ (cf. Figures
   \ref{F0}, \ref{F1}, \ref{F4}, etc.),  representing an ``almost
   monotone" (with just slightly oscillating tails at $y \to \pm \iy$) heteroclinic
path of equilibria 1 and 0, and overall ``dominant positivity" of
the fundamental solution $b(x,t)$ ($\int b(x,t)\, {\mathrm d}x
\equiv +1$), there is a sufficient probability that
 \be
 \label{maj55}
  J_2(x,t) > 0 \forA t \gg 1 \quad \mbox{and for large enough
  $x$,}
   \ee
at least in an ``a.a. sense", in a natural way. If \ef{maj55} is
true, then, by \ef{maj98}, the front propagation estimates
\ef{maj94}--\ef{maj96} actually hold.

The positivity \ef{maj55} of the integral in \ef{maj54} deserves
further analytic study, though it is very difficult, since the TW
profiles are not known explicitly. On the other hand, let us note
that the asymptotic positivity \ef{maj55} on the front curve
$\{x=x_f(t)\}$ can be studied numerically, but, without any doubt,
we believe that it actually takes place justifying the above TW
speed and $k$-estimates.


\end{small}
\end{appendix}

\end{document}